\documentclass[12pt]{article}
\usepackage[a4paper]{geometry}
\usepackage{ajc-new-CC-BY-ND}
\usepackage{tikz}
\usepackage{array}
\usepackage{multirow}
\usepackage{placeins}
\newcommand{\Asterisk}{\mathop{\scalebox{1.5}{\raisebox{0.0ex}{{\footnotesize{$\ast$}}}}}}%
\newcommand\unv[1]{
\phantom{\Asterisk} \hfill #1 \hfill \Asterisk
}

\tikzstyle{smallgraph}=[
  every node/.style={circle,draw,fill=black!30,inner sep=0pt,minimum width=4pt},
  every path/.append style={semithick}
]
\tikzstyle{largegraph}=[
  every node/.style={circle,draw,fill=black,inner sep=0pt,minimum width=6pt},
  every path/.append style={semithick}
]

\pgfmathsetmacro{\threevradius}{0.375*sqrt(2/3)}
\pgfmathsetmacro{\fourvradius}{0.375}
\pgfmathsetmacro{\fivevradius}{0.375}
\pgfmathsetmacro{\sixvspacing}{0.36}
\pgfmathsetmacro{\sevenvradius}{0.36}
\pgfmathsetmacro{\eightvradius}{0.5}
\pgfmathsetmacro{\topgspacing}{0.15}
\pgfmathsetmacro{\btmgspacing}{0.15}

\newcommand\threevertex{%
\foreach \n/\a in {1/90,2/210,3/330}{\node (\n) at (\a:\threevradius) {};}%
\path (1) -- ++(90:\topgspacing) node[draw=none,fill=none,minimum width=0pt] () {};%
\path (2) -- ++(270:\btmgspacing) node[draw=none,fill=none,minimum width=0pt] () {};%
}
\newcommand\fourvertex{%
\foreach \n/\a in {1/135,2/45,3/225,4/315}{\node (\n) at (\a:\fourvradius) {};}%
\path (1) -- ++(90:\topgspacing) node[draw=none,fill=none,minimum width=0pt] () {};%
\path (3) -- ++(270:\btmgspacing) node[draw=none,fill=none,minimum width=0pt] () {};%
}
\newcommand\fivevertex{%
\foreach \n/\a in {1/135,2/45,3/225,4/315}{\node (\n) at (\a:\fivevradius) {};} \node (5) at (0,0) {};%
\path (2) -- ++(90:\topgspacing) node[draw=none,fill=none,minimum width=0pt] () {};%
\path (3) -- ++(270:\btmgspacing) node[draw=none,fill=none,minimum width=0pt] () {};%
}
\newcommand\sixvertex{%
\node (1) at (0,1.2) {}; 
\def\lastn{1}\foreach \n/\a [remember=\n as \lastn] in%
    {2/0,3/315,4/225,5/180,6/135}{\path (\lastn) -- ++(\a:\sixvspacing) node (\n) {};}%
\path (1) -- ++(90:\topgspacing) node[draw=none,fill=none,minimum width=0pt] () {};%
\path (5) -- ++(270:\btmgspacing) node[draw=none,fill=none,minimum width=0pt] () {};%
}
\newcommand\sevenvertex{%
\foreach \n/\a in {1/180,2/135,3/90,4/45,5/0,6/300,7/240}{\node (\n) at (\a:\sevenvradius) {};}%
\path (3) -- ++(90:\topgspacing) node[draw=none,fill=none,minimum width=0pt] () {};%
\path (6) -- ++(270:\btmgspacing) node[draw=none,fill=none,minimum width=0pt] () {};%
}
\newcommand\eightvertex{%
\foreach \n/\a in {1/180,2/135,3/90,4/45,5/0,6/315,7/270,8/225}{\node (\n) at (\a:\eightvradius) {};}%
\path (3) -- ++(90:\topgspacing) node[draw=none,fill=none,minimum width=0pt] () {};%
\path (7) -- ++(270:\btmgspacing) node[draw=none,fill=none,minimum width=0pt] () {};%
}

\makeatletter
\newcommand\addcase[3]{\expandafter\def\csname\string#1@case@#2\endcsname{#3}}
\newcommand\makeswitch[2][]{%
  \newcommand#2[1]{%
    \ifcsname\string#2@case@##1\endcsname\csname\string#2@case@##1\endcsname\else#1\fi%
  }%
}
\makeatother

\newcommand{\ctikz}[1]{\ensuremath{\vcenter{\hbox{#1}}}}

\makeswitch[Missing graph]\ngraph

\addcase\ngraph{3.1}{
}

\def\fl[#1]{\left\lfloor\frac{#1}{2}\right\rfloor}
\def\flb[#1]{\left\lfloor\tfrac{#1}{2}\right\rfloor}

\usepackage{amsmath,amsfonts,latexsym,amssymb}

\usepackage{graphicx}

\usepackage[colorlinks=true,citecolor=blue,linkcolor=blue,urlcolor=blue]{hyperref}

\Volume{XY(2)}
\Year{2020}
\firstpage{1}
\lastpage{91}
\Received{today}
\Revised{tomorrow}
\runninghead{K. Clancy et al.\,/\,Australas. J. Combin. XY\,(2) (2020), 1--89}

\numberwithin{equation}{section}

\newtheorem{theorem}{Theorem}[section]
\newtheorem{corollary}[theorem]{Corollary}

\newtheorem{conjecture}[theorem]{Conjecture}

\newcommand{\arxiv}[1]{\href{http://arxiv.org/abs/#1}{\texttt{arXiv:#1}}}

\begin{document}



\title{A survey of graphs with known or bounded crossing numbers\\{\small (Version 2.0, Dec 8th 2021)}}


\author{Kieran Clancy \quad Michael Haythorpe\thanks{Corresponding author} \quad Alex Newcombe}
\Addr{College of Science and Engineering\\
Flinders University\\
1284 South Road, Tonsley 5042\\
Australia \\
{\tt kieran.clancy@flinders.edu.au}\\
{\tt michael.haythorpe@flinders.edu.au}\\
{\tt alex.newcombe@flinders.edu.au}}


\maketitle

\begin{abstract}
We present, to the best of the authors' knowledge, all known results for the (planar) crossing numbers of specific graphs and graph families. The results are separated into various categories: specifically, results for general graph families, results for graphs arising from various graph products, and results for recursive graph constructions.

\end{abstract}

\tableofcontents


\section{Introduction}

The {\em crossing number} of a graph is the minimum number of crossings over all possible {\em drawings} of $G$. There are many subtleties to consider with this statement, such as, what defines a crossing? What defines a drawing? In a broad survey on the variants of crossing numbers, Marcus Schaefer discusses these subtleties in excellent detail \cite{schaefersurvey}. Here, we are focused only on the standard crossing number in the plane, and as such, simplified definitions suffice. A graph $G$ has the vertex set $V(G)$ and edge set $E(G)$ and a {\em drawing} is a representation of $G$ in the plane. Vertices are represented as distinct points and each edge $e=\{u,v\}$ is represented as a continuous arc connecting the points associated with $u$ and $v$ in such a way that the interior of the arc does not contain any points associated with vertices. In addition, the interiors of the arcs are only allowed to intersect at a finite number of points and such that each intersection is strictly a crossing between the edges, as opposed to the edges touching and then not crossing. The intersections between the arcs are the {\em crossings} of the drawing and the {\em crossing number} of a graph is denoted by $cr(G)$ and is the minimum number of crossings over all possible drawings of $G$. In what follows, when no confusion is possible, we shall refer to the arcs and points given by a drawing as the `edges' and `vertices' of the drawing.



We do not attempt to survey the vast history of crossing numbers, which has been recounted brilliantly in several places, including \cite{schaefersurvey,schaefer2018,beinekewilson2010}, but we briefly mention three of the most influential lines of research which continue to inspire researchers today. The first line of research is the initial work into the crossing number of complete and complete bipartite graphs, which can be read about in \cite{beinekewilson2010}. The subsequent development of Zarankiewicz' Conjecture and the Harary-Hill Conjecture, which both remain largely unsolved, continue to propel crossing numbers to the forefront of research in topological graph theory. The second line of research is Frank Leighton's development of new techniques for bounding the crossing number, including the discovery of the famous Crossing Lemma \cite{leighton1983}. Additionally, Leighton's work made an important connection between crossing numbers and VLSI design, and this continues to attract valuable contributions from the computer science community. The third line of research is Garey and Johnson's proof that the {\em crossing number problem} is NP-complete \cite{gareyjohnson1983}. Specifically, the version of the problem described in \cite{gareyjohnson1983} asks whether a given graph has crossing number less than or equal to a given integer $k$. As has been common with many NP-complete problems, this continues to inspire research into the practical hardness of the problem as well as algorithmic approaches for computing crossing numbers.





To date, crossing numbers have only been determined for a small number of graph families. For some other graph families, bounds have been established for which the upper bounds are usually conjectured (but not proved) to be exact. In this survey, we summarise all such published results for the standard crossing number in the plane, along with references.

Our motivation for producing and maintaining this survey is threefold. The first motivation is that there appears to be many pockets of active research occurring in this area, but in many cases it seems that researchers in one pocket are unaware of similar research being conducted in another. As a consequence, a number of results have been proved multiple times by a number of different authors. Here, for the first time, we gather all of these results, from the early research into this topic by Richard Guy, to the extensive work on crossing numbers of graph products by Mari\'{a}n Kle\v{s}\v{c}, and the vast field of results published in the Chinese mathematical literature, together in one document. Whenever possible, we have attempted to give credit to the author who published the result first, as well as a summary of the partial results that led to larger results. In cases where two sets of authors independently published a result in the same year, we have credited both sets of authors.

The second motivation for this survey is to highlight the remaining gaps in the literature, where results remain to be determined. These are perhaps best illustrated by the various tables of results for Cartesian and join products of small graphs with paths, cycles, stars and discrete graphs, which have been compiled from dozens of individual publications, and for which some scant holes still remain. Our hope is that this will help researchers to focus on the remaining unsolved problems in this field.

The third motivation for this survey is to provide a comprehensive set of useful instances for benchmarking purposes. One of the few benchmark sets which has been repeatedly used in the crossing minimisation literature is known as the KnownCR instances, which were originally collected by Gutwenger \cite{gutwenger2010} in Section 4.3.2 of his thesis.  The KnownCR set of instances can now be expanded upon considerably for use in future research.


Although, to the best of our knowledge, no other similarly extensive survey of this kind exists, there are a number of other valuable resources regarding crossing numbers. In particular, we make note of a few here: Since 2011, Marcus Schaefer \cite{schaefersurvey} has maintained a dynamic survey of different variants of the crossing number. In 2018, Schaefer \cite{schaefer2018} also released an excellent book about various aspects of the crossing number problem.  For over a decade, Imrich Vrt'o \cite{vrto2014} has maintained a very extensive, but unannotated, bibliography of papers relating to crossing numbers; the latest update came in 2014. A fascinating discussion of the history of the problem was presented by Lowell Beineke and Robin Wilson in 2010 \cite{beinekewilson2010}. In 2005, Wynand Winterbach \cite{winterbach2005} gave an excellent summary of the state of research into crossing numbers at the time in his master's thesis. In 2010, Yuanqiu Huang and Jing Wang \cite{huangwang2010} published a survey paper in Chinese, with particular focus on results obtained by Chinese authors.

If one wishes to actually solve instances, Markus Chimani and Tilo Wiedera \cite{chimaniwiedera2016} produced a mixed-integer linear program in 2016 that is able to compute crossing numbers for small graphs, complete with a proof file, and produced an online interface, Crossing Number Web Compute \cite{chimaniwiederasite2016}, where researchers can submit their own instances. Alternatively, if a fast crossing minimisation heuristic is desired, we refer readers to the crossing minimisation heuristics built into the Open Graph Drawing Framework (OGDF) package \cite{chimanietal2014}, or alternatively to our algorithm, QuickCross \cite{clancyetal2019}, which is available for download at \href{http://fhcp.edu.au/quickcross}{http://fhcp.edu.au/quickcross}.

Care has been taken to ensure no symbol is used to represent two different concepts throughout this survey, even though in rare cases this means using non-standard notation. These are pointed out as they occur. A full glossary of symbols used in this survey is given in \hyperref[app-glossary]{Appendix \ref{app-glossary}}.

Wherever practical, we have attempted to verify every result in this survey in the following way: We have generated some moderate-sized instances from each family, and used QuickCross to attempt to find a drawing with the proposed number of crossings. In the vast majority of cases, the best result from QuickCross matched the proposed number of crossings. In cases where QuickCross was able to find a drawing with fewer crossings than proposed by a paper, that drawing itself is proof that the result is incorrect. In cases where QuickCross was only able to find drawings with more crossings than a result suggested, we searched for minimal counterexamples using Crossing Number Web Compute. We have collated all of these incorrect results in \hyperref[app-error]{Appendix \ref{app-error}}, not as a means of disparaging the authors, but rather to ensure that other researchers don't use these results as basis for subsequent proofs.


As new results are being published frequently, our intention is to keep this survey regularly updated, with the latest version always available at \arxiv{1901.05155}, and periodic updates published in Australasian Journal of Combinatorics. We welcome any correspondence alerting us to results which we have either neglected to include, or which have been discovered since the most recent update of this survey. As a general rule, we only include results which have appeared in a peer-reviewed journal other than in exceptional circumstances. Please send any such correspondence to the corresponding author.

\subsection{Asterisked results}

In an attempt to include as many known results as possible, we have considered results from any recognised journals or University periodicals. However, some of these journals impose no peer review, or that which does occur is inadequate. As such, the results contained within cannot be relied upon, either in their own right, or as a basis for subsequent proofs. Indeed, we have encountered many such cases where the proofs are either incorrect, or incomplete.

To address this, we have marked all results appearing within such journals with an asterisk. We emphasise that this determination is solely based on the journal in which the publication appears, rather than the quality of the publication itself. In the case that the result has been subsequently proved in a fully refereed journal, we have removed the asterisk and cited only the latter source. Again, marking results in this way is not intended to disparage the authors, but rather to highlight results which should be revisited and submitted to thorough peer review. Our hope is that, in time, we will be able to replace all of the asterisked results in this way.


\subsection{General results}

Throughout this survey, we will only include results for specific graphs or graph families. There are bounds for the crossing number of general graphs, which we discuss briefly here.

For graphs with sufficiently many edges, the {\em Crossing Lemma} provides a lower bound on the crossing number which depends upon a constant $c$. It was independently discovered by Leighton \cite{leighton1983} and Ajtai et al. (1982) \cite{ajtaietal1982}:

\begin{theorem}[Leighton, 1983 \cite{leighton1983}, Ajtai et al.,\ 1982 \cite{ajtaietal1982}]There is an absolute constant $c > 0$ such that for every graph $G$ with $n$ vertices and $m \geq 4n$ edges,
\[ cr(G) \geq \frac{cm^3}{n^2}. \]\end{theorem}

The Crossing Lemma is tight, other than for the choice of $c$, and it was originally shown that it holds for $c = \frac{1}{100}$. This was improved to $\frac{1}{64}$ by Chazelle, Sharir and Welzl in an email conversation summarised in \cite{aignerziegler2010}. Further improvements can be found if the number of edges in the graph is increased. Pach and T\'{o}th (1997) \cite{pachtoth1997} showed that $c$ can be increased to $1/33.75$ if $m \geq 7.5n$. Later, Pach et al.\ (2006) \cite{pachetal2006} improved this by showing that $c$ can be increased to $1024/31827$ (roughly $1/31.1$) if $m \geq 6.4375n$. Montaron (2005) \cite{montaron2005} determined various values of $c$ depending on the ratio of $m$ and $n$. Finally, using a different approach, Ackerman (2019) \cite{ackerman2019} showed that if $m \geq 6.95n$, $c$ can be increased to $\frac{1}{29}$.

%

In 2007, at an AMS special session in Chicago, Albertson stated a conjecture about graphs with chromatic number $n$:

\begin{conjecture}[Albertson, 2007] Consider any graph $G$ with chromatic number $n$. Then,
\[ cr(G) \geq cr(K_n). \]\end{conjecture}

The bound is obviously tight since $K_n$ has chromatic number $n$. The result is trivially true for $n \leq 4$, since $cr(K_n) = 0$ for these cases. For $n = 5$, it is equivalent to the four colour theorem. Oporowski and Zhao (2009) \cite{oporowskizhao2009} verified the case for $n = 6$. Albertson et al.\ (2010) \cite{albertsonetal2010} then further verified the conjecture for $n = 7, 8, 9, 10, 11, 12$. Bar\'{a}t and T\'{o}th (2010) \cite{barattoth2010} verified the cases $n = 13, 14, 15, 16$, and Ackerman (2019) \cite{ackerman2019} verified the cases $n = 17, 18$. For $n \geq 19$ the conjecture is still open, although Ackerman \cite{ackerman2019} did show that for $n = 19$, counterexamples could only exist for $G$ containing either 37 or 38 vertices.

\subsection{Crossing critical graphs}
A graph $G$ is {\em c-crossing-critical} if $cr(G)\geq c$, but every proper subgraph $H$ of $G$ has $cr(H)<c$. Crossing critical graphs are an important family of graphs for which the crossing number is known (or bounded) and they form a large topic on their own with independent theoretical motivations. We direct the interested reader to \cite{bokaletal2019} for an overview of the results from the most influential line of research in this area. We intend to include a more thorough discussion of crossing-critical graphs in a future iteration of this survey.

\section{Specific graphs and graph families}

\subsection{Complete multipartite graphs}
\begin{figure}[h!]\begin{center}\begin{tikzpicture}[largegraph,scale=0.5]
\foreach \n in {0,...,3}{
  \pgfmathtruncatemacro{\nextn}{2*\n-3}
  \node (a\n) at (0,\nextn) {};
}
\foreach \n in {0,...,4}{
  \pgfmathtruncatemacro{\nextn}{2*\n-4}
  \node (b\n) at (4,\nextn) {};
}
\foreach \n in {0,...,3}{
  \foreach \m in {0,...,4}{
    \draw (a\n)  -- (b\m);
  }
}\end{tikzpicture}$\;\;\;\;\;\;\;\;\;\;\;\;\;\;\;\;\;\;\;\;\;\;\;\;\;$
\begin{tikzpicture}[largegraph,scale=0.5]
\foreach \n in {0,...,3}{
  \pgfmathtruncatemacro{\nextn}{2*\n-3}
  \node (a\n) at (0,\nextn) {};
}
\foreach \n in {0,...,4}{
  \pgfmathtruncatemacro{\nextn}{2*\n-5}
  \node (b\n) at (\nextn,0) {};
}
\foreach \n in {0,...,3}{
  \foreach \m in {0,...,4}{
    \draw (a\n)  -- (b\m);
  }
}\end{tikzpicture}\caption{The complete bipartite graph $K_{4,5}$, with the second drawing in the style first described by Zarankiewicz.\label{fig-K45}}\end{center}\end{figure}
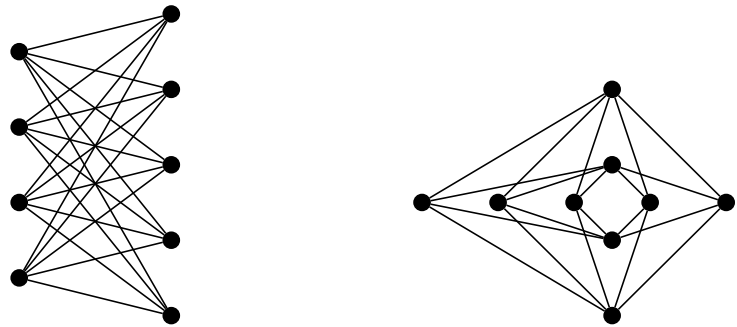

The complete $n$-partite graph $K_{m_1,m_2,\dots,m_n}$ is the graph on $\sum_{i=1}^n m_i$ vertices defined as follows. The vertices are partitioned into disjoint sets $V_1,V_2,\dots,V_n$, such that $|V_i| = m_i$ for $i=1,2,\dots n$. An edge exists between two vertices precisely when one vertex is in $V_i$ and the other is in $V_j$ where $i\neq j$. An example of the complete bipartite graph $K_{4,5}$ is displayed in \hyperref[fig-K45]{Figure \ref{fig-K45}} in two drawings; the latter drawing is optimal with respect to the number of crossings.

The crossing numbers for complete multi-partite graphs are only known for a small number of cases, however, a general upper bound is known and is due to Harborth \cite{harborth1970}. Harborth's upper bound has matched all of the known exact crossing numbers for these graphs thus far.

\begin{theorem}[Harborth,\ 1970 \cite{harborth1970}]For the complete $n$-partite graph $K_{x_1,x_2,\dots,x_n}$ define $s=\sum_{i=1}^n x_i$ and $c=\sum_{i=1}^n (x_i \mod 2)$.  Then the following holds:
\[ cr(K_{x_1,x_2,\dots,x_n}) \leq \frac{1}{8}\Bigg(\sum_{\substack{1\leq i<j<k<\ell \leq n}} \hspace{-0.4cm}3x_i x_j x_k x_\ell+3{ \flb[c] \choose 2} - \flb[c] \hspace{-0.5cm}\sum_{\substack{1\leq i<j \leq n \\ i=j=0 \pmod{2}}} \hspace{-0.4cm}x_i x_j \] \[ \;\;\;\;\, - \flb[c-1] \hspace{-0.4cm}\sum_{\substack{1\leq i<j \leq n \\ i\neq j\pmod{2}}}\hspace{-0.4cm} x_i x_j - \flb[c-2]  \hspace{-0.5cm}\sum_{\substack{1\leq i<j \leq n \\ i=j=1 \pmod{2}}}\hspace{-0.4cm} x_i x_j\Bigg) \] \[+ \sum_{i=1}^n \flb[x_i]\flb[x_i-1]\flb[s-x_i]\flb[s-x_i-1]-\hspace{-0.2cm}\sum_{\substack{1\leq i<j \leq n}}\hspace{-0.2cm}\flb[x_i]\flb[x_i-1]\flb[x_j] \lfloor \tfrac{x_j-1}{2}\rfloor .\]\label{thm-harborth}\end{theorem}

\subsubsection{Complete bipartite graphs}\label{sec-complete-bipartite}


One of the seminal results for crossing numbers was by Zarankiewicz, who attempted to solve the crossing number problem specifically for complete bipartite graphs; this special case is known as Tur\'{a}n's brick factory problem. Zarankiewicz claimed to have proved an exact bound \cite{zarankiewicz1955}, but his proof was subsequently found to contain an error \cite{guy1969}. Nonetheless, the result of Zarankiewicz still provides an upper bound:

\begin{theorem}[Zarankiewicz, 1955 \cite{zarankiewicz1955}]The crossing number of the complete bipartite graph $K_{m,n}$ is bounded as follows:
\[ cr(K_{m,n}) \leq \fl[n]\fl[n-1]\fl[m]\fl[m-1]. \]\label{thm-zarankiewicz}\end{theorem}

Due to isomorphism, it is clear that $cr(K_{m,n}) = cr(K_{n,m})$ and so in what follows, we will assume that $m \leq n$.

The upper bound has been shown to coincide with the true crossing number for some small cases. In Guy's \cite{guy1969} paper refuting Zarankiewicz's proof, it was shown that the result holds for $m \leq 4$. In 1971, Kleitman \cite{kleitman1971} verified it for $K_{5,n}$. In 1993, the cases of $K_{7,7}$ and $K_{7,9}$ were verified by Woodall \cite{woodall1993} to coincide with Zarankiewicz's formula. In \cite{guy1969}, it was shown that if the result holds for $K_{m,n}$, such that $m$ is odd, then it holds for $K_{m+1,n}$. Hence, Kleitman's result also verifies $K_{6,n}$, and Woodall's results also settle the cases $K_{7,8}$, $K_{7,10}$, $K_{8,8}$, $K_{8,9}$ and $K_{8,10}$. Despite detailing a flaw in Zarankiewicz's proof, Guy conjectured that the result would still hold in all cases and this has come to be known as Zarankiewicz' Conjecture:

\begin{conjecture}[Zarankiewicz, \cite{zarankiewicz1955}, Guy, 1969 \cite{guy1969}]\label{conj-zar}\hyperref[thm-zarankiewicz]{Theorem~\ref{thm-zarankiewicz}} holds with equality.\end{conjecture}

Several asymptotic lower bounds have been proved for complete bipartite graphs. The earliest results by Kleitman \cite{kleitman1971} imply asymptotic lower bounds, and these were improved upon in 2003 by Nahas \cite{nahas2003}. In 2006, de Klerk et al.\ \cite{deklerketal2006} used semi-definite programming methods to make significant progress, which were refined in 2007 by de Klerk et al.\ \cite{deklerketal2007}, who proved the following:


\begin{theorem}[de Klerk et al.,\ 2007 \cite{deklerketal2007}]For the complete bipartite graph $K_{m,n}$ with $m \geq 9$,
\[ \lim_{n \rightarrow \infty} \frac{cr(K_{m,n})}{\fl[n]\fl[n-1]\fl[m]\fl[m-1]} \geq \frac{0.8594m}{m-1}. \]\label{knm-lb}\end{theorem}

In 2013, Norin and Zwols \cite{norin2013} announced that the 0.8594 in \hyperref[knm-lb]{Theorem \ref{knm-lb}} could be replaced by 0.905, but this was never published.

In 2013, Christian et al.\ \cite{christianetal2013} gave a function $N_0(m)$ and showed that for each $m$, if \hyperref[conj-zar]{Conjecture~\ref{conj-zar}} can be confirmed for all $n \leq N_0(m)$, then it is true for all $n$. Hence for each $m$, in order to confirm \hyperref[conj-zar]{Conjecture~\ref{conj-zar}}, only finitely many cases need to be considered. Unfortunately, the function $N_0(m)$ is not practical even for small $m$.

\begin{theorem}[Christian et al.,\ 2013 \cite{christianetal2013}]Let $m$ be a positive integer. Define $N_0(m)$ as follows:
\[ N_0(m) := \left((m!)!\left(2\fl[m]\fl[m-1]\right)^{m!}\right)^4. \]

If \hyperref[conj-zar]{Conjecture~\ref{conj-zar}} is true for all $n \leq N_0(m)$, it is true for all $n$.\end{theorem}

\subsubsection{Complete tripartite graphs}\label{sec-complete_tripartite}

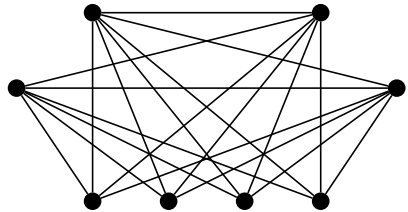
\begin{figure}[h!]\begin{center}\begin{tikzpicture}[largegraph,scale=0.5]
\node (a1) at (-5,3) {};
\node (a2) at (-3,5) {};
\node (b1) at (5,3) {};
\node (b2) at (3,5) {};
\node (c1) at (-3,0) {};
\node (c2) at (-1,0) {};
\node (c3) at (1,0) {};
\node (c4) at (3,0) {};

\draw (a1) -- (b1);
\draw (a1) -- (b2);
\draw (a2) -- (b1);
\draw (a2) -- (b2);
\foreach \n in {1,2}{
  \foreach \m in {1,2,3,4}{
    \draw (a\n) -- (c\m);
    \draw (b\n) -- (c\m);
  }
}\end{tikzpicture}\caption{The complete tripartite graph $K_{2,2,4}$.\label{fig-K224}}\end{center}\end{figure}

Complete tripartite graphs have their vertices partitioned into three sets. An example of $K_{2,2,4}$ is displayed in \hyperref[fig-K224]{Figure \ref{fig-K224}}.

The crossing number of $K_{a,b,n}$ has been determined for all cases where $a + b \leq 6$, except for $K_{3,3,n}$. The known cases are listed in \hyperref[tab-tripartite]{Table \ref{tab-tripartite}}, along with the publications where they were first discovered. Note that $K_{1,1,n}$ is planar. In a few cases, the verification came from the perspective of join products, since $K_{a,b} + D_n = K_{a,b,n}$, where $D_n$ is the discrete graph on $n$ vertices. See \hyperref[sec-join]{Section \ref{sec-join}} for more information and results on join products. Ho was the first author to explicitly consider $cr(K_{1,2,n})$, although it could be seen as a direct consequence of results from some other authors. First, $P_2 \Box S_n$ is a subdivision of $K_{1,2,n}$, and the crossing number of $P_2 \Box S_n$ was determined by Bokal (2007) \cite{bokal2007}. Second, $K_{1,2,n}$ is a subgraph of $P_2 + P_n$, and a supergraph of $K_{3,n}$, both of which were shown to have crossing number $\fl[n]\fl[n-1]$ by Kle\v{s}\v{c} (2007) \cite{klesc2007} and Guy (1969) \cite{guy1969} respectively. The result for $K_{1,4,n}$ was also independently proved by Ho (2008) \cite{ho2008}.

{\renewcommand{\arraystretch}{1.5}
\begin{table}[htbp]\begin{center}\caption{Crossing numbers of complete tripartite graphs.\label{tab-tripartite}}\smallskip\smallskip$\begin{array}{|c|l|l|}\hline
\rule{0pt}{2.3ex} \text{{\bf Graph family}}       & \text{{\bf Crossing number}} & \text{{\bf Publication}}     \\
\hline K_{1,2,n} & \fl[n]\fl[n-1] & \text{Ho (2008) \cite{ho2008}} \\
\hline K_{1,3,n} & 2\fl[n]\fl[n-1] + \fl[n] & \text{Asano (1986) \cite{asano1986}} \\
\hline K_{1,4,n} & 4\fl[n]\fl[n-1] + 2\fl[n]  & \text{Huang and Zhao (2008) \cite{huangzhao2008}} \\
\hline K_{1,5,n} & 6\fl[n]\fl[n-1] + 4\fl[n] \: \Asterisk& \text{Mei and Huang (2007) \cite{meihuang2007} } \Asterisk\\
\hline K_{2,2,n} & 2\fl[n]\fl[n-1]  & \text{Kle\v{s}\v{c} and Schr\"{o}tter (2011) \cite{klescschrotter2011}} \\
\hline K_{2,3,n} & 4\fl[n]\fl[n-1] + n & \text{Asano (1986) \cite{asano1986}} \\
\hline K_{2,4,n} & 6\fl[n]\fl[n-1] + 2n  & \text{Ho (2013) \cite{ho2013}} \\
\hline\end{array}$\end{center}\end{table}
{\renewcommand{\arraystretch}{1}





In 2004, Ho \cite{ho2004} showed that $cr(K_{3,3,n})$ can be determined if \hyperref[conj-zar]{Conjecture~\ref{conj-zar}} holds for $m = 7$ and $n \leq 20$. To date this is not known to be true, and so the crossing number of $K_{3,3,n}$ can only be given as a conjecture.

\begin{conjecture}[Ho, 2004 \cite{ho2004}]For $n \geq 1$, \[cr(K_{3,3,n}) = 6\fl[n]\fl[n-1] + 2\fl[3n] + 1.\]\end{conjecture}

There are some other results depending on the truth of \hyperref[conj-zar]{Conjecture~\ref{conj-zar}}. In particular, Huang and Zhao (2006) \cite{huangzhao2006,huangzhao2006_2} proved that $cr(K_{1,6,n})$ can be determined if \hyperref[conj-zar]{Conjecture~\ref{conj-zar}} holds for $K_{7,k}$ for all $k \geq 1$, and likewise that $cr(K_{1,8,n})$ can be determined if \hyperref[conj-zar]{Conjecture~\ref{conj-zar}} holds for $K_{9,k}$ for all $k \geq 1$. In 2008, Wang and Huang \cite{wanghuang2008} proved that $cr(K_{1,10,n})$ can be determined if \hyperref[conj-zar]{Conjecture~\ref{conj-zar}} holds for $K_{11,k}$ for all $k \geq 1$. Finally, Ho (2008) \cite{ho2008} generalised these results by showing the following.

\begin{theorem}[Ho, 2008 \cite{ho2008}] If \hyperref[conj-zar]{Conjecture~\ref{conj-zar}} is true for $K_{2M+1,n}$ then,
\[cr(K_{1,2M,n})=M^2 \fl[n+1]\fl[n]-M\fl[n].\]
\end{theorem}

Currently, \hyperref[conj-zar]{Conjecture~\ref{conj-zar}} is only known to hold for $M \leq 2$.

In 2017, Gethner et al.\ \cite{gethneretal2017} gave asymptotic lower bounds on the crossing number of the balanced complete tripartite graphs $K_{n,n,n}$.  Let $A(n,n,n)$ denote the right hand side of the inequality in \hyperref[thm-harborth]{Theorem \ref{thm-harborth}}, evaluated for $K_{n,n,n}$. Gethner et al.\ gave asymptotic lower bounds involving $A(n,n,n)$.

\begin{theorem}[Gethner et al.\, 2017 \cite{gethneretal2017}] For the balanced complete tripartite\\graphs $K_{n,n,n}$,
\[\lim_{n\rightarrow \infty} \frac{cr(K_{n,n,n})}{A(n,n,n)} \geq 0.666.\]
\end{theorem}

\subsubsection{Complete $4$-partite graphs}

In 2016, Shanthini and Babujee \cite{shanthinibabujee2016_2} $\Asterisk$ showed that the crossing number of $K_{1,1,m,n}$ can be expressed in terms of the crossing number of $K_{m+2,n+2}$:

\begin{theorem}[Shanthini and Babujee, 2016 \cite{shanthinibabujee2016_2}$\Asterisk$]\label{thm-K11mn} For $m,n \geq 1$,
\[ cr(K_{1,1,m,n}) = cr(K_{m+2,n+2}) + \fl[m]\fl[n] - mn. \: \Asterisk \]\end{theorem}

Since Zarankiewicz's conjecture is known to hold for $\min\{m,n\} \leq 6$, \hyperref[thm-K11mn]{Theorem~\ref{thm-K11mn}} settles $cr(K_{1,1,m,n})$ for $m \leq 4$. Each of those four cases had been previously settled. Specifically, the cases $m = 1, 2, 3$ were settled in 2007 by Qian and Huang \cite{qianhuang2007} $\Asterisk$, and the case $m = 4$ was settled in a separate 2016 paper by Shanthini and Babujee \cite{shanthinibabujee2016} $\Asterisk$.


There are also some results known for other $4$-partite graphs. In particular, He and Huang \cite{hehuang2007_2} $\Asterisk$ determined the crossing number of $K_{1,2,2,n}$ and Ho \cite{ho2008_2} determined the crossing number of $K_{2,2,2,n}$.

\begin{theorem}[He and Huang, 2007 \cite{hehuang2007_2} $\Asterisk$]For $n\geq 1$,
\[cr(K_{1,2,2,n}) = 4\fl[n]\fl[n-1] + n + \fl[n]. \: \Asterisk\]\end{theorem}

\begin{theorem}[Ho, 2008 \cite{ho2008_2}]For $n \geq 1$, \[cr(K_{2,2,2,n}) = 6\fl[n]\fl[n-1] + 3n.\]\end{theorem}

\subsubsection{Complete $5$-partite graphs}

In 2009, Ho \cite{ho2009} determined the crossing number for two $5$-partite graphs:

\begin{theorem}[Ho, 2009 \cite{ho2009}]For $n \geq 1$, the following hold:
\begin{align*}cr(K_{1,1,1,1,n}) &= 2\fl[n]\fl[n-1] + n,\\
cr(K_{1,1,1,2,n}) &= 4\fl[n]\fl[n-1] + 2n.\end{align*}\end{theorem}

\subsubsection{Complete $6$-partite graphs}

In 2008, L\"{u} and Huang \cite{luhuang2008} determined the crossing number of $K_{1,1,1,1,1,n}$:

\begin{theorem}[L\"{u} and Huang, 2008 \cite{luhuang2008}] For $n\geq 1$,
\[cr(K_{1,1,1,1,1,n}) = 4\fl[n]\fl[n-1] + 2n + \fl[n] + 1.\]\end{theorem}

\subsubsection{Complete bipartite graphs minus an edge}

In 2011, He et al.\ \cite{heetal2011} $\Asterisk$ considered the complete bipartite graph $K_{m,n}$ minus one edge, denoted as $K_{m,n} \setminus e$. Due to symmetry, it doesn't matter which edge is removed. They were able to settle the case when $m = 3$ or $m = 4$:

\begin{theorem}[He et al.,\ 2011 \cite{heetal2011} $\Asterisk$] For $n \geq 1$, the following hold:
\begin{align*}cr(K_{3,n} \setminus e) & = \phantom{2}\fl[n]\fl[n-1] - \fl[n-1], \: \Asterisk\\
cr(K_{4,n} \setminus e) &= 2\fl[n]\fl[n-1] - \fl[n-1]. \: \Asterisk \end{align*}\end{theorem}

In 2015, Chia and Lee \cite{chialee2015} gave a conjecture for larger $m$, and proved that in addition to $m = 3,4$, the conjecture is also true for the special case of $K_{5,5} \setminus e$:

\begin{conjecture}[Chia and Lee, 2015 \cite{chialee2015}] For $m,n \geq 1$,
\[ cr(K_{m,n} \setminus e) = cr(K_{m,n}) - \fl[m-1]\fl[n-1]. \]
\end{conjecture}

\subsection{Complete graphs}\label{sec-complete}

The complete graph $K_n$ is the graph on $n$ vertices containing an edge between every pair of vertices. Two examples of complete graphs, $K_6$ and $K_8$, are displayed in \hyperref[fig-complete]{Figure \ref{fig-complete}}.

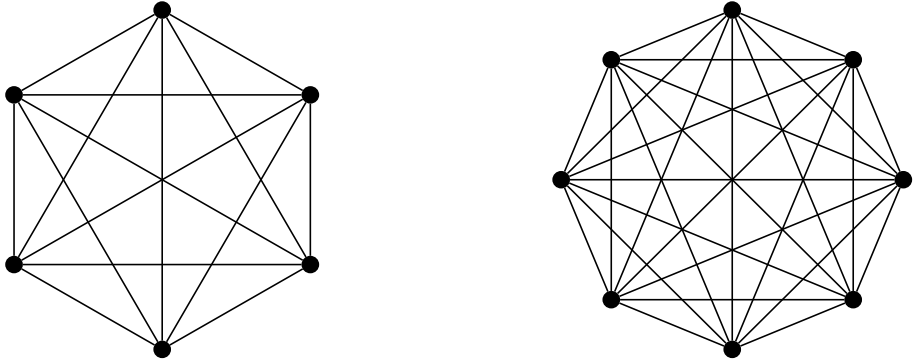
\begin{figure}[h!]\begin{center}\begin{tikzpicture}[largegraph]
\foreach \n in {0,...,5}{
  \node (\n) at (90+\n/6*360:2.25) {};
}
\foreach \n in {0,...,4}{
\pgfmathtruncatemacro{\nextn}{\n + 1}
  \foreach \m in {\nextn,...,5}{
    \draw (\n) -- (\m);
  }
}
\end{tikzpicture}$\;\;\;\;\;\;\;\;\;\;\;\;\;\;\;\;\;\;\;\;\;\;\;\;\;$
\begin{tikzpicture}[largegraph]
\foreach \n in {0,...,7}{
  \node (\n) at (90+\n/8*360:2.25) {};
}
\foreach \n in {0,...,6}{
\pgfmathtruncatemacro{\nextn}{\n + 1}
  \foreach \m in {\nextn,...,7}{
    \draw (\n) -- (\m);
  }
}
\end{tikzpicture}\caption{The complete graphs $K_6$ and $K_8$.\label{fig-complete}}\end{center}\end{figure}

Guy, in a 1960 paper \cite{guy1960}, describes how the problem of determining the crossing number for complete graphs was brought to his attention by Anthony Hill and Prof. C.A. Rogers, but indicates that Paul Erd\H{o}s had been looking at the problem for over twenty years at that point. Guy gave an upper bound which is conjectured to be equal to the crossing number, and showed that if the conjecture holds for odd $n$, then it also holds for $n+1$.

\begin{theorem}[Guy, 1960 \cite{guy1960}]The crossing number of the complete graph $K_n$ is bounded above as follows:
\[ cr(K_n) \leq \frac{1}{4}\fl[n]\fl[n-1]\fl[n-2]\fl[n-3]. \]\label{thm-guy}\end{theorem}

\begin{conjecture}[Guy, 1960 \cite{guy1960}, Harary and Hill, 1963 \cite{harary1963}]\hyperref[thm-guy]{Theorem~\ref{thm-guy}} holds with equality.\label{conj-guycomplete}\end{conjecture}

Guy showed that \hyperref[conj-guycomplete]{Conjecture~\ref{conj-guycomplete}} holds for $n \leq 6$. The same conjecture was also proposed by Harary and Hill \cite{harary1963} at around the same time and then again four years later by Saaty \cite{saaty1964}. In 1972, Guy \cite{guy1972} showed that \hyperref[conj-guycomplete]{Conjecture~\ref{conj-guycomplete}} holds for $n \leq 10$. This was the best known result for more than three decades, until Pan and Richter \cite{panrichter2007} proved in 2007 that the conjecture holds for $n = 11$, and hence for $n = 12$ as well. The remaining cases are still open.

The case of $n = 13$ has been considered closely. McQuillan et al.\ (2015) \cite{mcquillanetal2015} noted that a simple application of Kleitman's parity theorem \cite{kleitman1971} implies that $cr(K_{13})$ must be equal to one of the numbers $\{217, 219, 221, 223, 225\}$; McQuillan et al.\ then proved that $cr(K_{13}) \neq 217$. This result was further improved by \'{A}brego et al.\ (2015) \cite{abregoetal2015} who showed that $cr(K_{13}) \neq 219$ and $cr(K_{13}) \neq 221$. Hence, there are only two possibilities remaining for the crossing number of $K_{13}$; either 223 or 225.

A similar approach as the asymptotic lower bounds for complete bipartite graphs, provides asymptotic lower bounds for complete graphs. In 2019, Balogh et al. \cite{balogh2019} gave the current best version of this lower bound.

\begin{theorem}[Balogh et al.,\ 2019 \cite{balogh2019}]For the complete graph $K(n)$,
\[ \lim_{n \rightarrow \infty} \frac{cr(K_n)}{\frac{1}{4}\fl[n]\fl[n-1]\fl[n-2]\fl[n-3]} > 0.985. \]\label{knlb-thm}\end{theorem}

The constant in \hyperref[knlb-thm]{Theorem \ref{knlb-thm}} was an improvement on the previous work of Norin and Zwols \cite{norin2013} and de Klerk et al.\ \cite{deklerketal2007} who gave constants of 0.905 and 0.83 repectively.

\subsubsection{Complete graphs minus an edge}\label{sec-completeminus}

In 2007, Zheng et al.\ \cite{zhengetal2007}, considered the complete graph with a single edge removed, $K_n \setminus e$. They determined an upper bound for its crossing number, and conjectured that it would coincide with the crossing number. They also proved that the conjecture holds for $n \leq 8$.

\begin{theorem}[Zheng et al.,\ 2007 \cite{zhengetal2007}]For $n\geq1$,
\[ cr(K_n \setminus e) \leq \frac{1}{4}\fl[n+2]\fl[n-1]\fl[n-3]\fl[n-4], \]
with equality known to hold for $n \leq 8$.\label{thm-Knminuse}\end{theorem}

\begin{conjecture}[Zheng et al.,\ 2007 \cite{zhengetal2007}]\hyperref[thm-Knminuse]{Theorem~\ref{thm-Knminuse}} holds with equality.\label{conj-Knminuse}\end{conjecture}

In 2014, Ouyang et al.\ (2014) \cite{ouyangetal2014_2} proved that \hyperref[conj-Knminuse]{Conjecture~\ref{conj-Knminuse}} holds for even $n$ whenever \hyperref[conj-zar]{Conjecture~\ref{conj-zar}} holds for $n-1$. Hence, \hyperref[conj-Knminuse]{Conjecture~\ref{conj-Knminuse}} is currently known to hold for $n \leq 12$.

Chia and Lee \cite{chialee2015} independently discovered results equivalent to the results in \cite{ouyangetal2014_2} around the same time.

\subsubsection{Complete graphs minus a cycle}

In 1973, Guy and Hill \cite{guyhill1973} considered the crossing number of the complement of a cycle $\overline{C}_n$, that is, the complete graph $K_n$ with a simple cycle of length $n$ removed, which is defined for all $n \geq 3$. It is easy to check that $\overline{C}_n$ is planar for $n \leq 6$. They gave an upper bound for the crossing number of $\overline{C}_n$ in general, and conjectured that equality would hold. Guy and Hill also proved that equality does indeed hold for $n \leq 10$.

\begin{theorem}[Guy and Hill, 1973 \cite{guyhill1973}] For $n \geq 3$, the following holds:
\[ cr(\overline{C}_n) \leq \begin{cases}\frac{1}{64}(n-3)^2(n-5)^2, & \text{for odd $n$,}\\\frac{1}{64}n(n-4)(n-6)^2, & \text{for even $n$,}\end{cases} \]
with equality known to hold for $n \leq 10$.\label{thm-guyhill}\end{theorem}

\begin{conjecture}[Guy and Hill, 1973 \cite{guyhill1973}] \hyperref[thm-guyhill]{Theorem~\ref{thm-guyhill}} holds with equality for all $n \geq 3$.\end{conjecture}

From the proof of \hyperref[thm-guyhill]{Theorem~\ref{thm-guyhill}}, it is known that if the theorem holds for a particular odd $n$, then it holds for the next even $n$ as well.

Guy and Hill also determined a lower bound for the crossing number of $\overline{C}_n$.

\begin{theorem}[Guy and Hill, (1973) \cite{guyhill1973}]For $n \geq 9$, the following bounds hold:
\[ cr(\overline{C}_n) \geq \frac{\binom{n}{5}(n-15)(n-17)}{4 \binom{n-4}{3}} > \frac{1}{80}n(n+2)(n-9)(n-20). \]
\label{thm-guyhill2}
\end{theorem}

Note that \hyperref[thm-guyhill2]{Theorem~\ref{thm-guyhill2}} is trivial unless $n > 20$.

\subsection{Circulant graphs}\label{sec-circulant}

The circulant graph $Ci_n(L)$ is the graph on $n$ vertices such that, for $i = 1, \hdots, n$ and each $j$ appearing in the list $L$, vertex $i$ is adjacent to the $(i+j)$-th vertex and the $(i-j)$-th vertex (mod $n$). Two examples of circulant graphs, $Ci_{10}(\{1,3\})$ and $Ci_{11}(\{1,4,5\})$, are displayed in \hyperref[fig-circulant]{Figure \ref{fig-circulant}}. It is well known that $Ci_n(\{1,2\})$ is planar for all even $n$ and has crossing number 1 for all odd $n$. In what follows, we consider circulant graphs for various choices of $n$ and $L$.

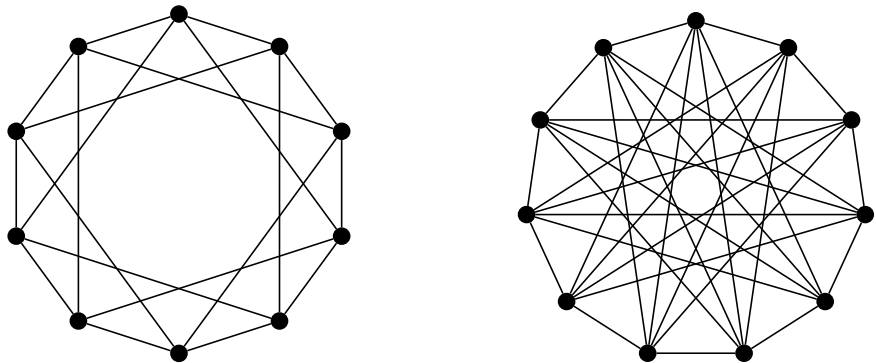
\begin{figure}[h!]\begin{center}\begin{tikzpicture}[largegraph]
\foreach \n in {0,...,9} \node (\n) at (90+\n/10*360:2.25) {};
\foreach \n in {0,...,9}{
  \pgfmathtruncatemacro{\plusone}{mod(\n+1,10)}
  \pgfmathtruncatemacro{\plusthree}{mod(\n+3,10)}
  \draw (\n) -- (\plusone);
  \draw (\n) -- (\plusthree);
}
\end{tikzpicture}\;\;\;\;\;\;\;\;\;\;\;\;\;\;\;\;\;\;
\begin{tikzpicture}[largegraph]
\foreach \n in {0,...,10} \node (\n) at (90+\n/11*360:2.25) {};
\foreach \n in {0,...,10}{
  \pgfmathtruncatemacro{\plusone}{mod(\n+1,11)}
  \pgfmathtruncatemacro{\plusfour}{mod(\n+4,11)}
  \pgfmathtruncatemacro{\plusfive}{mod(\n+5,11)}
  \draw (\n) -- (\plusone);
  \draw (\n) -- (\plusfour);
  \draw (\n) -- (\plusfive);
}
\end{tikzpicture}
\caption{The circulant graphs $Ci_{10}(\{1,3\})$ and $Ci_{11}(\{1,4,5\})$.\label{fig-circulant}}\end{center}\end{figure}


In 1986, Fiorini \cite{fiorini1986} investigated the crossing number of $Ci_n(\{1,3\})$, primarily as a vessel to establish the crossing number of certain generalized Petersen graphs. While he claimed to have determined the crossing number of $Ci_n(\{1,3\})$ for $n \geq 8$, his proof was later shown to contain an error. Nonetheless, his proof established a valid upper bound, and also established equality for $n = 8, 10, 12$:

\begin{theorem}[Fiorini, 1986 \cite{fiorini1986}]For $n \geq 8$, the following holds:
\[ cr(Ci_n(\{1,3\})) \leq \left\lfloor\frac{n}{3}\right\rfloor + (n \hspace*{-0.175cm}\mod 3), \]

holding with equality for $n = 8, 10, 12$.\end{theorem}

In 2004, Yang et al.\ \cite{yangetal2004} proved that equality does indeed hold for all $n \geq 8$:

\begin{theorem}[Yang et al.,\ 2004 \cite{yangetal2004}]For $n \geq 8$, the following holds:
\[ cr(Ci_n(\{1,3\})) = \left\lfloor\frac{n}{3}\right\rfloor + (n \hspace*{-0.175cm}\mod 3). \]\end{theorem}


In 2001, Yang and Zhao \cite{yangzhao2001} $\Asterisk$ determined the crossing number of $Ci_n(\{1,\lfloor n/2 \rfloor\})$ as follows:

\begin{theorem}[Yang and Zhao, 2001 \cite{yangzhao2001} $\Asterisk$]For $n \geq 6$, the following holds:
\[ cr(Ci_n(\{1,\lfloor n/2 \rfloor\})) = 1. \: \Asterisk \]\end{theorem}

For even $n$, $Ci_n(\{1,\lfloor n/2 \rfloor\})$ is isomorphic to the M\"{o}bius ladder graph on $n$ vertices, for which the crossing number was determined by Guy and Harary \cite{guyharary1967} in 1967.


In 2006, Lin et al.\ \cite{linetal2006} considered the circulant graph $Ci_n(\{1,\left\lfloor n/2 \right\rfloor -1\})$. They were able to determine the crossing number for even $n \geq 8$, and discover upper bounds for odd $n \geq 13$:

\begin{theorem}[Lin et al.,\ 2006 \cite{linetal2006}]For even $n \geq 8$, the following holds:

\[ cr(Ci_n(\{1,n/2 -1\})) = n/2. \]\end{theorem}


\begin{theorem}[Lin et al.,\ 2006 \cite{linetal2006}]For odd $n \geq 13$, the following holds:
\[ cr(Ci_n(\{1,(n-1)/2 -1\})) \leq \left\{\begin{array}{lcr}(n+1)/2, & \mbox{ for } n = 1 \pmod{8},\\
(n+3)/2, & \mbox{ for } n = 3 \pmod{8},\\
(n+3)/2, & \mbox{ for } n = 5 \pmod{8},\\
(n+1)/2, & \mbox{ for } n = 7 \pmod{8}.\end{array}\right. \]\end{theorem}


In 2005, Salazar \cite{salazar2005} proved that the crossing number of $Ci_n(\{1,k\})$ for $k \geq 5$ and $n \geq k^4$ can be bounded above and below and by functions of $n$ and $k$:

\begin{theorem}[Salazar, 2005 \cite{salazar2005}]For $k \geq 5$ and $n \geq k^4$, the following bounds hold:
\[  \left(1 - \frac{4}{k}\right)n + (4k^2 + 1 - k^3) \leq cr(Ci_n(\{1,k\})) \leq \left(1 - \frac{2}{k}\right)n + \frac{k^2 + k + 2}{2}. \]\end{theorem}


In 2005, Lin et al.\ \cite{linetal2005} considered the circulant graph $Ci_{mk}(\{1,k\})$, for integers $m \geq 3$, $k \geq 3$. They were able to determine the crossing number for $m = 3$, and established upper bounds for larger $m$:

\begin{theorem}[Lin et al.,\ 2005 \cite{linetal2005}]For $m,k \geq 3$, the following hold:
\begin{align*}&cr(Ci_{3k}(\{1,k\}) = k,\\
&cr(Ci_{4k}(\{1,k\}) \leq 2k + 1, \text{ and for $m\geq 5$,}\\
&cr(Ci_{mk}(\{1,k\}) \leq \min \{(m-2)(k+1) - 1, m(k-2)\}.\end{align*}
\end{theorem}


In 2005, Ma et al.\ \cite{maetal2005_2} considered the circulant graph $Ci_{2m+2}(\{1,m\})$ for $m \geq 2$:

\begin{theorem}[Ma et al.,\ 2005 \cite{maetal2005_2}]The following holds:
\[ cr(Ci_{2m+2}(\{1,m\})) = \begin{cases} 0, & \mbox{ for } m=2, \\ m+1, & \mbox{ for } m \geq 3. \end{cases}\]\end{theorem}


In 2007, Ho \cite{ho2007} considered the circulant graph $Ci_{3m+1}(\{1,m\})$ for $m \geq 3$:

\begin{theorem}[Ho, 2007 \cite{ho2007}]For $m \geq 3$, the following holds:
\[ cr(Ci_{3m+1}(\{1,m\})) = m+1. \]\end{theorem}


In 2008, Wang and Huang \cite{wanghuang2008_3} $\Asterisk$ considered the circulant graph $Ci_{3m-1}(\{1,m\})$ and derived the following upper and lower bounds, and conjectured that the upper bound would provide the correct crossing number.

\begin{theorem}[Wang and Huang, 2008 \cite{wanghuang2008_3} $\Asterisk$]For $m \geq 3$, the following holds:
\[ m \leq cr(Ci_{3m-1}(\{1,m\})) \leq m+1. \: \Asterisk \]\end{theorem}

\begin{conjecture}[Wang and Huang, 2008 \cite{wanghuang2008_3}]For $m \geq 3$, the following holds:
\[ cr(Ci_{3m-1}(\{1,m\})) = m+1. \]\end{conjecture}

\subsection{Generalized Petersen graphs}\label{sec-GP}

The generalized Petersen graphs were first studied by Coxeter \cite{coxeter1950} and later named by Watkins \cite{watkins1969}. The generalized Petersen graph $GP(n,k)$ is constructed by taking the union of the cycle graph $C_n$ and the circulant graph $Ci_n(k)$ (for which the definition is given in \hyperref[sec-circulant]{Section \ref{sec-circulant}}) where $n\geq3$, and then connecting the corresponding vertices in each by an edge. The result is a 3-regular graph containing $2n$ vertices, except in the special case when $n = 2k$ which corresponds to the M\"{o}bius ladder graph of size $n$ with the $n/2$ internal edges subdivided twice. Two examples of generalized Petersen graphs, $GP(9,2)$ and $GP(11,3)$, are displayed in \hyperref[fig-GP92]{Figure \ref{fig-GP92}}. It is common in the literature to use $P(n,k)$, rather than $GP(n,k)$, to denote generalized Petersen graphs. To avoid confusion with path graphs, we use the latter here.

\begin{figure}[h!]\begin{center}\begin{tikzpicture}[largegraph]
\foreach \n in {0,...,8}{
  \node (\n) at (90+\n/9*360:1.5) {};
  \node (o\n) at (90+\n/9*360:2.25) {};
  \draw (\n) -- (o\n);
}
\foreach \n in {0,...,8}{
  \pgfmathtruncatemacro{\nextn}{mod(\n+1,9)}
  \pgfmathtruncatemacro{\nextnextn}{mod(\n+2,9)}
  \draw (o\n) -- (o\nextn);
  \draw (\n) -- (\nextnextn);
}
\end{tikzpicture}\;\;\;\;\;\;\;\;\;\;\;\;\;\;\;\;\;\;
\begin{tikzpicture}[largegraph]
\foreach \n in {0,...,10}{
  \node (\n) at (90+\n/11*360:1.5) {};
  \node (o\n) at (90+\n/11*360:2.25) {};
  \draw (\n) -- (o\n);
}
\foreach \n in {0,...,10}{
  \pgfmathtruncatemacro{\nextn}{mod(\n+1,11)}
  \pgfmathtruncatemacro{\nextnextn}{mod(\n+3,11)}
  \draw (o\n) -- (o\nextn);
  \draw (\n) -- (\nextnextn);
}
\end{tikzpicture}\caption{The generalized Petersen graphs $GP(9,2)$ and $GP(11,3)$.\label{fig-GP92}}\end{center}\end{figure}
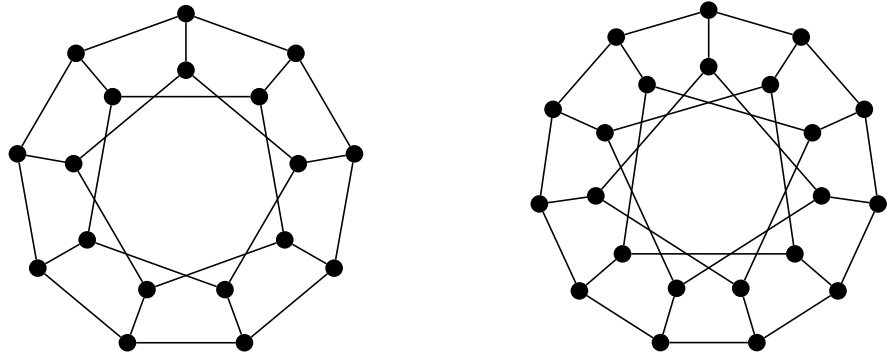

A number of results are known about the crossing numbers of certain subfamilies of generalized Petersen graphs. Due to symmetry, $GP(n,k)$ is isomorphic to $GP(n,n-k)$, and so we only need to consider $k \leq \fl[n]$. Some definitions of $GP(n,k)$ also require that $n \geq 3$. If we allow $n = 1, 2$, then $GP(1,k)$ corresponds to the path graph $P_1$, and $GP(2,k)$ to the cycle graph $C_4$, irrespective of the choice of $k$.

The earliest result in this subsection is due to Guy and Harary (1967) \cite{guyharary1967}, who showed that the crossing number of any M\"{o}bius ladder graph is 1. Since $GP(2k,k)$ is a subdivision of a M\"{o}bius ladder, we have the following result.

\begin{theorem}[Guy and Harary, 1967 \cite{guyharary1967}] For $k \geq 3$, the following holds:
\[ cr(GP(2k,k)) = 1. \]\end{theorem}

$GP(2k,k)$ is also known as the generalised Wagner graph, $V_{2k}$.

In 1981, Exoo et al.\ \cite{exooetal1981} considered $GP(n,k)$ for the special cases where $k = 1, 2$. They showed the former were always planar, and for $k = 2$ they gave the following theorem:


\begin{theorem}[Exoo et al.,\ 1981 \cite{exooetal1981}] The following holds:
\[ cr(GP(n,2)) = \begin{cases} 0, & \text{for even $n$,}\\ 0, & \text{for $n=3$,} \\ 2, & \text{for $n=5$,}\\ 3, & \text{for odd $n\geq7$.}\end{cases} \]\end{theorem}

In 1986, Fiorini \cite{fiorini1986} considered the crossing number of $GP(n,k)$ for the special case where $k = 3$, and settled the cases when $n = 0 \hspace*{-0.175cm}\mod 3$ and $n = 2 \hspace*{-0.175cm}\mod 3$. Fiorini also claimed to have solved the special case of $GP(10,3)$, but the proof was later refuted by McQuillan and Richter \cite{mcquillanrichter1992} in 1992, who also provided a nicer proof of the special case of $GP(8,3)$. They conjectured the crossing number for the case when $n = 1 \hspace*{-0.175cm}\mod 3$, which was later proved in 2002 by Richter and Salazar \cite{richtersalazar2002}, who also corrected some errors in Fiorini's proofs, to settle the case where $k = 3$. The results are summarised in the following theorem:


\begin{theorem}[Fiorini, (1986) \cite{fiorini1986}, Richter and Salazar, (2002) \cite{richtersalazar2002}]For $h \in \{0,1,2\}$, the following holds:
\[ cr(GP(3m+h,3)) = \begin{cases}2, & \text{for $h=0$, $m=3$},\\ m, & \text{for $h=0$, $m\geq4$},\\3, & \text{for $h=1$, $m=2$}, \\  m + 3, & \text{for $h=1$, $m\geq 3$},\\
m + 2, & \text{for $h=2$, $m \geq2$}.\end{cases} \]
\end{theorem}

For $k = 4$, Fiorini \cite{fiorini1986} also considered the crossing number of the special case of $GP(4n,4)$, while in 1997 Sara\v{z}in \cite{sarazin1997} settled the special case of $GP(10,4)$.

\begin{theorem}[Fiorini, 1986 \cite{fiorini1986}] For $n \geq 4$, the following holds:
\[ cr(GP(4n,4)) = 2n. \]\end{theorem}

\begin{theorem}[Sara\v{z}in, 1997 \cite{sarazin1997}] $ cr(GP(10,4)) = 4. $ \end{theorem}

In 2009, Lin et al.\ \cite{linetal2009} then settled every remaining graph for $n \leq 16$ (see the upcoming \hyperref[tab-GP]{Table \ref{tab-GP}}), and gave some conjectures for two remaining cases:

\begin{conjecture}[Lin et al.,\ 2009 \cite{linetal2009}] For $k \geq 3$, the following hold:
\[ cr(GP(4k+2,2k)) = 2k+1, \]
\[ cr(GP(4k+2,4)) = 2k+2. \]\end{conjecture}

%

The known crossing numbers of $GP(n,k)$ for $n \leq 17$ are summarised in \hyperref[tab-GP]{Table \ref{tab-GP}} and \hyperref[app-gpref]{Appendix \ref{app-gpref}} lists the publications where these results were first proved.

\setlength\tabcolsep{1 pt}
\begin{table}[h!]\centering
\caption{Crossing numbers for $GP(n,k)$, $5 \leq n \leq 17$. For $n \leq4$, $GP(n,k)$ is planar. \hyperref[app-gpref]{Appendix \ref{app-gpref}} lists the publications where these results were first proved. The numbers in this table (other than for $n = 17$) were first collated by Lin et al.\ \cite{linetal2009}.}\smallskip\smallskip\label{tab-GP}
\begin{tabular}{|c|c|c|c|c|c|c|c|c|c|c|c|c|c|}\hline
\multicolumn{1}{|r|}{\;\;\; \bf $n$ \;} & \multirow{2}{*}{\: \bf 5 \:} & \multirow{2}{*}{\: \bf 6 \:} & \multirow{2}{*}{\: \bf 7 \:} & \multirow{2}{*}{\: \bf 8 \:} & \multirow{2}{*}{\: \bf 9 \:} & \multirow{2}{*}{\: \bf 10 \:} & \multirow{2}{*}{\: \bf 11 \:} & \multirow{2}{*}{\: \bf 12 \:} & \multirow{2}{*}{\: \bf 13 \:} & \multirow{2}{*}{\: \bf 14 \:} & \multirow{2}{*}{\: \bf 15 \:} & \multirow{2}{*}{\: \bf 16 \:} & \multirow{2}{*}{\: \bf 17 \:} \\
\multicolumn{1}{|c|}{\bf $k$\;\;} & & & & & & & & & & & & & \\
\hline \rule{0pt}{2.5ex} {\bf 1}   &  0  &  0  &  0  &  0  &  0  &  0  &  0  &  0  &  0  &  0  &  0  &  0  &  0\\
\hline \rule{0pt}{2.5ex} {\bf 2}   &  2  &  0  &  3  &  0  &  3  &  0  &  3  &  0  &  3  &  0  &  3  &  0  &  3\\
\hline \rule{0pt}{2.5ex} {\bf 3}   &    &  1  &  3  &  4  &  2  &  6  &  5  &  4  &  7  &  6  &  5  &  8  &  7\\
\hline \rule{0pt}{2.5ex} {\bf 4}   &    &    &    &  1  &  3  &  4  &  5  &  4  &  7  &  8  &  10  &  8  &  \\
\hline \rule{0pt}{2.5ex} {\bf 5}   &    &    &    &    &    &  1  &  3  &  8  &  9  &  6  &  5  &  8  &  \\
\hline \rule{0pt}{2.5ex} {\bf 6}   &    &    &    &    &    &    &    &  1  &  3  &  7  &  10  &  12  &  7\\
\hline \rule{0pt}{2.5ex} {\bf 7}   &    &    &    &    &    &    &    &    &    &  1  &  3  &  9  &  \\
\hline \rule{0pt}{2.5ex} {\bf 8}   &    &    &    &    &    &    &    &    &    &    &    &  1  &  3\\
\hline \end{tabular}
\end{table}
\setlength\tabcolsep{6 pt}

There are also some results for arbitrarily large $k$. In 2005, Ma et al.\ \cite{maetal2005} settled the case for $GP(2k+1,k)$ for $k \geq 3$:

\begin{theorem}[Ma et al.,\ 2005 \cite{maetal2005}]For $k \geq 3$, the following holds:
\[ cr(GP(2k+1,k)) = 3. \]\end{theorem}

In 2003, Fiorini and Gauci \cite{fiorinigauci2003} settled the case for $GP(3k,k)$ for $k \geq 4$:

\begin{theorem}[Fiorini and Gauci, 2003 \cite{fiorinigauci2003}]For $k \geq 4$, the following holds:
\[ cr(GP(3k,k)) = k. \]\end{theorem}

In 2019, Gauci and Xuereb settle the cases for $GP(3k-1,k$) and $GP(3k+1,k)$ for $k \geq 3$.

\begin{theorem}[Gauci and Xuereb, 2019 \cite{gaucixuereb2019}]For $k \geq 3$, the following hold:
\[ cr(GP(3k-1,k)) = k+1,\]
\[ cr(GP(3k+1,k)) = k+3.\]\end{theorem}

%
%

For general $GP(n,k)$ with $k \geq 5$, there are also lower and upper bounds known. In particular, in 2004, Pinontoan and Richter \cite{pinontoanrichter2004} determined bounds for $k \geq 6$ and $n \geq 2k+1$:

\begin{theorem}[Pinontoan and Richter, 2004 \cite{pinontoanrichter2004}]For $k \geq 6$, $n \geq 2k+1$, there exists a nonnegative constant $c_k$ such that:
\[ \frac{nk}{3} - c_k \leq cr(GP(n,k)) \leq 2n - \frac{4n}{k} + c_k. \]\end{theorem}

Also, in 2005, Salazar \cite{salazar2005} derived bounds via a different approach for $n \geq k$ and $k \geq 5$:

\begin{theorem}[Salazar, 2005 \cite{salazar2005}]For $n \geq k \geq 5$, the following bounds hold:
\[  \frac{2}{5}\left[\left(1 - \frac{4}{k}\right)(n-k^4)\right] + (4k^2 + 1 - k^3) \leq cr(GP(n,k)) \leq \left(2 - \frac{2}{k}\right)n + \frac{k^2 + k + 2}{2}. \]\end{theorem}

\subsection{Path powers}\label{sec-pathpowers}

Consider the path graph $P_n$ on $n+1$ vertices. The graph $P_n^k$, called the $k$-power of the graph $P_n$, is a graph on the same vertex set as $P_n$. An edge $\{a,b\}$ exists in $P_n^k$ if and only if the distance between $a$ and $b$ on $P_n$ is at most $k$. Two examples of path powers, $P_6^2$ and $P_6^3$, are displayed in \hyperref[fig-pathpowers]{Figure \ref{fig-pathpowers}}.

\begin{figure}[h!]\begin{center}\begin{tikzpicture}[largegraph,scale=1]
\foreach \n in {0,...,6}{
  \node (\n) at (\n,0) {};
}

\foreach \n in {0,...,5}{
  \pgfmathtruncatemacro{\nextn}{\n + 1}
  \draw (\n) -- (\nextn);
}

\foreach \n in {0,...,4}{
  \pgfmathtruncatemacro{\nextnn}{\n + 2}
  \draw[bend left=90,looseness=0.6] (\n) to (\nextnn);
}
\end{tikzpicture}$\;\;\;\;\;\;\;\;\;\;\;\;\;\;\;\;\;$
\begin{tikzpicture}[largegraph,scale=1]
\foreach \n in {0,...,6}{
  \node (\n) at (\n,0) {};
}

\foreach \n in {0,...,5}{
  \pgfmathtruncatemacro{\nextn}{\n + 1}
  \draw (\n) -- (\nextn);
}

\foreach \n in {0,...,4}{
  \pgfmathtruncatemacro{\nextnn}{\n + 2}
  \draw[bend left=90,looseness=0.6] (\n) to (\nextnn);
}
\foreach \n in {0,...,3}{
  \pgfmathtruncatemacro{\nextnnn}{\n + 3}
  \draw[bend left=90,looseness=0.6] (\n) to (\nextnnn);
}
\end{tikzpicture}\caption{The path powers $P_6^2$ and $P_6^3$.\label{fig-pathpowers}}\end{center}\end{figure}
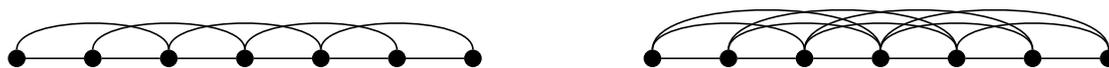

The crossing number of $P_n^k$ is known for some values of $k$. First, if $k \leq 3$, then $P_n^k$ is planar. If $k = n$, then $P_n^k$ is isomorphic to $K_n$; see \hyperref[sec-complete]{Section \ref{sec-complete}}. If $k = n-1$, then $P_n^k$ is isomorphic to $K_n$ minus an edge; see \hyperref[sec-completeminus]{Section \ref{sec-completeminus}}.

At this point the only other case which has been fully settled is $k = 4$, which was considered in 1993 by Harary and Kainen \cite{hararykainen1993}:

\begin{theorem}[Harary and Kainen,\ 1993 \cite{hararykainen1993}]For $n \geq 4$, the following holds:
\[ cr(P_n^4) = n-3. \]\end{theorem}

In 1999, Harary et al.\ \cite{hararyetal1999} extended this result by providing lower and upper bounds for $P_n^5$.

\begin{theorem}[Harary et al.,\ 1999 \cite{hararyetal1999}]For $n \geq 6$, the following bounds hold:
\[ 2n-9 \leq cr(P_n^5) \leq 4n - 21. \]\label{thm-Pn5}\end{theorem}

In 2009, Zheng et al.\ \cite{zhengetal2009} gave some upper bounds for higher powers, and conjectured that they would coincide with the exact crossing number. These results are listed in \hyperref[thm-pathpowers-table]{Table \ref{thm-pathpowers-table}}.

{\renewcommand{\arraystretch}{1.5}
\begin{table*}[htbp]\begin{center}\caption{Results for paths of higher powers, from \cite{zhengetal2009}.}\label{thm-pathpowers-table}\smallskip\smallskip$\begin{array}{|l|l|l|l|}
\hline cr(P_6^5)=3 & cr(P_7^5) = 6 & cr(P_8^5) = 9& cr(P_n^5) \leq 4n-23, \text{ for } n \geq 9\\
\hline cr(P_7^6) = 9 & cr(P_8^6) = 15 & cr(P_9^6) \leq 22& cr(P_n^6) \leq 8n-51, \text{ for } n \geq 10\\
\hline cr(P_8^7)= 18  & cr(P_9^7)\leq 30 & cr(P_{10}^7) \leq 42& cr(P_{11}^7) \leq 57\\
\hline & & & cr(P_n^7) \leq 15n-109, \text{ for } n \geq 12\\
\hline\end{array}$\end{center}\end{table*}
{\renewcommand{\arraystretch}{1}


\begin{conjecture}[Zheng et al.,\ 2009 \cite{zhengetal2009}]All upper bounds in \hyperref[thm-pathpowers-table]{Table~\ref{thm-pathpowers-table}} hold with equality.\label{conj-pathpowers}\end{conjecture}

\subsection{Kn\"{o}del graphs}

The Kn\"{o}del graph $\mathcal{W}_{\Delta,n}$, for even $n \geq 2$ and $1 \leq \Delta \leq \lfloor\log_2{n}\rfloor$ is the graph on $n$ vertices which is defined as follows. The vertices are labelled $(i,j)$ where $i = 1, 2$ and $0 \leq j \leq \frac{n}{2} - 1$. For every $j$, there is an edge between $(1,j)$ and every vertex $(2,(j + 2^p - 1)\mod (n/2))$ for $p = 0, \hdots, \Delta - 1$. An example of the Kn\"{o}del graph $\mathcal{W}_{3,10}$ is displayed in \hyperref[fig-knodel]{Figure \ref{fig-knodel}}.

\begin{figure}[h!]\begin{center}
\begin{tikzpicture}[largegraph]
\foreach \n in {0,...,4}{
  \node (a\n) at (\n,2) {};
  \node (b\n) at (\n,0) {};
}
\foreach \n in {0,...,4} \foreach \k in {0,1,2}{
  \pgfmathtruncatemacro{\m}{mod(\n+2^\k-1,5)}
  \draw (a\n) -- (b\m);
}
\end{tikzpicture}
\caption{The Kn\"{o}del graph $\mathcal{W}_{3,10}$.\label{fig-knodel}}\end{center}\end{figure}
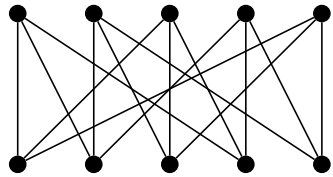

Zheng et al.\ \cite{zhengetal2008_2} considered the special case of $\mathcal{W}_{3,n}$, for even $n \geq 8$, and determined the crossing number in all cases:

\begin{theorem}[Zheng et al.,\ 2008 \cite{zhengetal2008_2}]The following holds:
\[ cr(\mathcal{W}_{3,n}) = \begin{cases}0, & \mbox{ for } n=8, \\ 1, & \mbox{ for } n=10, \\ \left\lfloor\frac{n}{6}\right\rfloor + \frac{(n\hspace*{-0.25cm}\mod 6)}{2}, & \mbox{ for even } n\geq12. \end{cases}\]\end{theorem}

\subsection{Flower Snarks}

Snarks are simple, connected, bridgeless cubic graphs with edge chromatic number equal to 4. Most definitions also demand that Snarks have girth at least five. The Flower Snarks are an infinite family of Snarks discovered by Isaacs \cite{isaacs1975}. The following construction is valid for any $n \geq 3$, but only odd $n \geq 5$ yield Flower Snarks. The results encompass both constructions from odd and even $n$ and collectively we shall denote the graphs as $I_n$.  In particular, the graph $I_n$ on $4n$ vertices is produced by first taking $n$ copies of $K_{1,3}$. Suppose that in copy $i$ of $K_{1,3}$ the three vertices with degree 1 are labelled $a_i$, $b_i$ and $c_i$. The copies are joined together via a cycle $a_1, a_2, \hdots, a_n, a_1$, and a cycle $b_1, b_2, \hdots, b_n, c_1, c_2, \hdots, c_n, b_1$.  Two examples of Flower Snarks, $I_5$ and $I_7$, are displayed in \hyperref[fig-flower]{Figure \ref{fig-flower}}.

\begin{figure}[h!]\begin{center}
\begin{tikzpicture}[largegraph,rotate=216]
\foreach \n in {0,...,4}{
  \pgfmathtruncatemacro{\ang}{90+\n*360/5}
  \begin{scope}[rotate=\ang]
    \node (\n) at (2,0) {};
    \path (\n) -- ++(180:1) node (b\n) {};
    \path (\n) -- ++(60:1) node (c\n) {};
    \path (\n) -- ++(300:1) node (d\n) {};
    \draw (\n) -- (b\n);
    \draw (\n) -- (c\n);
    \draw (\n) -- (d\n);
  \end{scope}
}
\foreach \n in {0,...,4}{
  \pgfmathtruncatemacro{\nextn}{mod(\n+1,5)}
  \draw (b\n) -- (b\nextn);
}
\foreach \n in {0,...,3}{
  \pgfmathtruncatemacro{\nextn}{mod(\n+1,5)}
  \draw[bend right=70] (c\n) to (c\nextn);
}
\draw (c4) -- (d0);
\foreach \n in {0,...,3}{
  \pgfmathtruncatemacro{\nextn}{mod(\n+1,5)}
  \draw[bend right=70] (d\n) to (d\nextn);
}
\draw[bend right=80] (d4) to (c0);
\end{tikzpicture}\;\;\;\;\;\;\;\;\;\;\;\;\;\;\;
\begin{tikzpicture}[largegraph,rotate=206]
\foreach \n in {0,...,6}{
  \pgfmathtruncatemacro{\ang}{90+\n*360/7}
  \begin{scope}[rotate=\ang]
    \node (\n) at (2,0) {};
    \path (\n) -- ++(180:.8) node (b\n) {};
    \path (\n) -- ++(60:.8) node (c\n) {};
    \path (\n) -- ++(300:.8) node (d\n) {};
    \draw (\n) -- (b\n);
    \draw (\n) -- (c\n);
    \draw (\n) -- (d\n);
  \end{scope}
}
\foreach \n in {0,...,6}{
  \pgfmathtruncatemacro{\nextn}{mod(\n+1,7)}
  \draw (b\n) -- (b\nextn);
}
\foreach \n in {0,...,5}{
  \pgfmathtruncatemacro{\nextn}{mod(\n+1,7)}
  \draw[bend right=70] (c\n) to (c\nextn);
}
\draw (c6) -- (d0);
\foreach \n in {0,...,5}{
  \pgfmathtruncatemacro{\nextn}{mod(\n+1,7)}
  \draw[bend right=70] (d\n) to (d\nextn);
}
\draw[bend right=80] (d6) to (c0);
\end{tikzpicture}
\caption{The Flower Snarks $I_5$ and $I_7$.\label{fig-flower}}\end{center}\end{figure}
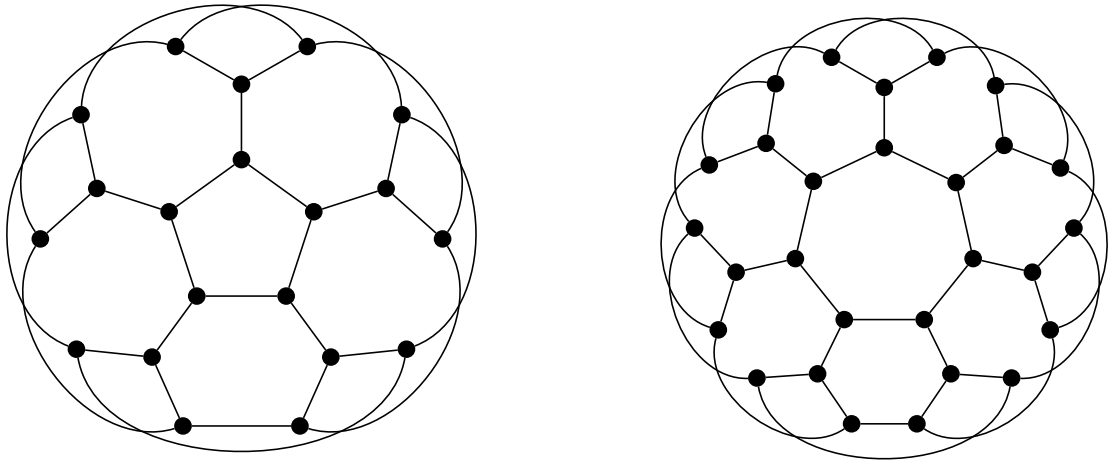

Zheng et al.\ \cite{zhengetal2008_3} determined the crossing numbers of the graphs $I_n$ (for both odd and even $n$):

\begin{theorem}[Zheng et al.,\ 2008 \cite{zhengetal2008_3}] For $n\geq 3$, the following holds:
\[ cr(I_n) = \begin{cases} n-1, & \mbox{ for } 3\leq n \leq 5, \\ n, & \mbox{ for } n \geq 6. \end{cases} \]\end{theorem}

\subsection{Hexagonal graph $H_{3,n}$}

The Hexagonal graph $H_{3,n}$ is defined as follows for $n \geq 2$. Take the union of three cycles of length $2n$, and label their vertices $a_i$, $b_i$ and $c_i$ respectively, for $i = 1, \hdots, 2n$. Then add the edges $\{a_{2i-1},b_{2i-1}\}, \{b_{2i},c_{2i}\}, \{a_{2i},c_{2i-1}\}$, for $i = 1, \hdots, n$. An example of the Hexagonal graph $H_{3,3}$ is displayed in \hyperref[fig-hex]{Figure \ref{fig-hex}}.

%
%

\begin{figure}[h!]\begin{center}
\begin{tikzpicture}[largegraph]
\foreach \n in {0,...,5}{
  \node (a\n) at (\n*1.5,1) {};
  \node (b\n) at (\n*1.5,0) {};
  \node (c\n) at (\n*1.5,-1) {};
}
\foreach \n in {0,...,4}{
  \pgfmathtruncatemacro{\nextn}{\n+1}
  \draw (a\n) -- (a\nextn);
  \draw (b\n) -- (b\nextn);
  \draw (c\n) -- (c\nextn);
}
\draw[bend left=90,looseness=0.2] (a0) to (a5);
\draw[bend right=90,looseness=0.2] (c0) to (c5);
\draw[bend left=135,looseness=1.3] (b0) to (b5);
\foreach \n in {0,...,2}{
  \pgfmathtruncatemacro{\aindex}{2*\n}
  \pgfmathtruncatemacro{\bindex}{2*\n}
  \draw (a\aindex) -- (b\bindex);
  \pgfmathtruncatemacro{\bindex}{2*\n+1}
  \pgfmathtruncatemacro{\cindex}{2*\n+1}
  \draw (b\bindex) -- (c\cindex);
  \pgfmathtruncatemacro{\cindex}{2*\n}
  \pgfmathtruncatemacro{\aindex}{2*\n+1}
  \draw (c\cindex) -- (a\aindex);
}
\end{tikzpicture}
\caption{The Hexagonal graph $H_{3,3}$, in an optimal drawing.\label{fig-hex}}\end{center}\end{figure}
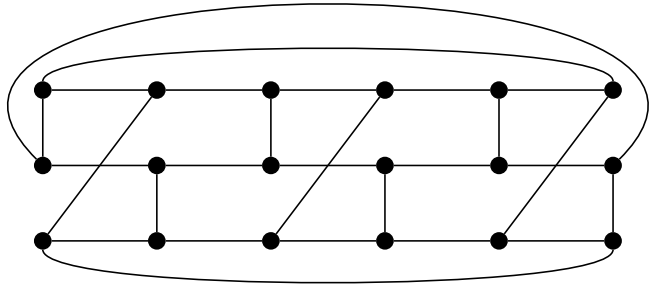

In 2019, Wang et al.\ \cite{wangetal2019} considered $H_{3,n}$ and determined its crossing number for all $n \geq 2$:

\begin{theorem}[Wang et al.\ (2019) \cite{wangetal2019}] For $n \geq 2$, the following holds:
\[ cr(H_{3,n}) = n. \]\end{theorem}

\section{Cartesian products of graphs}\label{sec-cartesian}

The Cartesian product of two graphs $G$ and $H$ is written as $G \Box H$. The result is a graph with vertex set $V(G) \times V(H)$, and edges between vertices $(u,u')$ and $(v,v')$ if and only if either $u = v$ and $(u',v') \in E(H)$, or $u' = v'$ and $(u,v) \in E(G)$. An example of the Cartesian product of two paths, $P_3 \Box P_4$, is displayed in \hyperref[fig-p3boxp4]{Figure \ref{fig-p3boxp4}}.

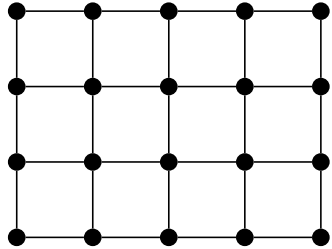
\begin{figure}[h!]\begin{center}
\begin{tikzpicture}[largegraph,scale=0.5]
\foreach \n in {0,...,4}{
  \node (a\n) at (2*\n,6) {};
  \node (b\n) at (2*\n,4) {};
  \node (c\n) at (2*\n,2) {};
  \node (d\n) at (2*\n,0) {};
}
\foreach \n in {0,...,3}{
  \pgfmathtruncatemacro{\m}{\n + 1}
  \draw (a\n) -- (a\m);
  \draw (b\n) -- (b\m);
  \draw (c\n) -- (c\m);
  \draw (d\n) -- (d\m);
}
\foreach \n in {0,...,4}{
  \draw (a\n) -- (b\n);
  \draw (b\n) -- (c\n);
  \draw (c\n) -- (d\n);
}

\end{tikzpicture}\caption{The Cartesian product $P_3 \Box P_4$.\label{fig-p3boxp4}}\end{center}\end{figure}

The Cartesian product was the first graph product for which the crossing number was investigated in detail, beginning in the 1970s and continuing to this day. Investigations typically fall into two categories; either the Cartesian product of two families of graphs of unbounded size is considered, or the Cartesian product of specific graphs with a family of graphs of unbounded size. In the latter case, researchers all over the world have dedicated themselves to finding the crossing numbers of Cartesian products involving the various connected graphs of small orders, each often involving an ad hoc proof or technique unique to that graph.

If a graph $G$ is disconnected, then $G \Box H$ is equivalent to taking the union of the Cartesian products of each component of $G$ with $H$. Hence, the crossing number of $G \Box H$ is simply equal to the sum of the Cartesian products of $H$ with each of its connected components. To that end, in what follows, all results are for connected graphs.

\subsection{Cartesian products of cycles, paths and stars}

A number of results have been determined for the Cartesian product of families of graphs. In particular, much attention has been paid to Cartesian products involving cycles, paths, and stars. The star graph $S_n$ is simply the complete bipartite graph $K_{1,n}$. Since $K_{1,1} = P_1$ and $K_{1,2} = P_2$, we only consider $S_n$ for $n \geq 3$. It is easy to check that the Cartesian product of a cycle with a path, or a path with a path, will result in a planar graph. It can also easily be checked that the Cartesian product of two stars $S_n \Box S_m$ is isomorphic to a subdivision of the complete tripartite graph $K_{1,m,n}$, which is discussed in \hyperref[sec-complete_tripartite]{Section \ref{sec-complete_tripartite}}. We now consider the three remaining cases.

\subsubsection{Two cycles}

The first publication to consider the Cartesian product of two cycles was by Harary et al.\ \cite{hararykainen1973} in 1973:

\begin{conjecture}[Harary et al.,\ 1973 \cite{hararykainen1973}] For $n \geq m \geq 3$, the following holds:
\[ cr(C_m \Box C_n) = (m-2)n. \]\label{conj-cyccyc}\end{conjecture}

\hyperref[conj-cyccyc]{Conjecture~\ref{conj-cyccyc}} has been verified for $m \leq 7$. Most of the proofs first rely on finding the crossing number for the special case when $n = m$ and then extending it for $n > m$. The first such result was by Ringeisen and Beineke \cite{ringeisenbeineke1978} who verified \hyperref[conj-cyccyc]{Conjecture~\ref{conj-cyccyc}} for $m = 3$ in 1978.

In 1970, Eggleton and Guy \cite{eggletonguy1970} proved that $cr(C_4 \Box C_4) = 8$ but the paper contained a (separate) error and was never published. In 1980, Beineke and Ringeisen \cite{beinekeringeisen1980} verified \hyperref[conj-cyccyc]{Conjecture~\ref{conj-cyccyc}} for $m = 4$ but relied on Eggleton and Guy's result which had still not appeared in the literature. A proof for $cr(C_4 \Box C_4) = 8$ was finally published in 1995 by Dean and Richter \cite{deanrichter1995}.

In 1995, Richter and Thomassen \cite{richterthomassen1995} verified \hyperref[conj-cyccyc]{Conjecture~\ref{conj-cyccyc}} for $C_5 \Box C_5$, and the general case was subsequently verified for $m = 5$ by Kle\^{s}\^{c} et al.\ \cite{klescetal1996} in 1996.

In 1996, Anderson et al.\ \cite{andersonetal1996} verified \hyperref[conj-cyccyc]{Conjecture~\ref{conj-cyccyc}} for $C_6 \Box C_6$, and the general case was subsequently verified for $m = 6$ by Richter and Salazar \cite{richtersalazar2001} in 2001, although the bulk of the work for doing so was first detailed in 1997 in Salazar's PhD thesis \cite{salazar1997}.

In 1996, Anderson et al.\ \cite{andersonetal1996_2} verified \hyperref[conj-cyccyc]{Conjecture~\ref{conj-cyccyc}} for $C_7 \Box C_7$, and the general case was subsequently verified for $m = 7$ by Adamsson and Richter \cite{adamssonrichter2004} in 2004. Adamsson and Richter also claimed (providing only a sketch of the proof) that the conjecture can be shown to hold for $m = 8$ if it can first be shown that $cr(C_8 \Box C_8) = 48$.

In 2004, Glebsky and Salazar \cite{glebskysalazar2004} provided a breakthrough by showing that for each $m$, the conjecture must be true for all but a finite number of cases.

\begin{theorem}[Glebsky and Salazar, 2004 \cite{glebskysalazar2004}] For $n \geq m(m+1)$ and $m \geq 3$, the following holds:
\[ cr(C_m \Box C_n) = (m-2)n. \]\end{theorem}

Some lower bounds have also been developed for the Cartesian product of two cycles. In 1973, Harary et al.\ \cite{hararykainen1973} showed that $cr(C_m \Box C_n) \geq m$ and asked if this could be improved. In 1995, Shahrokhi et al.\ \cite{shahrokhietal1995} showed that $cr(C_m \Box C_n) \geq \frac{mn}{90}$, and that $cr(C_m \Box C_n) \geq \frac{mn}{6}$ if $n = m$ or $n = m+1$. In 1998, Shahrokhi et al.\ \cite{shahrokhietal1998} improved their result further to show that $cr(C_m \Box C_n) \geq \frac{mn}{9}$, and that $cr(C_m \Box C_n) \geq \frac{3mn}{5}$ if $n \leq 5(m-1)/4$. In 2000, the lower bound was improved again to $cr(C_m \Box C_n) \geq \frac{(m-2)n}{3}$ by Salazar \cite{salazar2000}. Finally, this was further considered in 2004 by Salazar and Ugalde \cite{salazarugalde2004}, who gave what is currently the best asymptotic lower bound:

\begin{theorem}[Salazar and Ugalde, 2004 \cite{salazarugalde2004}] For every $\epsilon > 0$, there is an $N_\epsilon$ with the following property. For all $n \geq m \geq N_\epsilon$, the following holds:
\[ cr(C_m \Box C_n) \geq (0.8 - \epsilon)mn. \]\end{theorem}

\subsubsection{Paths and stars}

This case has now been completely settled. The first results were due to Jendrol' and \v{S}cerbov\'{a} \cite{jendrolscerbova1982} in 1982, where they determined an upper bound for the crossing number of $S_m \Box P_n$, and conjectured that it held with equality. Their conjecture would ultimately prove to be correct. They also verified their conjecture for $m = 3$, as well as the special case of $m = 4, n = 2$. In 1991, Kle\v{s}\v{c} \cite{klesc1991} verified the conjecture for $m = 4$. Finally, in 2007, Bokal verified the conjecture for all cases:

\begin{theorem}[Bokal, 2007 \cite{bokal2007}] For $m \geq 3$ and $n \geq 1$, the following holds:
\[ cr(S_m \Box P_n) = (n-1) \fl[m] \fl[m-1]. \]\end{theorem}

\subsubsection{Cycles and stars}\label{sec-cyclestar}

Results for $S_m \Box C_n$ are only known for $\min(m,n) \leq 4$. For arbitrarily large cycles, the crossing numbers of $S_3 \Box C_n$ and $S_4 \Box C_n$ are due to Jendrol' and \v{S}cerbov\'{a} (1982) \cite{jendrolscerbova1982} and Kle\v{s}\v{c} (1991) \cite{klesc1991} respectively.

\begin{theorem}[Jendrol' and \v{S}cerbov\'{a}, 1982 \cite{jendrolscerbova1982}] The following holds:
\[ cr(S_3 \Box C_n) = \begin{cases} 1, & \mbox{ for } n=3, \\ 2, & \mbox{ for } n=4, \\4, & \mbox{ for } n=5, \\n, & \mbox{ for }n \geq 6. \end{cases}\]\end{theorem}

\begin{theorem}[Kle\v{s}\v{c}, 1991 \cite{klesc1991}] The following holds:
\[ cr(S_4 \Box C_n) = \begin{cases} 2, & \mbox{ for } n=3, \\ 4, & \mbox{ for } n=4, \\8, & \mbox{ for } n=5, \\ 2n, & \mbox{ for }n \geq 6. \end{cases}\]\end{theorem}

For arbitrarily large stars, the crossing numbers of $S_m \Box C_3$ and $S_m \Box C_4$ are known. The former does not appear to have been explicitly stated anywhere, but it is trivial to exhibit drawings of $S_m \Box C_3$ with $\fl[m]\fl[m-1]$ crossings, which is also a lower bound because it contains $S_m \Box P_2$ as a subgraph. The latter case was settled by Kle\v{s}\v{c} (1994) \cite{klesc1994}.

\begin{theorem}[Kle\v{s}\v{c}, 1994 \cite{klesc1994}] For $m \geq 1$, the following holds:
\[cr(S_m \Box C_4) = 2\left\lfloor\frac{(m-1)^2}{4}\right\rfloor .\]\end{theorem}

\subsection{Cartesian products of paths with other graphs}

\subsubsection{Paths and cycles with one or two extra edges}

In 2007, Yuan and Huang \cite{yuanhuang2007} $\Asterisk$ considered the Cartesian product of a path with a graph consisting of a cycle with one or two extra edges added. Consider the cycle graph $C_m$ on $m \geq 5$ vertices. From this, construct a new graph $H$, by adding a chord between two vertices of distance two. Similarly, consider a new graph $B$, constructed by taking $C_m$ and adding in two chords, whose endpoints have distance at least two.

\begin{theorem}[Yuan and Huang, 2007 \cite{yuanhuang2007} $\Asterisk$]For $n \geq 1$, the following hold:
\[ cr(P_n \Box H) = n-1, \: \Asterisk \]
\[ cr(P_n \Box B) = 2n-2. \: \Asterisk \]\end{theorem}

%
%

\subsubsection{Paths and double cones}

The join product $G + H$ (considered further in \hyperref[sec-join]{Section \ref{sec-join}}) is equal to the union of $G$ and $H$, plus edges linking every vertex of $G$ to every vertex of $H$. A generalisation of the wheel graph is the cone graph $W_{l,n}$, which is equal to the join product $C_m + D_l$, where $D_l$ is the discrete graph on $l$ isolated vertices with zero edges. Then, $W_{1,n}$ is simply the wheel graph. In this context, $W_{2,n}$ is referred to as the double cone. An example of the double cone $W_{2,6}$ is displayed in \hyperref[fig-wheel]{Figure \ref{fig-wheel}}. In 2011, Zheng et al.\ \cite{zhengetal2011} considered the Cartesian product of double cones with paths:

\begin{theorem}[Zheng et al.,\ 2011 \cite{zhengetal2011}]For $m \geq 3$ and $n \geq 1$, the following holds:
\[ cr(W_{2,m} \Box P_n) = 2n\fl[m]\fl[m-1] + 2n. \]\end{theorem}

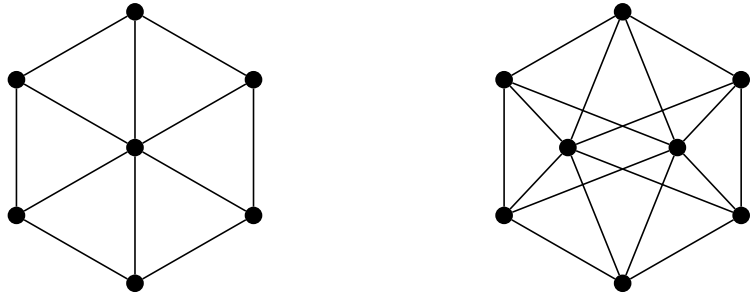
\begin{figure}[h!]\begin{center}\begin{tikzpicture}[largegraph,scale=0.8]
\foreach \n in {0,...,5}{
  \node (\n) at (90+\n/6*360:2.25) {};
}
\node (6) at (0,0) {};

\foreach \n in {0,...,5}{
  \pgfmathtruncatemacro{\nextn}{\n + 1}
  \draw (\n) -- (\nextn);
}
\foreach \n in {0,...,4}{
  \draw (\n) -- (6);
}
\draw (0) -- (5);
\end{tikzpicture}$\;\;\;\;\;\;\;\;\;\;\;\;\;\;\;\;\;\;\;\;\;\;\;\;\;$
\begin{tikzpicture}[largegraph,scale=0.8]
\foreach \n in {0,...,5}{
  \node (\n) at (90+\n/6*360:2.25) {};
}
\node (6) at (-0.9,0) {};
\node (7) at (0.9,0) {};

\foreach \n in {0,...,5}{
  \pgfmathtruncatemacro{\nextn}{\n + 1}
  \draw (\n) -- (\nextn);
  \draw (\n) -- (7);
}
\foreach \n in {0,...,4}{
  \draw (\n) -- (6);
}
\draw (0) -- (5);
\end{tikzpicture}\caption{The wheel graph $W_6$ and the double cone $W_{2,6}$.\label{fig-wheel}}\end{center}\end{figure}

\subsubsection{Paths and complete graphs}

Consider the Cartesian product of the complete graph and the path graph, $K_n \Box P_m$. It is clear that the graph is planar for $n = 3$, so we now consider $n > 3$.

The cases for $n = 4$ and $n = 5$ were settled by Kle\v{s}\v{c} \cite{klesc1994,klesc1999} in 1994 and 1999 respectively:

\begin{theorem}[Kle\v{s}\v{c}, 1994 \cite{klesc1994}]For $m \geq 1$, $cr(K_4 \Box P_m) = 2m$.\end{theorem}

\begin{theorem}[Kle\v{s}\v{c}, 1999 \cite{klesc1999}]For $m \geq 1$, $cr(K_5 \Box P_m) = 6m$.\end{theorem}

In 2007, Zheng et al.\ \cite{zhengetal2007} settled the case for $n = 6$, and gave an upper bound the general case which they conjectured would hold with equality:

\begin{theorem}[Zheng et al.,\ 2007 \cite{zhengetal2007}]For $m \geq 1$, $cr(K_6 \Box P_m) = 15m + 3$.\end{theorem}

\begin{theorem}[Zheng et al.,\ 2007 \cite{zhengetal2007}]For $n \geq 4$ and $m \geq 1$,
\[ cr(K_n \Box P_m) \leq \frac{1}{4}\fl[n+1]\fl[n-1]\fl[n-2]\left(m\fl[n+4] + \fl[n-4]\right). \]\label{thm-pathcomplete}\end{theorem}

\begin{conjecture}[Zheng et al.,\ 2007 \cite{zhengetal2007}] \hyperref[thm-pathcomplete]{Theorem \ref{thm-pathcomplete}} holds with equality.\label{conj-pathcomplete}\end{conjecture}

In 2014, Ouyang et al.\ \cite{ouyangetal2014_2} confirmed that \hyperref[conj-pathcomplete]{Conjecture~\ref{conj-pathcomplete}} holds for $n \leq 10$.

\subsubsection{Paths and circulant graphs}

Consider the path graph on $n+1$ vertices $P_n$, and the circulant graph $Ci_m(\{1,k\})$. The latter is sometimes referred to as $C(m,k)$. There has been some effort to find the crossing number of $P_n \Box Ci_m(1,k)$, with the cases for $m \leq 6$ being covered in \hyperref[sec-cart3]{Sections \ref{sec-cart3}--}\ref{sec-cart6}. Here, we include the known results for larger $m$.

\begin{theorem}[Yuan et al.,\ 2008 \cite{yuanhuangetal2008} $\Asterisk$] For $n \geq 1$, $cr(P_n \Box Ci_7(1,2)) = 8n. \: \Asterisk$\end{theorem}

\begin{theorem}[Yuan et al.,\ 2008 \cite{yuanetal2008}] For $n \geq 1$, $cr(P_n \Box Ci_8(1,2)) = 8n.$\end{theorem}

\begin{theorem}[Wang and Ma, 2017 \cite{wangma2017} $\Asterisk$] For $n \geq 1$, $cr(P_n \Box Ci_8(1,4)) = 9n - 1. \: \Asterisk$\end{theorem}

\begin{theorem}[Yuan et al.,\ 2013 \cite{yuanetal2013} ] For $n \geq 1$, $cr(P_n \Box Ci_9(1,2)) = 10n. $\end{theorem}

\begin{theorem}[Yuan et al.,\ 2009 \cite{yuanetal2009}] For $n \geq 1$, $cr(P_n \Box Ci_{10}(1,2)) = 10n.$\end{theorem}

Yuan et al.\ \cite{yuanetal2009} also claim, without providing a proof, that $cr(P_n \Box Ci_{11}(1,2)) = cr(P_n \Box Ci_{12}(1,2)) = 12n$. Later, a separate paper was published proving the latter result:

\begin{theorem}[Yuan and Huang, 2011 \cite{yuanhuang2011} ] For $n \geq 1$, $cr(P_n \Box Ci_{12}(1,2)) = 12n.$ \end{theorem}

\subsubsection{Paths and complete multipartite graphs}\label{sec-pathmulti}

There are a few results relating to the Cartesian products of paths and complete multipartite graphs. \hyperref[thm-Tang]{Theorem \ref{thm-Tang}} was also rediscovered a year later by Zhang et al.\ \cite{zhengetal2008}.

\begin{theorem}[Tang et al.,\ 2007 \cite{tangetal2007}]For $m \geq 2$ and $n \geq 1$,
\[ cr(P_n \Box K_{2,m}) = 2n\fl[m]\fl[m-1]. \]\label{thm-Tang}\end{theorem}

\begin{theorem}[Ouyang et al.,\ 2014 \cite{ouyangetal2014}]For $m \geq 2$  and $n \geq 1$,
\[ cr(P_n \Box K_{1,1,m}) = 2n\fl[m]\fl[m-1] + (n-1)\fl[m]. \]\end{theorem}

\subsection{Cartesian products of cycles with other graphs}

\subsubsection{Cycles and 2-powers of paths}

For a description of 2-powers of paths, see \hyperref[sec-pathpowers]{Section \ref{sec-pathpowers}}. In 2012, Kle\v{s}\v{c} and Kravecov\'{a} \cite{klesckravecova2012} considered the Cartesian product of $P^2_n$ and cycles. They conjectured the crossing number, and proved it for the special case where the cycle has three vertices:

\begin{conjecture}[Kle\v{s}\v{c} and Kravecov\'{a}, 2012 \cite{klesckravecova2012}]For $n \geq 2$ and $m \geq 3$,
\[ cr(P^2_n \Box C_m) = m(n-1). \]\end{conjecture}

\begin{theorem}[Kle\v{s}\v{c} and Kravecov\'{a}, 2012 \cite{klesckravecova2012}]For $n \geq 2$, \[cr(P^2_n \Box C_3) = 3n - 3.\]\end{theorem}

Later in 2012, Kravecov\'{a} and Petrillov\'{a} \cite{kravecova2012} also proved the result for $C_4$:

\begin{theorem}[Kravecov\'{a} and Petrillov\'{a}, 2012 \cite{kravecova2012}]For $n \geq 2$, \[cr(P^2_n \Box C_4) = 4n - 4.\]\end{theorem}

\subsection{Cartesian products of stars with other graphs}\label{sec-startree}

\subsubsection{Stars and trees}

In 2007, Bokal \cite{bokal2007_2} showed that the crossing number of $S_n \Box\, T$, where $T$ is a tree, can be written in terms of the crossing numbers of complete tripartite graphs, as follows:

\begin{theorem}[Bokal, 2007 \cite{bokal2007_2}]Consider a tree $T$, and define $d(v)$ to be the degree of vertex $v$ in $T$. Then,
\[ cr(S_n \Box\, T) = \sum_{v \in V(T)} cr(K_{1,d(v),n}). \]\end{theorem}

At the time Bokal published his result, $cr(K_{1,m,n})$ was only known for $m \leq 3$. In the ensuing years, it was determined for $m = 4, 5$ (see \hyperref[sec-complete_tripartite]{Section \ref{sec-complete_tripartite}}). Hence, the following corollary can be stated.

\begin{corollary}Consider a tree $T$ with maximum degree 5. Let $n_i$ be the number of vertices of degree $i$ contained in $T$, and set $a = n_2 + 2n_3 + 4n_4 + 6n_5$ and $b = n_3 + 2n_4 + 4n_5$. Then,
\[ cr(S_n \Box\, T) = a\fl[n]\fl[n-1] + b\fl[n]. \]\end{corollary}

\subsubsection{Stars and sunlet graphs}

The sunlet graph $\mathcal{S}_n$ is the graph on $2n$ vertices constructed by taking the cycle $C_n$ and adding a pendant edge to each vertex. Two examples of sunlet graphs, $\mathcal{S}_6$ and $\mathcal{S}_9$ are displayed in \hyperref[fig-sunlet]{Figure \ref{fig-sunlet69}}.

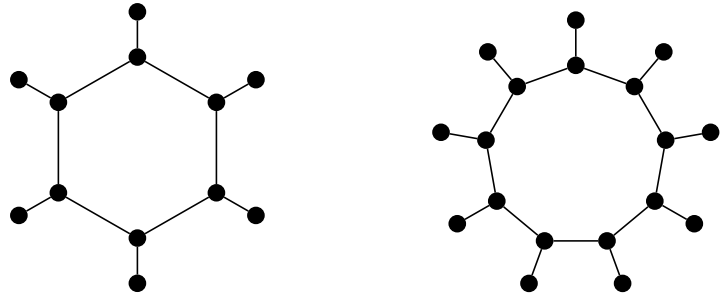
\begin{figure}[h!]\begin{center}\begin{tikzpicture}[largegraph,scale=0.8]
\foreach \n in {0,...,5}{
  \node (\n) at (90+\n/6*360:1.5) {};
  \node (o\n) at (90+\n/6*360:2.25) {};
  \draw (\n) -- (o\n);
}
\foreach \n in {0,...,5}{
  \pgfmathtruncatemacro{\nextn}{mod(\n+1,6)}
  \pgfmathtruncatemacro{\nextnextn}{mod(\n+2,6)}
  \draw (\n) -- (\nextn);
}
\end{tikzpicture}\;\;\;\;\;\;\;\;\;\;\;\;\;\;\;\;\;\;
\begin{tikzpicture}[largegraph,scale=0.8]
\foreach \n in {0,...,8}{
  \node (\n) at (90+\n/9*360:1.5) {};
  \node (o\n) at (90+\n/9*360:2.25) {};
  \draw (\n) -- (o\n);
}
\foreach \n in {0,...,8}{
  \pgfmathtruncatemacro{\nextn}{mod(\n+1,9)}
  \pgfmathtruncatemacro{\nextnextn}{mod(\n+2,9)}
  \draw (\n) -- (\nextn);
}
\end{tikzpicture}\caption{The sunlet graphs $\mathcal{S}_6$ and $\mathcal{S}_9$.\label{fig-sunlet69}}\end{center}\end{figure}

In 2019, Haythorpe and Newcombe \cite{haythorpenewcombe2019} considered the Cartesian product of a sunlet graph $\mathcal{S}_n$ and a star graph $S_m$. They were able to determine the crossing number for $m \leq 3$ and gave an upper bound for larger $m$ which they conjecture coincides with the crossing number.

\begin{theorem}[Haythorpe and Newcombe, 2019 \cite{haythorpenewcombe2019}]For $n \geq 3$ and $m \geq 1$, the following holds:
 \[cr(\mathcal{S}_n \Box\, S_m) = \begin{cases}0, & \mbox{ for } m=1, \\n, & \mbox{ for } m=2, \\ 3n, & \mbox{ for } m=3.\end{cases}\]
Additionally, for $m\geq 4$,
$cr(\mathcal{S}_n \Box\, S_m) \leq \frac{nm(m-1)}{2}.$\label{thm-sunlet}\end{theorem}

\begin{conjecture}[Haythorpe and Newcombe, 2019 \cite{haythorpenewcombe2019}]The upper bound in \hyperref[thm-sunlet]{Theorem~\ref{thm-sunlet}} holds with equality.\end{conjecture}

\subsection{Cartesian products of other graph families}

\subsubsection{Complete graphs and complete bipartite graphs}\label{sec-cartesian-complete}

Zheng et al.\ \cite{zhengetal2008} proved the following bounds related to Cartesian products involving complete graphs with cycles, and complete bipartite graphs with paths:

\begin{theorem}[Zheng et al.,\ 2008 \cite{zhengetal2008}] The following hold:

\begin{enumerate}\item $cr(K_m \Box C_n) \geq n \cdot cr(K_{m+2})$, for $n \geq 3$, $m \geq 5$.
\item $cr(K_m \Box C_n) \leq \frac{n}{4} \fl[m+2] \fl[m+1] \fl[m] \fl[m-1]$, for $m = 5,6,7$, $n \geq 3$, and for $m \geq 8$ with even $n \geq 4$.
\item $cr(K_{m,l} \Box P_n) \leq (n-1)(\fl[m+2] \fl[m+1] \fl[l+2] \fl[l+1] - ml) \\[0.25em]\hspace*{8em}+ 2(\fl[m+1] \fl[m] \fl[l+1]\fl[l]  - \fl[m] \fl[l] )$, for $m,l \geq 2$ and $n \geq 1$.\end{enumerate}\end{theorem}

They also showed that equality holds for item 2 for some small values of $m$:

\begin{theorem}[Zheng et al.,\ 2008 \cite{zhengetal2008}] For $n \geq 1, m \geq 2$,

$cr(K_m \Box C_n) = \left\{\begin{array}{ccl}9n & & m = 5, n \geq 3,\\
18n & & m = 6, n \geq 3,\\
36n & & m = 7, n \geq 3,\\
60n & & m = 8, n \geq 4 \mbox{ and even},\\
100n & & m = 9, n \geq 4 \mbox{ and even},\\
150n & & m = 10, n \geq 4 \mbox{ and even}.\end{array}\right.$
\end{theorem}


\subsubsection{Wheels and trees}

Consider the wheel graph $W_n$ on $n+1$ vertices. In 2017, Kle\v{s}\v{c} et al.\ \cite{klescetal2017} considered the Cartesian product of wheel graphs with trees of maximum degree no larger than 5:

\begin{theorem}[Kle\v{s}\v{c} et al.,\ 2017 \cite{klescetal2017}]Let $W_n$ be the wheel graph on $n+1$ vertices, and $T$ be a tree with maximum degree $\triangle(T) \leq 5$. Let $n_i$ be the number of vertices of degree $i$ in $T$, and set $a = n_1 + n_2 + 2n_3 + 2n_4 + 3n_5$ and $b = n_2 + 2n_3 + 4n_4 + 6n_5$ and $c = n_3 + 2n_4 + 4n_5$. Then, for $n \geq 3$,
\[ cr(W_n \Box\, T) = a + b\fl[n]\fl[n-1] + c\fl[n]. \]\label{thm-wheeltree}\end{theorem}

This result extended an earlier result of Cartesian products of wheels with subcubic trees by Bokal \cite{bokal2007_2} in 2007.

Wang and Huang \cite{wanghuangtoappear} gave a similar result, except for the situation when the tree can have any degree, but the wheel is restricted to maximum degree 5.

\begin{theorem}[Wang and Huang,\ (to appear) \cite{wanghuangtoappear}]Let $W_n$ be the wheel\\graph on $n+1$ vertices, for $n = 3, 4, 5$, and $T$ be a tree with maximum degree $\triangle(T)$. Let $n_i$ be the number of vertices of degree $i$ in $T$. Then,
\[ cr(W_n \Box T) = \begin{cases}\sum_{i=1}^{\triangle(T)} n_i(2\fl[i]\fl[i-1] + i) & \text{for $n = 3$,}\\
\sum_{i=1}^{\triangle(T)} n_i(4\fl[i]\fl[i-1] + i + \fl[i]) & \text{for $n = 4$,}\\
\sum_{i=1}^{\triangle(T)} n_i(6\fl[i]\fl[i-1] + i + 3\fl[i]) & \text{for $n = 5$.}
\end{cases} \]\end{theorem}

\subsection{Cartesian products with 3-vertex graphs}\label{sec-cart3}

There are only two connected non-isomorphic graphs on 3 vertices; the path $P_2$ and the cycle $C_3$. Although no author made a point of considering Cartesian products with 3-vertex graphs specifically, the results presented in \hyperref[tab-3vcart]{Table~\ref{tab-3vcart}} are all known from other works, as follows.

The Cartesian products of paths and paths, or paths and cycles, always produce planar graphs. The result for $P_2 \Box S_n$ is due to Bokal (2007) \cite{bokal2007}. The result for $C_3 \Box C_n$ is due to Ringeisen and Beineke (1978) \cite{ringeisenbeineke1978}. The result for $C_3 \Box S_n$ was discussed in \hyperref[sec-cyclestar]{Section \ref{sec-cyclestar}}.

\begin{table}[h!]\begin{center}\caption{Crossing numbers of Cartesian products of 3-vertex graphs with paths, cycles and stars. The results for $P_n$ are for $n \geq 1$, and the results for $C_n$ and $S_n$ are for $n \geq 3$.\label{tab-3vcart}}\smallskip\smallskip$\begin{array}{|c|c|c|c|c|}\hline
\rule{0pt}{2.3ex}& G & cr(G \Box P_n) & cr(G \Box C_n) & cr(G \Box S_n)\\
\hline \ctikz{\ngraph{3.3}} & P_2 & 0      & 0                & \fl[n]\fl[n-1] \\
\hline \ctikz{\ngraph{3.4}} & C_3 & 0      & n                & \fl[n]\fl[n-1] \\
\hline\end{array}$\end{center}\end{table}

\subsection{Cartesian products with 4-vertex graphs}\label{sec-cart4}

There are six connected non-isomorphic graphs on 4 vertices. We denote the graphs on 4 vertices as $G^4_i$, and note that there are eleven such graphs. However, only six of these are connected, which correspond to $i = 5, 7, 8, 9, 10, 11$ as displayed in \hyperref[tab-4vcart]{Table \ref{tab-4vcart}}. In 1980, Beineke and Ringeisen \cite{beinekeringeisen1980} considered the crossing numbers of the Cartesian product of each of them with the cycle graph on $n$ vertices. They were able to settle all cases except for $G^4_7 \Box C_n$, which was subsequently settled in 1982 by Jendrol' and \v{S}cerbov\'{a} \cite{jendrolscerbova1982}. The latter authors also settled $G^4_7 \Box P_n$. The result for $G^4_7 \Box C_n$ is only valid when $n \geq 6$. For smaller values of $n$, we have $cr(G^4_7 \Box C_3) = 1$, $cr(G^4_7 \Box C_4) = 2$ and $cr(G^4_7 \Box C_5) = 4$.

In 1994, Kle\v{s}\v{c} \cite{klesc1994} then derived crossing numbers of the remaining Cartesian products of each of the six connected four-vertex graphs with path graphs and star graphs. The results are summarised in \hyperref[tab-4vcart]{Table \ref{tab-4vcart}}.

\begin{table}[h!]\begin{center}\caption{Crossing numbers of Cartesian products of 4-vertex graphs with paths, cycles and stars. The results for $P_n$ are for $n \geq 1$, and the results for $C_n$ and $S_n$ are for $n \geq 3$.\label{tab-4vcart}}\smallskip\smallskip$\begin{array}{|c|c|c|c c|c|}\hline
\rule{0pt}{2.3ex} i        & G^4_i                 &  cr(G^4_i \Box P_n)  & \multicolumn{2}{c|} {cr(G^4_i \Box C_n) } &cr(G^4_i \Box S_n) \\
\hline 5 & \ctikz{\ngraph{4c.2}} & 0                                          & 0 &                                         & 2\fl[n]\fl[n-1]                   \\
\hline 7 & \ctikz{\ngraph{4c.1}} & n-1                                        & \;\;\;\;\; n \;\; (n \geq 6), &\begin{array}{lcl} 1 & & (n=3)\\2 & & (n=4)\\4 & & (n=5)\end{array}              & 2\fl[n]\fl[n-1] + \fl[n]                    \\
\hline 8 & \ctikz{\ngraph{4c.4}} & 0                                          & \;\;\;\; 2n \;\; (n \geq 3), &\begin{array}{lcl} 4 & & (n=3)\end{array}                                         & 2\fl[n]\fl[n-1]                           \\
\hline 9 & \ctikz{\ngraph{4c.3}} & n-1                                        & n &                                         & 2\fl[n]\fl[n-1] + \fl[n]                  \\
\hline 10 & \ctikz{\ngraph{4c.5}} & n-1                                        & 2n &                                        & 2\fl[n]\fl[n-1] + \fl[n]                  \\
\hline 11 & \ctikz{\ngraph{4c.6}} & 2n                                         & 3n &                                        & 2\fl[n]\fl[n-1] + 2n                      \\
\hline\end{array}$\end{center}\end{table}

\subsection{Cartesian products with 5-vertex graphs}\label{sec-cart5}

There are 21 connected non-isomorphic graphs on 5 vertices. We denote the connected graphs on 5 vertices as $G^5_i$, as defined in \hyperref[tab-5vcart]{Table \ref{tab-5vcart}}. For most of them, the crossing numbers of their Cartesian products with cycles, stars and paths are known. The results are summarised in \hyperref[tab-5vcart]{Table \ref{tab-5vcart}}, with the graphs indexed in the order originally designated by Kle\v{s}\v{c} in \cite{klesc2001_2}. The main contributions are discussed below and a list of publications where each result was first proved is displayed in \hyperref[app-c5]{Appendix \ref{app-c5}}.

\begin{table}[htbp]\begin{center}\vspace*{-0.7cm}\caption{Crossing numbers of Cartesian products of 5-vertex graphs with paths, cycles and stars. Unless otherwise stated, the results for $P_n$ are for $n \geq 1$, and the results for $C_n$ and $S_n$ are for $n \geq 3$. Empty cells imply the crossing number is not known.\label{tab-5vcart}}\smallskip\smallskip\scalebox{0.95}{$\begin{array}{|c|c|c|cc|c|}\hline
\rule{0pt}{2.3ex} i         & G^5_i                  & cr(G^5_i \Box P_n) & \multicolumn{2}{c|}{cr(G^5_i \Box C_n)}                                                     & cr(G^5_i \Box S_n)           \\
\hline 1  & \ctikz{\ngraph{5.1}}   & 0                                & \;\;\; 0                                   &                                                                       & 3 \fl[n]\fl[n-1]                    \\
\hline 2  & \ctikz{\ngraph{5.2}}   & 2n-2                             & \;\;\;\;\; 2n \;\; (n \geq 6),             & \begin{array}{lcl} 2 & & (n=3)\\4 & & (n=4)\\8 & & (n=5)\end{array}   & n(n-1)                              \\
\hline 3  & \ctikz{\ngraph{5.3}}   & n-1                              & \;\;\;\;\; n \;\; (n \geq 6),              & \begin{array}{lcl} 1 & & (n=3)\\2 & & (n=4)\\4 & & (n=5)\end{array}   & 3 \fl[n] \fl[n-1] + \fl[n]          \\
\hline 4  & \ctikz{\ngraph{5.4}}   & n-1                              & \;\;\; n                                   &                                                                       & 3 \fl[n] \fl[n-1] + \fl[n]          \\
\hline 5  & \ctikz{\ngraph{5.5}}   & n-1                              & \;\;\; n                                   &                                                                       &                                     \\
\hline 6  & \ctikz{\ngraph{5.6}}   & 2n-2                             & \;\;\;\;\; 2n \;\; (n \geq 6),             & \begin{array}{lcl} 4 & & (n=3)\\6 & & (n=4)\\9 & & (n=5)\end{array}   & n(n-1)                              \\
\hline 7  & \ctikz{\ngraph{5.7}}   & n-1                              & \;\;\;\;\; 2n \;\; (n \geq 4),             & \begin{array}{lcl} 4 & & (n=3)\end{array}                             & \unv{4 \fl[n] \fl[n-1] + \fl[n] }         \\
\hline 8  & \ctikz{\ngraph{5.8}}   & 0                                & \;\;\;\;\; 3n \;\; (n \geq 5),             & \begin{array}{ccl} 5 & & (n=3)\\10 & & (n=4)\end{array}               &                                     \\
\hline 9  & \ctikz{\ngraph{5.9}}   & 2n-2                             & \;\;\; 2n                                  &                                                                       & n(n-1)                              \\
\hline 10 & \ctikz{\ngraph{5.10}}  & 2n                               &                                            & \begin{array}{lcl} 9 & & (n=3)\end{array}                             & 4 \fl[n] \fl[n-1] + 2n              \\
\hline 11 & \ctikz{\ngraph{5.11}}  & 2n-2                             & \;\;\;\;\; 3n \;\; (n \geq 4),             & \begin{array}{lcl} 7 & & (n=3)\end{array}                             & n(n-1)                              \\
\hline 12 & \ctikz{\ngraph{5.12}}  & 2n-2                             & \;\;\; 2n                                  &                                                                       &                            \\
\hline 13 & \ctikz{\ngraph{5.13}}  & n-1                              & \;\;\;\;\; 3n \;\; (n \geq 4),             & \begin{array}{lcl} 7 & & (n=3)\end{array}                             &           \\
\hline 14 & \ctikz{\ngraph{5.14}}  & 2n-2                             & \;\;\; 3n                                  &                                                                       & n(n-1)                              \\
\hline 15 & \ctikz{\ngraph{5.15}}  & 3n-1                             &                                            &                                                                       & 4 \fl[n] \fl[n-1] + 2n + \fl[n]     \\
\hline 16 & \ctikz{\ngraph{5.16}}  & 3n-1                             & \;\; 2\left(n + \fl[n+1]\right)            &                                                                       &                                     \\
\hline 17 & \ctikz{\ngraph{5.17}}  & 2n                               &                                            &                                                                       & 4 \fl[n] \fl[n-1] + 2n              \\
\hline 18 & \ctikz{\ngraph{5.18}}  & 3n-1                             &                                            &                                                                       & 4 \fl[n] \fl[n-1] + 2n + \fl[n]     \\
\hline 19 & \ctikz{\ngraph{5.19}}  & 3n-1                             &                                            &                                                                       & 4 \fl[n] \fl[n-1] + 2n + \fl[n]     \\
\hline 20 & \ctikz{\ngraph{5.20}}  & 4n                               &                                            &                                                                       &\unv{ 4 \fl[n] \fl[n-1] + 4n}              \\
\hline 21 & \ctikz{\ngraph{5.21}}  & 6n                               & \;\;\; 9n                                  &                                                                       & 4 \fl[n] \fl[n-1] + \fl[n] + 5n + 1 \\
\hline\end{array}$}\end{center}\end{table}

For paths, the crossing number has been determined for all 21 graphs. The majority of the results were first determined in Kle\v{s}\v{c} (2001) \cite{klesc2001_2}. In particular, $cr(G^5_i \Box P_n)$ was first determined in \cite{klesc2001_2} for $i = 3, 4, 5, 6, 7, 9, 11, 13, 14, 17, 19, 20$. In addition, it is easy to check that $G^5_1 \Box P_n$ and $G^5_8 \Box P_n$ are planar.

For cycles, the crossing number has been determined for 15 graphs, with the cases $i = 10, 15, 17, 18, 19, 20$ still unknown. The majority of cases were settled in Kle\v{s}\v{c} (2001) \cite{klesc2001_2} $(i = 4, 5, 9, 12)$ and Kle\v{s}\v{c} (2005) \cite{klesc2005} $(i = 3, 6, 7, 13, 14)$, although each of the latter were first stated (without proof) in \cite{klesc2001_2}. For $G^5_9 \Box C_n$, the result was claimed in \cite{klesc2001_2} to hold only for $n \geq 6$, but was then later shown to hold for all $n \geq 3$ in Kle\v{s}\v{c} (2005) \cite{klesc2005}. In addition, Kle\v{s}\v{c} (2002) \cite{klesc2002} also determined an upper bound for $G^5_{10} \Box C_n$.

\begin{theorem}[Kle\v{s}\v{c}, 2002 \cite{klesc2002}]$cr(G^5_{10} \Box C_3) = 9$, and for $n\geq4$, \[cr(G^5_{10} \Box C_n) \leq 4n.\]\end{theorem}

For stars, the crossing number has been determined for 16 graphs, with the cases $i = 5, 8, 12, 13, 16$ still unknown.  The result for $G^5_2 \Box S_n$ was first determined by Huang and Zhao (2008) \cite{huangzhao2008} and independently by Ho (2008) \cite{ho2008} by recognizing it as a subdivision of $K_{1,4,n}$.

\subsection{Cartesian products with 6-vertex graphs}\label{sec-cart6}

In this section, we gather the known crossing numbers for Cartesian products of various 6-vertex graphs with paths, cycles and stars. There are, up to isomorphism, 112 connected graphs on 6 vertices. In order to refer to these graphs in a consistent manner, we give each 6-vertex graph a label, in the order originally designated by Harary \cite{harary1969}. We denote the graphs on 6 vertices as $G^6_i$, as defined for all indices, in \hyperref[app-g6]{Appendix \ref{app-g6}} and also repeated for the graphs of interest in the Tables of this section.  The list in \hyperref[app-g6]{Appendix \ref{app-g6}} includes disconnected graphs, which brings the total to 156 graphs, but the disconnected graphs are only used in \hyperref[sec-join]{Section \ref{sec-join}}.

\vspace*{-0.15cm}\subsubsection{Paths}\label{subsec-pathc6}

In 2013, Kle\v{s}\v{c} and Petrillov\'{a} \cite{klescpetrillova2013_2} gave a summary of known results, including the crossing numbers of Cartesian products of path graphs with forty different graphs on 6 vertices. The majority of those results were first determined in \cite{klescpetrillova2013_2}, including $cr(G^6_i \Box P_n)$ for $i =$ 26, 27, 28, 29, 41, 43, 44, 45, 46, 47, 48, 53, 54, 59, 60, 61, 64, 66, 72, 73, 74, 77, 79, 80, 83, 85, 86, 94, 104, 111, 121. For $G^6_{121}$, the result relied on a previous paper which has not undergone peer review, and so we mark that result with an asterisk. These, along with the remaining settled cases are displayed in \hyperref[tab-6vcart_path]{Table \ref{tab-6vcart_path}}. A list of publications where each result was first proved is displayed in \hyperref[app-c6p]{Appendix \ref{app-c6p}}. In total, the crossing number of $G^6_i \Box P_n$ has been settled for 52 of the 6-vertex graphs to date, with an additional 7 results claimed in papers which have not undergone adequate peer review. The latter are marked with asterisks.

Note that all results are for $n \geq 1$ except the following cases:

\begin{itemize} \item $G^6_{62} \Box P_n$ : the result only holds for $n \geq 3$. We also have $cr(G^6_{62} \Box P_1) = 0$ and $cr(G^6_{62} \Box P_2) = 2$ \cite{petrillova2019}.
\end{itemize}

\vspace*{-0.15cm}\subsubsection{Cycles}

To date, the crossing number has been determined for 31 cases, with the majority of these being collated by Dra\v{z}ensk\'{a} and Kle\v{s}\v{c} (2011) \cite{drazenskaklesc2011} and Clancy et al (to appear) \cite{clancyetaltoappear}.  These results are listed in \hyperref[tab-6vcart_cycle]{Table \ref{tab-6vcart_cycle}} and a list of publications where each result was first proved is displayed in \hyperref[app-c6c]{Appendix \ref{app-c6c}}.

Note that for $i = 26, 27, 28$ the proof given in \cite{drazenskaklesc2007} claims the result holds for $n \geq 3$, but in fact it holds only for $n \geq 6$. We have provided the correct numbers for $n = 3, 4, 5$. See Section \ref{app-error} for more details.



\vspace*{-0.15cm}\subsubsection{Stars}

The crossing numbers for the Cartesian product of many 6-vertex graphs and stars were provided in 2013 by Kle\v{s}\v{c} and Schr\"{o}tter \cite{klescschrotter2013} who made an attempt to gather the results known to them at the time. In that paper, they determined the crossing number of $G^6_i \Box S_n$, for seventeen graphs, however some of them had been previously determined. Only the cases $i =$ 27, 31, 43, 47, 48, 53, 59, 72, 73, 77, 79, 80, 104 were newly settled in \cite{klescschrotter2013}, although the case of $i = 27$ could also be seen as a corollary of Bokal (2007) \cite{bokal2007_2} and Huang and Zhao (2008) \cite{huangzhao2008}.  In total, the crossing number of $G^6_i \Box S_n$ has been settled for 21 graphs to date and these are displayed in \hyperref[tab-6vcart_star]{Table \ref{tab-6vcart_star}}. A further 8 cases have had results claimed in papers which have not undergone adequate peer review; these are marked by asterisks. A list of publications where each result was first proved is displayed in \hyperref[app-c6s]{Appendix \ref{app-c6s}}.

\vspace*{-0.5cm}\begin{table}[htbp]\begin{center}\caption{Known crossing numbers of Cartesian products of 6-vertex graphs with paths. All results are for $n \geq 1$ except where noted in Section \ref{subsec-pathc6}.\label{tab-6vcart_path}}\smallskip\smallskip\scalebox{1.0}{$\begin{array}{|c|c|c||c|c|c||c|c|c|}\hline
\rule{0pt}{2.3ex} i & G^6_i 	& cr(G^6_i \Box P_n) 	& i & G^6_i 	& cr(G^6_i \Box P_n) 	& i & G^6_i 	& cr(G^6_i \Box P_n)	\\
\hline 25&\ctikz{\ngraph{6.25}}  	& 0	        &61& \ctikz{\ngraph{6.61}}  	& 2n      	&89& \ctikz{\ngraph{6.89}}  	& 3n - 3        \\
\hline 26&\ctikz{\ngraph{6.26}}  	& n-1     	&62& \ctikz{\ngraph{6.62}}  	& 3n - 5  	&90& \ctikz{\ngraph{6.90}}  	& 3n - 3        \\
\hline 27&\ctikz{\ngraph{6.27}}  	& 2n - 2  	&63& \ctikz{\ngraph{6.63}}  	& 2n - 2  	&91& \ctikz{\ngraph{6.91}}  	& 3n - 1        \\
\hline 28&\ctikz{\ngraph{6.28}}  	& n - 1   	&64& \ctikz{\ngraph{6.64}}  	& 2n - 2  	&93 &\ctikz{\ngraph{6.93}}  	& \unv{4n} 	    \\
\hline 29&\ctikz{\ngraph{6.29}}  	& 2n - 2  	&65& \ctikz{\ngraph{6.65}}  	& 3n - 3  	&94 &\ctikz{\ngraph{6.94}}  	& 2n - 2 	    \\
\hline 31&\ctikz{\ngraph{6.31}}  	& 4n - 4  	&66& \ctikz{\ngraph{6.66}}  	& 2n - 2  	&103 &\ctikz{\ngraph{6.103}} 	& 6n - 2 	    \\
\hline 40&\ctikz{\ngraph{6.40}}  	& 0	        &68& \ctikz{\ngraph{6.68}}  	& 3n - 1  	&104 &\ctikz{\ngraph{6.104}} 	& 4n - 4 	   	\\
\hline 41&\ctikz{\ngraph{6.41}}  	& n - 1  	&70& \ctikz{\ngraph{6.70}}  	& 3n - 3  	&109 &\ctikz{\ngraph{6.109}} 	& \unv{4n}      \\
\hline 42&\ctikz{\ngraph{6.42}}  	& 2n - 4  	&71& \ctikz{\ngraph{6.71}}  	& 3n - 1  	&111 &\ctikz{\ngraph{6.111}} 	& 3n - 1 	    \\
\hline 43&\ctikz{\ngraph{6.43}}  	& n - 1   	&72& \ctikz{\ngraph{6.72}}  	& 4n - 4  	&113 &\ctikz{\ngraph{6.113}} 	& 4n - 4 	   	\\
\hline 44&\ctikz{\ngraph{6.44}}  	& 2n - 2 	&73& \ctikz{\ngraph{6.73}}  	& 4n - 4  	&119 &\ctikz{\ngraph{6.119}} 	& \unv{7n - 1} 	\\
\hline 45&\ctikz{\ngraph{6.45}}  	& 2n - 2  	&74& \ctikz{\ngraph{6.74}}  	& 2n - 2  	&120 &\ctikz{\ngraph{6.120}} 	& 3n - 3     	\\
\hline 46&\ctikz{\ngraph{6.46}}  	& n - 1   	&75& \ctikz{\ngraph{6.75}}  	& 2n      	&121 &\ctikz{\ngraph{6.121}} 	& \unv{4n}      \\
\hline 47&\ctikz{\ngraph{6.47}}	    & 2n - 2  	&77& \ctikz{\ngraph{6.77}}  	& 2n - 2  	&125 &\ctikz{\ngraph{6.125}} 	& 5n - 3 	   	\\
\hline 48&\ctikz{\ngraph{6.48}}	    & 4n - 4  	&79& \ctikz{\ngraph{6.79}}  	& 4n - 4  	&146 &\ctikz{\ngraph{6.146}} 	& \unv{5n - 1} 	\\
\hline 51&\ctikz{\ngraph{6.51}}	    & 3n - 3  	&80& \ctikz{\ngraph{6.80}}  	& 4n - 4  	&152 &\ctikz{\ngraph{6.152}} 	& \unv{6n}      \\
\hline 53&\ctikz{\ngraph{6.53}}	    & 2n - 2  	&83& \ctikz{\ngraph{6.83}}  	& 2n - 2  	&154 &\ctikz{\ngraph{6.154}} 	& 9n - 1 	   	\\
\hline 54&\ctikz{\ngraph{6.54}}	    & 2n - 2  	&85& \ctikz{\ngraph{6.85}}  	& 2n      	&155 &\ctikz{\ngraph{6.155}} 	& \unv{12n}     \\
\hline 59&\ctikz{\ngraph{6.59}}	    & 2n - 2  	&86& \ctikz{\ngraph{6.86}}  	& 3n - 1  	&156 &\ctikz{\ngraph{6.156}} 	& 15n + 3	    \\
\hline 60&\ctikz{\ngraph{6.60}}     & n - 1     &87& \ctikz{\ngraph{6.87}}  	& 3n - 1  	&    &                          &               \\
\hline
\end{array}$}\end{center}\end{table}

\begin{table}[htbp]\begin{center}\caption{Crossing numbers of Cartesian products of 6-vertex graphs with cycles. \label{tab-6vcart_cycle}}\scalebox{0.825}{$\begin{array}{|c|c|c c||c|c|c c|}
\hline\rule{0pt}{2.3ex} i & G^6_i & cr(G^6_i \Box C_n) & \mbox{(small cases)} & i & G^6_i & cr(G^6_i \Box C_n) & \mbox{(small cases)} \\
\hline 25                 & \ctikz{\ngraph{6.25}}  & \begin{array}{rl}0  	& \end{array}              &                                                                   &60                 & \ctikz{\ngraph{6.60}}  & \begin{array}{rl}4n    & (n \geq 6)\end{array}    & \begin{array}{rl} 8  & (n=3)\\16  & (n=4)\\20  & (n=5)\end{array}  \\
\hline 26                 & \ctikz{\ngraph{6.26}}  & \begin{array}{rl}\;\;n & (n \geq 6)\end{array}    & \begin{array}{rl} 1  & (n=3)\\2   & (n=4)\\4   & (n=5)\end{array} &63                 & \ctikz{\ngraph{6.63}}  & \begin{array}{rl}2n    & (n \geq 4)\end{array}    & \begin{array}{rl} 6  & (n=3)\end{array}                            \\
\hline 27                 & \ctikz{\ngraph{6.27}}  & \begin{array}{rl}2n    & (n \geq 6)\end{array}    & \begin{array}{rl} 2  & (n=3)\\4   & (n=4)\\8   & (n=5)\end{array} &64                 & \ctikz{\ngraph{6.64}}  & \begin{array}{rl}2n    & (n \geq 6)\end{array}    & \begin{array}{rl}6  & (n=3)\\8  & (n=4)\\10  & (n=5)\end{array}    \\
\hline 28                 & \ctikz{\ngraph{6.28}}  & \begin{array}{rl}\;\;n & (n \geq 6)\end{array}    & \begin{array}{rl} 1  & (n=3)\\2   & (n=4)\\4   & (n=5)\end{array} &66                 & \ctikz{\ngraph{6.66}}  & \begin{array}{rl}3n    & (n \geq 5)\end{array}    & \begin{array}{rl} 7  & (n=3)\\12  & (n=4)\end{array}               \\
\hline 29                 & \ctikz{\ngraph{6.29}}  & \begin{array}{rl}2n    & (n \geq 6)\end{array}    & \begin{array}{rl} 2  & (n=3)\\4   & (n=4)\\8   & (n=5)\end{array} &67                 & \ctikz{\ngraph{6.67}}  & \begin{array}{rl}3n    & (n \geq 4)\end{array}    & \begin{array}{rl} \\7  & (n=3)\\  & \end{array}                    \\
\hline 40                 & \ctikz{\ngraph{6.40}}  & \begin{array}{rl}4n    & (n \geq 6)\end{array}    & \begin{array}{rl} 6  & (n=3)\\12  & (n=4)\\18  & (n=5)\end{array} &70                 & \ctikz{\ngraph{6.70}}  & \begin{array}{rl}3n    & (n \geq 5)\end{array}    & \begin{array}{rl} 7 & (n=3)\\12  & (n=4)\end{array}                \\
\hline 41                 & \ctikz{\ngraph{6.41}}  & \begin{array}{rl}3n    & (n \geq 5)\end{array}    & \begin{array}{rl} 5  & (n=3)\\10  & (n=4)\end{array}              &75                 & \ctikz{\ngraph{6.75}}  & \begin{array}{rl}2n    & (n \geq 4)\end{array}    & \begin{array}{rl} 6  & (n=3)\end{array}                            \\
\hline 42                 & \ctikz{\ngraph{6.42}}  & \begin{array}{rl}2n    & (n \geq 4)\end{array}    & \begin{array}{rl} 4  & (n=3)\end{array}                           &77                 & \ctikz{\ngraph{6.77}}  & \begin{array}{rl}2n    & (n \geq 6)\end{array}    & \begin{array}{rl} 6 & (n=3)\\8  & (n=4)\\10  & (n=5)\end{array}    \\
\hline 43                 & \ctikz{\ngraph{6.43}}  & \begin{array}{rl}\;\;n & (n \geq 3)\end{array}    &                                                                   &78                 & \ctikz{\ngraph{6.78}}  & \begin{array}{rl}3n    & (n \geq 6)\end{array}    & \begin{array}{rl} 7  & (n=3)\\10  & (n=4)\\14  & (n=5)\end{array}  \\
\hline 44                 & \ctikz{\ngraph{6.44}}  & \begin{array}{rl}2n    & (n \geq 4)\end{array}    & \begin{array}{rl} 4  & (n=3)\end{array}                           &83                 & \ctikz{\ngraph{6.83}}  & \begin{array}{rl}4n    & (n \geq 6)\end{array}    & \begin{array}{rl} 10  & (n=3)\\16  & (n=4)\\20  & (n=5)\end{array} \\
\hline 46                 & \ctikz{\ngraph{6.46}}  & \begin{array}{rl}\;\;n & (n \geq 3)\end{array}    &                                                                   &90                 & \ctikz{\ngraph{6.90}}  & \begin{array}{rl}4n    & (n \geq 6)\end{array}    & \begin{array}{rl} 11 & (n=3)\\16  & (n=4)\\20  & (n=5)\end{array}  \\
\hline 47                 & \ctikz{\ngraph{6.47}}  & \begin{array}{rl}2n    & (n \geq 6)\end{array}    & \begin{array}{rl} 4  & (n=3)\\6  & (n=4)\\9  & (n=5)\end{array}   &92                 & \ctikz{\ngraph{6.92}}  & \begin{array}{rl}3n    & (n \geq 4)\end{array}    & \begin{array}{rl} 9  & (n=3)\end{array}                            \\
\hline 49                 & \ctikz{\ngraph{6.49}}  & \begin{array}{rl}2n    & (n \geq 4)\end{array}    & \begin{array}{rl} 4  & (n=3)\end{array}                           &98                 & \ctikz{\ngraph{6.98}}  & \begin{array}{rl}3n    & (n \geq 5)\end{array}    & \begin{array}{rl} 9  & (n=3)\\12 & (n=4)\end{array}                \\
\hline 53                 & \ctikz{\ngraph{6.53}}  & \begin{array}{rl}2n    & (n \geq 6)\end{array}    & \begin{array}{rl} 4  & (n=3)\\6  & (n=4)\\9  & (n=5)\end{array}   &113                & \ctikz{\ngraph{6.113}} & \begin{array}{rl}4n    & (n \geq 3)\end{array}    &                                                                    \\
\hline 54                 & \ctikz{\ngraph{6.54}}  & \begin{array}{rl}2n    & (n \geq 6)\end{array}    & \begin{array}{rl} 4  & (n=3)\\6  & (n=4)\\9  & (n=5)\end{array}   &156                & \ctikz{\ngraph{6.156}} & \begin{array}{rl}18n   & (n \geq 3)\end{array}    &                                                                    \\
\hline 59                 & \ctikz{\ngraph{6.59}}  & \begin{array}{rl}4n    & (n \geq 6)\end{array}    & \begin{array}{rl} 8  & (n=3)\\16  & (n=4)\\20  & (n=5)\end{array} &		            &						 &						                                                                                                 & \\
\hline\end{array}$}\end{center}\end{table}

\renewcommand{\arraystretch}{1.7}
\begin{table}[htbp]\begin{center}\vspace*{-0.5cm}\caption{Crossing numbers of Cartesian products of 6-vertex graphs with stars. \label{tab-6vcart_star}}\smallskip\smallskip $\begin{array}{|c|c|c||c|c|c|}\hline i & G^6_i & cr(G^6_i \Box S_n)  & i & G^6_i & cr(G^6_i \Box S_n) \\
\hline \rule{0pt}{2.3ex} 25     & \ctikz{\ngraph{6.25}}  & 4\fl[n]\fl[n-1]                  & 77     & \ctikz{\ngraph{6.77}}  & 4\fl[n]\fl[n-1] + 2\fl[n]              \\
\hline \rule{0pt}{2.3ex} 26     & \ctikz{\ngraph{6.26}}  & 4\fl[n]\fl[n-1] + \fl[n]         & 79     & \ctikz{\ngraph{6.79}}  & 6\fl[n]\fl[n-1] + 4\fl[n]              \\
\hline \rule{0pt}{2.3ex} 27     & \ctikz{\ngraph{6.27}}  & 5\fl[n]\fl[n-1] + 2\fl[n]        & 80     & \ctikz{\ngraph{6.80}}  & 6\fl[n]\fl[n-1] + 4\fl[n]              \\
\hline \rule{0pt}{2.3ex} 28     & \ctikz{\ngraph{6.28}}  & 4\fl[n]\fl[n-1] + \fl[n]         & 85     & \ctikz{\ngraph{6.85}}  & \unv{6\fl[n]\fl[n-1] + 2n}             \\
\hline \rule{0pt}{2.3ex} 29     & \ctikz{\ngraph{6.29}}  & 4\fl[n]\fl[n-1] + 2\fl[n]        & 93     & \ctikz{\ngraph{6.93}}  & \unv{6\fl[n]\fl[n-1] + 4n}             \\
\hline \rule{0pt}{2.3ex} 31     & \ctikz{\ngraph{6.31}}  & 6\fl[n]\fl[n-1] + 4\fl[n]        & 94     & \ctikz{\ngraph{6.94}}  & 6\fl[n]\fl[n-1] + 2\fl[n]              \\
\hline \rule{0pt}{2.3ex} 43     & \ctikz{\ngraph{6.43}}  & 4\fl[n]\fl[n-1] + \fl[n]         & 104    & \ctikz{\ngraph{6.104}} & 6\fl[n]\fl[n-1] + 4\fl[n]              \\
\hline                   47     & \ctikz{\ngraph{6.47}}  & 5\fl[n]\fl[n-1] + 2\fl[n]        & 111    & \ctikz{\ngraph{6.111}} & \unv{6\fl[n]\fl[n-1] + 2\fl[n] + 2n }  \\
\hline \rule{0pt}{2.3ex} 48     & \ctikz{\ngraph{6.48}}  & 6\fl[n]\fl[n-1] + 4\fl[n]        & 120    & \ctikz{\ngraph{6.120}} & \unv{6\fl[n]\fl[n-1] + 3\fl[n]}        \\
\hline                   53     & \ctikz{\ngraph{6.53}}  & 4\fl[n]\fl[n-1] + 2\fl[n]        & 124    & \ctikz{\ngraph{6.124}} & \unv{6\fl[n]\fl[n-1] + 2n + 3\fl[n]}   \\
\hline \rule{0pt}{2.3ex} 59     & \ctikz{\ngraph{6.59}}  & 6 \fl[n]\fl[n-1] + 2\fl[n]       & 125    & \ctikz{\ngraph{6.125}} & 6\fl[n]\fl[n-1] + 3\fl[n] + 2n         \\
\hline                   61     & \ctikz{\ngraph{6.61}}  & \unv{6\fl[n]\fl[n-1] + 2n}       & 130    & \ctikz{\ngraph{6.130}} & \unv{6\fl[n]\fl[n-1] + 4n }            \\
\hline \rule{0pt}{2.3ex} 62     & \ctikz{\ngraph{6.62}}  & 5\fl[n]\fl[n-1] + 2\fl[n]        & 137    & \ctikz{\ngraph{6.137}} & \unv{6\fl[n]\fl[n-1] + 4n }            \\
\hline \rule{0pt}{2.3ex} 72     & \ctikz{\ngraph{6.72}}  & 6\fl[n]\fl[n-1] + 4\fl[n]        & 152    & \ctikz{\ngraph{6.152}} & 6\fl[n]\fl[n-1] + 6n                   \\
\hline \rule{0pt}{2.3ex} 73     & \ctikz{\ngraph{6.73}}  & 6\fl[n]\fl[n-1] + 4\fl[n]        &        &                        &                                        \\
\hline
\end{array}$\end{center}\end{table}

\null
\vfill
\newpage

\subsection{Cartesian products with graphs on 7 or more vertices}

\subsubsection{7-vertex graphs with paths, cycles and stars}

In 2005, He and Huang \cite{hehuang2005} $\Asterisk$ considered six graphs on seven vertices and determined the crossing number of their Cartesian products with paths. The results are summarised in \hyperref[tab-7vcart_path]{Table \ref{tab-7vcart_path}}, with those graphs labelled $G^7_1$ to $G^7_6$. The Cartesian product of two additional graphs with paths were considered by Liu et al.\ (2012) \cite{liuetal2012} $\Asterisk$ and Ding et al.\ (2018) \cite{dingetal2018}; those graphs are labelled $G^7_7$ and $G^7_8$ respectively. There are also a number of known results for Cartesian products of named 7-vertex graphs with paths, cycles and stars. These are provided in \hyperref[tab-cart7named]{Table \ref{tab-cart7named}} along with a list of the publications where the result was proved.

\begin{table}[htbp]\begin{center}\caption{Crossing numbers of Cartesian products of 7-vertex graphs with paths. All results are for $n \geq 1$.\label{tab-7vcart_path}}\smallskip\smallskip$\begin{array}{|c|c|c||c|c|c||c|c|c|}\hline \rule{0pt}{2.3ex} i & G^7_i & cr(G^7_i \Box P_n) & i & G^7_i & cr(G^7_i \Box P_n) & i & G^7_i & cr(G^7_i \Box P_n)\\
\hline 1   &  \ctikz{\ngraph{7p.1}}  & \unv{3n - 3}   & 4  &  \ctikz{\ngraph{7p.4}}  & \unv{n - 1}    & 7  &  \ctikz{\ngraph{7p.7}}  & \unv{5n - 1} \\
\hline 2   &  \ctikz{\ngraph{7p.2}}  & \unv{4n - 4}   & 5  &  \ctikz{\ngraph{7p.5}}  & \unv{2n}       & 8  &  \ctikz{\ngraph{7p.8}}  & 4n - 4 \\
\hline 3   &  \ctikz{\ngraph{7p.3}}  & \unv{2n - 2}   & 6  &  \ctikz{\ngraph{7p.6}}  & \unv{2n - 2}   &    &                         &        \\
\hline
\end{array}$\end{center}\end{table}

{\renewcommand{\arraystretch}{1.5}
\begin{table}[htbp]\begin{center}\caption{Crossing numbers of Cartesian products of named 7-vertex graphs with paths, cycles and stars.\label{tab-cart7named}}\smallskip\smallskip$\begin{array}{|l|l|l|}\hline
\rule{0pt}{2.3ex} \text{{\bf Graph family} }       & \text{{\bf Crossing number}} & \text{{\bf Publication}}     \\
\hline S_6 \Box P_n & 6n - 6 & \text{Bokal (2007) \cite{bokal2007}} \\
\hline Ci_{7}(\{1,2\}) \Box P_n & 8n \: \Asterisk  & \text{Yuan et al.\ (2008) \cite{yuanhuangetal2008}} \: \Asterisk\\
\hline W_6 \Box P_n & 7n - 5  & \text{Bokal (2007) \cite{bokal2007_2}} \\
\hline K_{2,5} \Box P_n & 8n  & \text{Tang et al.\ (2007) \cite{tangetal2007}} \\
\hline K_{1,1,5} \Box P_n & 10n-2 & \text{Ouyang et al.\ (2014) \cite{ouyangetal2014}}\\
\hline K_7 \Box P_n & 30n+6  & \text{Ouyang et al.\ (2014) \cite{ouyangetal2014_2}} \\
\hline P_6 \Box S_n & 5\fl[n]\fl[n-1] & \text{Bokal (2007) \cite{bokal2007}}\\
\hline C_7 \Box C_n & 5n & \text{Adamsson and Richter (2004) \cite{adamssonrichter2004}}\\
\hline K_7 \Box C_n & 36n & \text{Zheng et al.\ (2008) \cite{zhengetal2008}}\\
\hline\end{array}$\end{center}\end{table}
{\renewcommand{\arraystretch}{1}




Additionally, the crossing number of the Cartesian product of a star with any 7-vertex tree except $S_6$ is known due to Bokal (2007) \cite{bokal2007_2} (see \hyperref[sec-startree]{Section \ref{sec-startree}}).

\subsubsection{8-vertex graphs with paths, cycles and stars}

The results in \hyperref[tab-8vcart_path]{Table \ref{tab-8vcart_path}} are due to Yuan and Huang (2007) \cite{yuanhuang2007} $\Asterisk$, with the exception of $G^8_5$ which is due to Ding et al.\ (2018) \cite{dingetal2018}. There are also a number of known results for Cartesian products of named 8-vertex graphs with paths, cycles and stars.  These are provided in \hyperref[tab-cart8named]{Table \ref{tab-cart8named}} along with a list of the publications where the result was proved.

\begin{table}[htbp]\begin{center}\caption{Crossing numbers of Cartesian products of 8-vertex graphs with paths. All results are for $n \geq 1$.\label{tab-8vcart_path}}\smallskip\smallskip$\begin{array}{|c|c|c||c|c|c||c|c|c|}\hline \rule{0pt}{2.3ex} i & G^8_i & cr(G^8_i \Box P_n) & i & G^8_i & cr(G^8_i \Box P_n) & i & G^8_i & cr(G^8_i \Box P_n) \\
\hline 1 &    \ctikz{\ngraph{8p.1}}  & \unv{n - 1}    & 3 &    \ctikz{\ngraph{8p.3}}  & \unv{2n - 2} & 5 & \ctikz{\ngraph{8p.5}} & 4n - 4  \\
\hline 2 &    \ctikz{\ngraph{8p.2}}  & \unv{2n - 2}   & 4 &    \ctikz{\ngraph{8p.4}}  & \unv{2n - 2} & & &  \\
\hline
\end{array}$\end{center}\end{table}

{\renewcommand{\arraystretch}{1.5}
\begin{table}[htbp]\begin{center}\caption{Crossing numbers of Cartesian products of named 8-vertex graphs with paths, cycles and stars.\label{tab-cart8named}}\smallskip\smallskip$\begin{array}{|l|l|l|}\hline
\rule{0pt}{2.3ex} \text{{\bf Graph family} }                & \text{{\bf Crossing number}} & \text{{\bf Publication}}    \\
\hline S_7 \Box P_n & 9n - 9 & \text{Bokal (2007) \cite{bokal2007}} \\
\hline W_7 \Box P_n & 10n - 8  & \text{Bokal (2007) \cite{bokal2007_2}} \\
\hline Ci_{8}(\{1,2\}) \Box P_n & 8n  & \text{Yuan et al.\ (2008) \cite{yuanetal2008}} \\
\hline Ci_{8}(\{1,4\}) \Box P_n & 9n-1 \: \Asterisk  & \text{Wang and Ma (2017) \cite{wangma2017}} \: \Asterisk\\
\hline GP(4,1) \Box P_n & 8n \: \Asterisk & \text{Yuan and Huang (2011) \cite{yuanhuang2011_2}} \: \Asterisk\\
\hline K_{2,6} \Box P_n & 12n  & \text{Tang et al.\ (2007) \cite{tangetal2007}} \\
\hline K_{1,1,6} \Box P_n & 15n-3 & \text{Ouyang et al.\ (2014) \cite{ouyangetal2014}}\\
\hline K_8 \Box P_n & 54n+18 & \text{Ouyang et al.\ (2014) \cite{ouyangetal2014_2}}\\
\hline P_7 \Box S_n & 6\fl[n]\fl[n-1] & \text{Bokal (2007) \cite{bokal2007}}\\
\hline\end{array}$\end{center}\end{table}
{\renewcommand{\arraystretch}{1}

%
%

Additionally, the crossing number of the Cartesian product of a star with any 8-vertex tree with maximum degree 5 is known due to Bokal (2007) \cite{bokal2007_2} (see \hyperref[sec-startree]{Section \ref{sec-startree}}).

\subsection{Variants of toroidal grid graphs}

Toroidal grid graphs are another name for the Cartesian product of two cycles. In 2002, Foley et al.\ \cite{foleyetal2002} considered the crossing numbers of two variants of toroidal grid graphs; specifically, the twisted toroidal grid graph $\mathcal{T}(m,n)$ and the crossed toroidal grid graphs $\mathcal{X}(m,n)$. The twisted toroidal graph $\mathcal{T}(3,n)$ is equivalent to $C_3 \Box C_n$, except vertices 1, 2, 3 in the final cycle link to vertices 3, 1, 2 in the first cycle, respectively. The crossed toroidal graph $\mathcal{X}(3,n)$ is equivalent to $C_3 \Box C_n$, except vertices 1, 2, 3 in the final cycle link to vertices 1, 3, 2 in the first cycle, respectively. Examples of $\mathcal{T}(3,5)$ and $\mathcal{X}(3,5)$ are displayed in \hyperref[fig-toroidal]{Figure \ref{fig-toroidal}}.


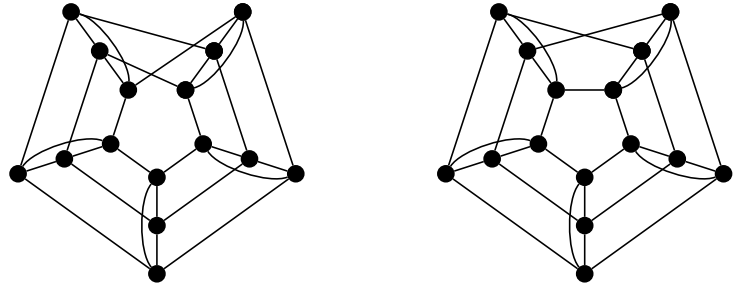
\begin{figure}[h!]\begin{center}
\begin{tikzpicture}[largegraph,scale=0.8]
\foreach \n in {0,...,5}{
  \foreach \l in {1,2,3}{
    \node (\l\n) at (\n*360/5+54: 0.8*\l cm) {};
  }
}
\foreach \n in {0,...,4}{
  \draw (1\n) -- (2\n);
  \draw (2\n) -- (3\n);
  \draw[bend left=40,looseness=0.6] (3\n) to (1\n);
}

\foreach \n in {1,...,4}{
  \foreach \l in {1,2,3}{
    \pgfmathtruncatemacro{\nextn}{\n+1}
    \draw (\l\n) -- (\l\nextn);
  }
}
\draw (30) to (11);
\draw (20) to (31);
\draw (10) to (21);
\end{tikzpicture}\;\;\;\;\;\;\;\;\;\;\;\;\;\;
\begin{tikzpicture}[largegraph,scale=0.8]
\foreach \n in {0,...,5}{
  \foreach \l in {1,2,3}{
    \node (\l\n) at (\n*360/5+54: 0.8*\l cm) {};
  }
}
\foreach \n in {0,...,4}{
  \draw (1\n) -- (2\n);
  \draw (2\n) -- (3\n);
  \draw[bend left=40,looseness=0.6] (3\n) to (1\n);
}

\foreach \n in {1,...,4}{
  \foreach \l in {1,2,3}{
    \pgfmathtruncatemacro{\nextn}{\n+1}
    \draw (\l\n) -- (\l\nextn);
  }
}
\draw (30) to (21);
\draw (20) to (31);
\draw (10) to (11);
\end{tikzpicture}
\caption{The twisted toroidal grid graph $\mathcal{T}(3,5)$ and the crossed toroidal grid graph $\mathcal{X}(3,5)$, respectively.\label{fig-toroidal}}\end{center}\end{figure}

\begin{theorem}[Foley et al.,\ 2002 \cite{foleyetal2002}]For $n \geq 3$, the following holds:
\[ cr(\mathcal{T}(3,n)) = cr(\mathcal{X}(3,n)) = n. \]\end{theorem}

\section{Join products of graphs}\label{sec-join}

The join product of two graphs $G$ and $H$, denoted $G + H$ and sometimes simply called the join of $G$ and $H$, is equal to the union of $G$ and $H$, plus edges linking every vertex of $G$ to every vertex of $H$. An example of the join product of two paths, $P_3 + P_4$, is displayed in \hyperref[fig-join]{Figure \ref{fig-join}} in two drawings, with the latter drawing being crossing-optimal.

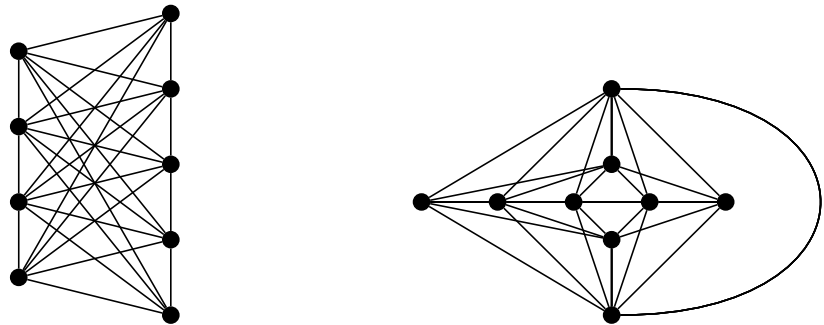
\begin{figure}[h!]\begin{center}$\;\;\;\;\;\;\;\;\;\;\;\;\;\;$\begin{tikzpicture}[largegraph,scale=0.5]
\foreach \n in {0,...,3}{
  \pgfmathtruncatemacro{\nextn}{2*\n-3}
  \node (a\n) at (0,\nextn) {};
}
\foreach \n in {0,...,4}{
  \pgfmathtruncatemacro{\nextn}{2*\n-4}
  \node (b\n) at (4,\nextn) {};
}
\foreach \n in {0,...,3}{
  \foreach \m in {0,...,4}{
    \draw (a\n)  -- (b\m);
  }
}
\foreach \n in {0,...,2}{
  \pgfmathtruncatemacro{\nextn}{\n + 1}
  \draw (a\n) -- (a\nextn);
}
\foreach \n in {0,...,3}{
  \pgfmathtruncatemacro{\nextn}{\n + 1}
  \draw (b\n) -- (b\nextn);
}

\end{tikzpicture}$\;\;\;\;\;\;\;\;\;\;\;\;\;\;\;\;\;\;\;\;\;\;\;\;\;$
\begin{tikzpicture}[largegraph,scale=0.5]
\foreach \n in {0,...,3}{
  \pgfmathtruncatemacro{\nextn}{2*\n-3}
  \node (a\n) at (0,\nextn) {};
}
\foreach \n in {0,...,4}{
  \pgfmathtruncatemacro{\nextn}{2*\n-5}
  \node (b\n) at (\nextn,0) {};
}
\foreach \n in {0,...,3}{
  \foreach \m in {0,...,4}{
    \draw (a\n)  -- (b\m);
  }
\draw (a0) -- (a1);
\draw (a2) -- (a3);
\draw[bend right=90,looseness=3] (a0) to (a3);

\draw (b0) -- (b1);
\draw (b1) -- (b2);
\draw (b2) -- (b3);
\draw (b3) -- (b4);

}\end{tikzpicture}\caption{Two drawings of the join product $P_3 + P_4$, with the second being optimal.\label{fig-join}}\end{center}\end{figure}

In effect, the join product connects $G$ and $H$ by the edges of the complete bipartite graph $K_{|G|,|H|}$. Hence, the crossing number of a join product is bounded below by the crossing number of the corresponding complete bipartite graph. Accordingly, the crossing number of the complete bipartite graph has appeared in all known results for join products to date. Within this section we will adopt the common practice of defining $Z(m,n)$ to be the conjectured crossing number of $K_{m,n}$ (see \hyperref[sec-complete-bipartite]{Section \ref{sec-complete-bipartite}}).

\[ Z(m,n) := \fl[m]\fl[m-1]\fl[n]\fl[n-1]. \]

Recall from \hyperref[sec-complete-bipartite]{Section \ref{sec-complete-bipartite}} that $Z(m,n)$ is only known to coincide with $cr(K_{m,n})$ if $\min\{m,n\} \leq 6$, or for the special cases when
$7 \leq \min\{m,n\} \leq 8$ and $7 \leq \max\{m,n\} \leq 10$. The crossing numbers of families of graphs resulting from join products have, to date, only been computed exactly when at least one graph involved has at most six vertices. Hence, this section will be divided into subsections corresponding to the size of the fixed graph in the join product. We highlight a few important points regarding this section before continuing.

Unlike for Cartesian products, the join product always results in a connected graph, even if one or both of the graphs involved is disconnected. Hence, we will also consider disconnected graphs in this section. In addition to cycles and paths which were considered extensively in \hyperref[sec-cartesian]{Section \ref{sec-cartesian}}, it is also common to consider the crossing number of join products of graphs with the discrete graph. The discrete graph $D_n$ is the graph with $n$ isolated vertices and no edges. It is also often denoted $nK_1$. To date it has been rare to consider the crossing number of join products involving arbitrarily large stars.

A number of interesting graphs can be viewed as being the result of join products. Most notably, complete multipartite graphs can be viewed as resulting from join products, in the following way. Consider the complete $k$-partite graph $K_{a_1,a_2,\hdots,a_k}$. Then $K_{a_1,a_2,\hdots,a_k} + D_n$ is isomorphic to the complete $(k+1)$-partite graph $K_{a_1,a_2,\hdots,a_k,n}$. Hence, a number of results for join products can be taken from the various publications on crossing numbers of complete multipartite graphs, and vice versa. We include such results in what follows, even if they were not originally presented as join product results.

Finally, when referring to join products, it is common in the literature to use the notation $P_n$ to refer to the path graph on $n$ vertices; this is contrary to the more standard usage of $P_n$ to refer to the path graph on $n+1$ vertices. The reason for this is that the order of each input graph is an important variable in join products, as the crossing number will inevitably contain $Z(m,n)$ for input graphs of size $m$ and $n$. Despite this common practice in the literature, for the sake of consistency we will maintain our notation from the rest of this survey, and use $P_{n-1}$ to refer to the path graph on $n$ vertices.

\subsection{Join products with 3-vertex graphs}

There are four graphs, up to isomorphism, on three vertices; see \hyperref[tab-3vjoin]{Table \ref{tab-3vjoin}}. The crossing numbers of the join products with discrete graphs, paths and cycles has been found for each of them. Each of the results for join products with $P_{n-1}$ and $C_n$ was found by Kle\v{s}\v{c} (2007) \cite{klesc2007}. Since every join product with a three vertex graph contains $K_{3,n}$ and hence a lower bound is $Z(3,n)$, and $D_n$ is a subgraph of $P_{n-1}$, the results for join products with $D_n$ can be seen as immediate corollaries.

\begin{table}[h!]\begin{center}\caption{Crossing numbers of joins of 3-vertex graphs with discrete graphs, paths and cycles. The results for $P_{n-1}$ are for $n \geq 2$ and the results for $C_n$ are for $n \geq 3$.\label{tab-3vjoin}}\smallskip\smallskip
$\begin{array}{|c|c|c|c|c|}\hline
\rule{0pt}{2.3ex} i        & G^3_i                & cr(G^3_i + D_n) & cr(G^3_i + P_{n-1}) & cr(G^3_i + C_n)\\
\hline 1 & \ctikz{\ngraph{3.1}} & Z(3,n)                & Z(3,n)                    & Z(3,n)               \\
\hline 2 & \ctikz{\ngraph{3.2}} & Z(3,n)                & Z(3,n)                    & Z(3,n)               \\
\hline 3 & \ctikz{\ngraph{3.3}} & Z(3,n)                & Z(3,n)                    & Z(3,n) + 1           \\
\hline 4 & \ctikz{\ngraph{3.4}} & Z(3,n)                & Z(3,n) + 1                & Z(3,n) + 2           \\
\hline\end{array}$\end{center}\end{table}

\subsection{Join products with 4-vertex graphs}

There are eleven graphs, up to isomorphism, on four vertices; see \hyperref[tab-4vjoin]{Table \ref{tab-4vjoin}}. The crossing numbers of the join products with discrete graphs, paths and cycles has been found for each of them.

\begin{table}[h!]\begin{center}\caption{Crossing numbers of joins of 4-vertex graphs with discrete graphs, paths and cycles. The results for $D_n$ are for $n \geq 1$, the results for $P_{n-1}$ are for $n \geq 2$, and the results for $C_n$ are for $n \geq 3$.\label{tab-4vjoin}}\smallskip\smallskip
$\begin{array}{|c|c|c|c|c|}\hline
\rule{0pt}{2.3ex} i         & G^4_i                  & cr(G^4_i + D_n) & cr(G^4_i + P_{n-1}) & cr(G^4_i + C_n)   \\
\hline 1  & \ctikz{\ngraph{4a.1}}  & Z(4,n)                & Z(4,n)                    & Z(4,n)                  \\
\hline 2  & \ctikz{\ngraph{4a.2}}  & Z(4,n)                & Z(4,n)                    & Z(4,n)                  \\
\hline 3  & \ctikz{\ngraph{4a.3}}  & Z(4,n)                & Z(4,n)                    & Z(4,n)                  \\
\hline 4  & \ctikz{\ngraph{4a.4}}  & Z(4,n)                & Z(4,n)                    & Z(4,n) + 1              \\
\hline 5  & \ctikz{\ngraph{4a.5}}  & Z(4,n)                & Z(4,n)                    & Z(4,n) + 1              \\
\hline 6  & \ctikz{\ngraph{4a.6}}  & Z(4,n) + \fl[n]       & Z(4,n) + \fl[n]           & Z(4,n) + \fl[n] + 2     \\
\hline 7  & \ctikz{\ngraph{4a.7}}  & Z(4,n) + \fl[n]       & Z(4,n) + \fl[n]           & Z(4,n) + \fl[n] + 2     \\
\hline 8  & \ctikz{\ngraph{4a.8}}  & Z(4,n)                & Z(4,n) + 1                & Z(4,n) + 2              \\
\hline 9  & \ctikz{\ngraph{4a.9}}  & Z(4,n) + \fl[n]       & Z(4,n) + \fl[n]           & Z(4,n) + \fl[n] + 2     \\
\hline 10 & \ctikz{\ngraph{4a.10}} & Z(4,n) + \fl[n] & Z(4,n) + \fl[n] + 1       & Z(4,n) + \fl[n] + 3     \\
\hline 11 & \ctikz{\ngraph{4a.11}} & Z(4,n) + n            & Z(4,n) + n + 1            & Z(4,n) + n + 4          \\
\hline\end{array}$\end{center}\end{table}

Each of the results for join products with $P_{n-1}$ and $C_n$ were first shown by Kle\v{s}\v{c} (2007) \cite{klesc2007}. The results for join products with $D_n$ were first shown by Kle\v{s}\v{c} and Schr\"{o}tter (2011) \cite{klescschrotter2011}, with three exceptions. First, $G^4_1 + D_n = K_{4,n}$, for which the crossing number was first determined by Guy (1969) \cite{guy1969}. Second, $G^4_7 + D_n = K_{1,3,n}$, for which the crossing number was first determined by Asano (1986) \cite{asano1986}. Finally, $G^4_{11} + D_n = K_{1,1,1,1,n}$, for which the crossing number was first determined by Ho (2009) \cite{ho2009}.

\subsection{Join products with 5-vertex graphs}

We first consider the 21 connected graphs on five vertices; see \hyperref[tab-5vjoin]{Table \ref{tab-5vjoin}}. We use the same graph identifiers as first given by Kle\v{s}\v{c} \cite{klesc2001_2} to describe the 5-vertex graphs. For many of the graphs, the crossing numbers of their join products with discrete graphs, paths and cycles have been determined. A list of publications where each result was first proved is displayed in \hyperref[app-j5]{Appendix \ref{app-j5}}, except for the cases $G^5_{1} + D_n$ and $G^5_{8} + D_n$ which are immediately corollaries of Kleitman \cite{kleitman1971}. It is worth noting that $n(n-1)$, a common expression in \hyperref[tab-5vjoin]{Table \ref{tab-5vjoin}}, is equal to $Z(5,n) + 2\fl[n]$.

\begin{table}[htbp]\begin{center}\caption{Crossing numbers of joins of connected 5-vertex graphs with discrete graphs, paths and cycles. The results for $D_n$ are for $n \geq 1$, the results for $P_{n-1}$ are for $n \geq 2$, and the results for $C_n$ are for $n \geq 3$. Empty cells imply that the crossing number has not yet been determined.\label{tab-5vjoin}}\smallskip\smallskip
\scalebox{0.92}{$\begin{array}{|c|c|c|c|c|}\hline
\rule{0pt}{2.3ex} i         & G^5_i                      & cr(G^5_i + D_n)    & cr(G^5_i + P_{n-1})          & cr(G^5_i + C_n)   \\
\hline 1  & \ctikz{\ngraph{5.1}}       & Z(5,n)                   & Z(5,n)                             & Z(5,n) + 1              \\
\hline 2  & \ctikz{\ngraph{5.2}}       & n(n-1)                   & \unv{n(n-1)}  & n(n-1) + 2              \\
\hline 3  & \ctikz{\ngraph{5.3}}       & Z(5,n) + \fl[n]          & Z(5,n) + \fl[n]                    & Z(5,n) + \fl[n] + 1     \\
\hline 4  & \ctikz{\ngraph{5.4}}       & Z(5,n) + \fl[n]          & Z(5,n) + \fl[n]                    & Z(5,n) + \fl[n] + 1 \\
\hline 5  & \ctikz{\ngraph{5.5}}       & Z(5,n) + \fl[n]          & Z(5,n) + \fl[n]                    & Z(5,n) + \fl[n] + 2     \\
\hline 6  & \ctikz{\ngraph{5.6}}       & n(n-1)                   & n(n-1)                             & \unv{n(n-1) + 2   }           \\
\hline 7  & \ctikz{\ngraph{5.7}}       & Z(5,n) + \fl[n]        & \unv{Z(5,n) + \fl[n] + 1  }      & Z(5,n) + \fl[n] + 2     \\
\hline 8  & \ctikz{\ngraph{5.8}}       & Z(5,n)                    & Z(5,n) + 1                         & Z(5,n) + 2              \\
\hline 9  & \ctikz{\ngraph{5.9}}       & n(n-1)                         & n(n-1)                             & \unv{n(n-1) + 2  }       \\
\hline 10 & \ctikz{\ngraph{5.10}}      & Z(5,n) + n               & Z(5,n) + n + 1                     & Z(5,n) + n + 3          \\
\hline 11 & \ctikz{\ngraph{5.11}}      & n(n-1)                   & \unv{n(n-1) + 1}                         & \unv{n(n-1) + 3 } \\
\hline 12 & \ctikz{\ngraph{5.12}}      & \unv{n(n-1)}            & n(n-1)                             &                         \\
\hline 13 & \ctikz{\ngraph{5.13}}      & Z(5,n) + \fl[n]          & Z(5,n) + \fl[n] + 1                & Z(5,n) + \fl[n] + 2     \\
\hline 14 & \ctikz{\ngraph{5.14}}      & n(n-1)                   & \unv{n(n-1) + 1 }                        & \unv{n(n-1) + 3 } \\
\hline 15 & \ctikz{\ngraph{5.15}}      & \unv{Z(5,n) + n + \fl[n] }   & Z(5,n) + n + \fl[n] + 2            & Z(5,n) + n + \fl[n] + 4 \\
\hline 16 & \ctikz{\ngraph{5.16}}      & Z(5,n) + n + \fl[n]      &  & \\
\hline 17 & \ctikz{\ngraph{5.17}}      & Z(5,n) + n               & Z(5,n) + n + 1                     &                         \\
\hline 18 & \ctikz{\ngraph{5.18}}      & Z(5,n) + n + \fl[n]  & \unv{Z(5,n) + n + \fl[n] + 2}      & Z(5,n) + n + \fl[n] + 4 \\
\hline 19 & \ctikz{\ngraph{5.19}}      & Z(5,n) + n + \fl[n]  & Z(5,n) + n + \fl[n] + 1     & Z(5,n) + n + \fl[n] + 4 \\
\hline 20 & \ctikz{\ngraph{5.20}}      & Z(5,n) + 2n              & Z(5,n) + 2n + 2                    &                         \\
\hline 21 & \ctikz{\ngraph{5.21}}      & Z(5,n) + 2n + \fl[n] + 1 & Z(5,n) + 2n + \fl[n] + 4           &                         \\
\hline\end{array}$}\end{center}\end{table}

In addition, in 2018, Su \cite{su2018_2} gave a conjecture about $cr(G^5_{21} + C_n)$:

\begin{conjecture}[Su, 2018 \cite{su2018_2}] For $n\geq 3$,
\[ cr(G^5_{21} + C_n) = Z(5,n) + 2n + \fl[n] + 7. \]\end{conjecture}

\subsubsection{Disconnected 5-vertex graphs}

In 2014, Li \cite{li2014} $\Asterisk$ considered the disconnected graph constructed by taking the union of $C_4$ and one isolated vertex:

\begin{theorem}[Li, 2014 \cite{li2014} $\Asterisk$]Let $G$ be $C_4 \cup K_1$, then the following hold:
\begin{align*}cr(G + D_n) & = Z(5,n) + \fl[n],\mbox{ for } n \geq 1, \: \Asterisk\\
cr(G + P_{n-1}) & = Z(5,n) + \fl[n] + 1,\mbox{ for } n \geq 2, \: \Asterisk\\
cr(G + C_n) & = Z(5,n) + \fl[n] + 2,\mbox{ for } n \geq 3. \: \Asterisk\end{align*}\end{theorem}

The result for $cr(G + D_n)$ was independently confirmed in 2018 by Ding and Huang \cite{dinghuang2018}.

In 2019, Sta\v{s} \cite{stas2019} considered the disconnected graph constructed by taking the union of $G^4_9$ with an isolated vertex.

\begin{theorem}[Sta\v{s}, 2019 \cite{stas2019}]Let $G$ be $G^4_9 \cup K_1$, then for $n\geq 1$, the following holds:
\[ cr(G + D_n) = Z(5,n) + \fl[n].\]\end{theorem}

Kle\v{s}\v{c} and Sta\v{s} \cite{klescstastoappear} also considered $G^4_9 \cup K_1$ as well as another disconnected graph constructed by taking the union of $G^4_7$ with an isolated vertex.

\begin{theorem}[Kle\v{s}\v{c} and Sta\v{s}, to appear \cite{klescstastoappear}]Let $G_1$ be $G^4_7 \cup K_1$ and $G_2$ be $G^4_9 \cup K_1$, then the following hold:
\begin{align*}cr(G_1 + D_n) = cr(G_2 + D_n) & = Z(5,n) + \fl[n],\mbox{ for } n \geq 1,\\
cr(G_1 + P_{n-1}) = cr(G_2 + P_{n-1}) & = Z(5,n) + \fl[n],\mbox{ for } n \geq 2,\\
cr(G_1 + C_n) = cr(G_2 + C_n) & = Z(5,n) + \fl[n] + 1,\mbox{ for } n \geq 3.\end{align*}\end{theorem}

In 2020, Sta\v{s} \cite{stas2020_4} considered the disconnected graph constructed by taking the union of $K_4$ with an isolated vertex.

\begin{theorem}[Sta\v{s}, 2020 \cite{stas2020_4}]Let $G$ be $K_4 \cup K_1$, then for $n\geq 1$, the following holds:
\[ cr(G + D_n) = Z(5,n) + \fl[n] + n.\]\end{theorem}

%
%
%

\subsection{Join products with 6-vertex graphs}

So far, join products involving 6-vertex graphs are only known for some cases. Specifically, for the graphs displayed in \hyperref[tab-6vjoin]{Table \ref{tab-6vjoin}}. A list of publications where each result was first proved is displayed in \hyperref[app-j6]{Appendix \ref{app-j6}}, except for the cases $G^6_{25} + D_n$ and $G^6_{40} + D_n$ which are immediately corollaries of Kleitman \cite{kleitman1971}. As in \hyperref[sec-cart6]{Section \ref{sec-cart6}}, we use the graph indices from \hyperref[app-g6]{Appendix \ref{app-g6}} to denote each graph for which a result has been determined.

\begin{table}[h!]\centering\vspace*{-0.3cm}\caption{Crossing numbers of joins of particular 6-vertex graphs with discrete graphs, paths and cycles. The results for $D_n$ are for $n \geq 1$, the results for $P_{n-1}$ are for $n \geq 2$, and the results for $C_n$ are for $n \geq 3$. Empty cells imply that the crossing number has not yet been determined.\label{tab-6vjoin}}\vspace{0.3cm}
\resizebox{\textwidth}{!}{$\begin{array}{|c|c|c|c|c|}\hline
\rule{0pt}{2.3ex} i & G^6_i & \hspace*{1.2cm} cr(G^6_i + D_n) \hspace*{1.2cm} & \hspace*{1.2cm} cr(G^6_i + P_{n-1}) \hspace*{1.2cm} & \hspace*{1.2cm} cr(G^6_i + C_n) \hspace*{1.2cm}\\
\hline 20 &  \ctikz{\ngraph{6.25}}      & Z(6,n) + \fl[n]       & Z(6,n) + \fl[n] + 1      &                          \\
\hline 25 &  \ctikz{\ngraph{6.25}}      & Z(6,n)                & Z(6,n)                   & Z(6,n) + 1               \\
\hline 31 &  \ctikz{\ngraph{6.31}}      & \unv{Z(6,n) + 4\fl[n]}& \unv{Z(6,n) + 4\fl[n]}   & Z(6,n) + 4\fl[n] + 3     \\
\hline 35 &  \ctikz{\ngraph{6.35}}      & Z(6,n) + 2\fl[n]      &                          &                          \\
\hline 40 &  \ctikz{\ngraph{6.40}}      & Z(6,n)                & Z(6,n) + 1               & Z(6,n) + 2               \\
\hline 41 &  \ctikz{\ngraph{6.41}}      & Z(6,n) + \fl[n]       & Z(6,n) + \fl[n] + 1      &                           \\
\hline 44 &  \ctikz{\ngraph{6.44}}      & Z(6,n) + 2\fl[n]      & \unv{Z(6,n) + 2\fl[n]}      &                          \\
\hline 45 &  \ctikz{\ngraph{6.45}}      & Z(6,n) + 2\fl[n]      &                             &                          \\
\hline 48 &  \ctikz{\ngraph{6.48}}      &                       & \unv{Z(6,n) + 4\fl[n]  }  &                          \\
\hline 49 &  \ctikz{\ngraph{6.49}}      & Z(6,n) + 2\fl[n]      & Z(6,n) + 2\fl[n]         & Z(6,n) + 2\fl[n] + 2     \\
\hline 59 &  \ctikz{\ngraph{6.59}}      & \unv{Z(6,n) + 2\fl[n] } & Z(6,n) + 2\fl[n] + 1     &                          \\
\hline 60 &  \ctikz{\ngraph{6.60}}      & Z(6,n) + \fl[n]       & Z(6,n) + \fl[n] + 1      & Z(6,n) + \fl[n] + 2      \\
\hline 61 &  \ctikz{\ngraph{6.61}}      & Z(6,n) + n            & Z(6,n) + n + 1           & Z(6,n) + n + 3           \\
\hline 66 &  \ctikz{\ngraph{6.66}}      & Z(6,n) + 2\fl[n]      &                          &                          \\
\hline 72 &  \ctikz{\ngraph{6.72}}      &                       & \unv{Z(6,n) + 4\fl[n]}         &                          \\
\hline 73 &  \ctikz{\ngraph{6.73}}      &                       & \unv{Z(6,n) + 4\fl[n]  }       &                          \\
\hline 74 &  \ctikz{\ngraph{6.74}}      & Z(6,n) + 2\fl[n]      &                          &                          \\
\hline 79 &  \ctikz{\ngraph{6.79}}      &                       & \unv{Z(6,n) + 4\fl[n]  }   &                          \\
\hline 83 &  \ctikz{\ngraph{6.83}}      & Z(6,n) + 2\fl[n]      & Z(6,n) + 2\fl[n] + 1     & Z(6,n) + 2\fl[n] + 2     \\
\hline 84 &  \ctikz{\ngraph{6.84}}      & Z(6,n) + n + \fl[n]   &                          &                          \\
\hline 85 &  \ctikz{\ngraph{6.85}}      & \unv{Z(6,n) + n }     &                          &                          \\
\hline 89 &  \ctikz{\ngraph{6.89}}      & Z(6,n) + 3\fl[n]      &                          &                          \\
\hline 93 &  \ctikz{\ngraph{6.93}}      & Z(6,n) + 2n           &                          &                          \\
\hline 94 &  \ctikz{\ngraph{6.94}}      & Z(6,n) + 2\fl[n]      & Z(6,n) + 2\fl[n] + 1     & Z(6,n) + 2\fl[n] + 3     \\
\hline
\end{array}$}\end{table}

\begin{table}[h!]Table \ref{tab-6vjoin} (continued):
\begin{center}\resizebox{\textwidth}{!}{$\begin{array}{|c|c|c|c|c|}
\hline 103 & \ctikz{\ngraph{6.103}}     & \unv{Z(6,n) + 2\fl[n] + 2n} & \unv{Z(6,n) + 2\fl[n] + 2n + 2}&                          \\
\hline 109 & \ctikz{\ngraph{6.109}}     & Z(6,n) + 2n           & Z(6,n) + 2n + 1          & Z(6,n) + 2n + 3          \\
\hline 111 & \ctikz{\ngraph{6.111}}     & \unv{Z(6,n) + n + \fl[n] }  & \unv{Z(6,n) + n + \fl[n] + 1}  & \unv{Z(6,n) + n + \fl[n] + 3 } \\
\hline 120 & \ctikz{\ngraph{6.120}}     & Z(6,n) + 3\fl[n]      & Z(6,n) + 3\fl[n] + 2     & Z(6,n) + 3\fl[n] + 4     \\
\hline 124 & \ctikz{\ngraph{6.124}}      & \unv{Z(6,n) + n + 3\fl[n]}  &                          &                          \\
\hline 125 &  \ctikz{\ngraph{6.125}}      & Z(6,n) + n + 3\fl[n]  & \unv{Z(6,n) + n + 3\fl[n] + 1} &                          \\
\hline 130 &  \ctikz{\ngraph{6.130}}      & \unv{Z(6,n) + 2n}           &                          & \hspace*{3.92cm}         \\
\hline 137 &  \ctikz{\ngraph{6.137}}     & \unv{Z(6,n) + 2n  } &                          &                          \\
\hline 152 &  \ctikz{\ngraph{6.152}}      & Z(6,n) + 3n           &                          &                          \\
\hline\end{array}$}\end{center}\end{table}

\null
\vfill
\newpage
\null
\vfill

\subsection{Other join products}

\subsubsection{Triangular snakes}

The triangular snake graph $TS_n$ is the graph with $n$ vertices, for odd $n$, defined as follows. Start by taking the path graph $P_{n-1}$ and add the edges $\{2i-1,2i+1\}$ for $i = 1, \hdots, \frac{n-1}{2}$. An example of $TS_{11}$ is displayed in \hyperref[fig-ts13]{Figure \ref{fig-ts13}}.

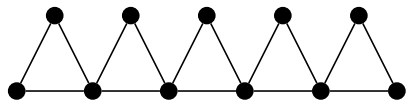
\begin{figure}[h!]\begin{center}
\begin{tikzpicture}[largegraph,scale=0.5]
\foreach \n in {0,...,5}{
  \pgfmathtruncatemacro{\m}{2*\n - 1}
  \node (a\n) at (\m,0) {};
}
\foreach \n in {0,...,4}{
  \pgfmathtruncatemacro{\m}{2*\n}
  \node (b\n) at (\m,2) {};
}
\foreach \n in {0,...,4}{
  \pgfmathtruncatemacro{\m}{\n + 1}
  \draw (a\n) -- (b\n);
  \draw (a\n) -- (a\m);
  \draw (a\m) -- (b\n);
}
\end{tikzpicture}\caption{The triangular snake $TS_{11}$.\label{fig-ts13}}\end{center}\end{figure}

Rajan et al.\ \cite{rajanetal2012} considered the join product of $TS_n$ with discrete graphs, paths, and cycles with at most six vertices:

\begin{theorem}[Rajan et al.,\ 2012 \cite{rajanetal2012}]For odd $n \geq 5$ and $m\leq6$, the following hold:
\begin{align*}cr(TS_n + D_m) & = Z(n,m) + \fl[n]\fl[m],\\
cr(TS_n + P_{m-1}) & = Z(n,m) + \fl[n]\fl[m],\\
cr(TS_n + C_m) & = Z(n,m) + \fl[n]\fl[m] + 2.\end{align*}\end{theorem}

\subsubsection{Cycles and wheels}

In 2014, Yue et al.\ \cite{yueetal2014} considered the join product of the wheel graph $W_m$ with the cycle graph $C_n$, and gave a conjecture as to its crossing number. The result for $m=3$ was already known from Kle\v{s}\v{c} (2007) \cite{klesc2007}. Stas and Valiska \cite{stasvaliska2021} subsequently confirmed the result holds for $\min\{m,n\} \leq 4$.

\begin{conjecture}[Yue et al.,\ 2014 \cite{yueetal2014}, Stas and Valiska, 2021 \cite{stasvaliska2021}]For $m,n \geq 3$, the following holds:
\[ cr(W_m + C_n) = Z(m+1,n) + \fl[m]\fl[m-1]\fl[n] + \left\lceil\frac{m}{2}\right\rceil + \left\lceil\frac{n}{2}\right\rceil + 2, \]
with the conjecture known to hold for $m = 3$ and $m = 4$.\end{conjecture}

%
%
%

\section{Other kinds of graph products}

\subsection{Strong products}

The strong product of two graphs $G$ and $H$, denoted $G \boxtimes H$, is the graph with the vertex set $V(G) \times V(H)$ and edge set $\{\big((u,v),(x,y)\big)$ : $u = x$ and $(v,y) \in E(H)$, or $v = y$ and $(u,x) \in E(G)$, or $(u,x) \in E(G)$ and $(v,y) \in E(H)$\}.  An example of $P_3 \boxtimes P_4$ is displayed in \hyperref[fig-p3xp4]{Figure \ref{fig-p3xp4}}.

\subsubsection{Two paths}

In 2013, Kle\v{s}\v{c} et al.\ \cite{klescetal2013} considered the strong product of two path graphs, $P_n \boxtimes P_m$ for $n,m \geq 2$. They first determined the crossing number of $P_n \boxtimes P_2$, and then proposed a conjecture about the crossing number of $P_n \boxtimes P_m$, which was subsequently proved by Ma \cite{ma2017} to be correct in all cases except $P_3 \boxtimes P_3$.

\begin{figure}[h!]\begin{center}
\begin{tikzpicture}[largegraph,scale=0.67]
\foreach \n in {0,...,4}{
  \node (a\n) at (2*\n,6) {};
  \node (b\n) at (2*\n,4) {};
  \node (c\n) at (2*\n,2) {};
  \node (d\n) at (2*\n,0) {};
}
\foreach \n in {0,...,3}{
  \pgfmathtruncatemacro{\m}{\n + 1}
  \draw (a\n) -- (a\m);
  \draw (b\n) -- (b\m);
  \draw (c\n) -- (c\m);
  \draw (d\n) -- (d\m);
  \draw (a\n) -- (b\m);
  \draw (a\m) -- (b\n);
  \draw (b\n) -- (c\m);
  \draw (b\m) -- (c\n);
  \draw (c\n) -- (d\m);
  \draw (c\m) -- (d\n);
}
\foreach \n in {0,...,4}{
  \draw (a\n) -- (b\n);
  \draw (b\n) -- (c\n);
  \draw (c\n) -- (d\n);
}

\end{tikzpicture}\caption{The strong product $P_3 \boxtimes P_4$.\label{fig-p3xp4}}\end{center}\end{figure}
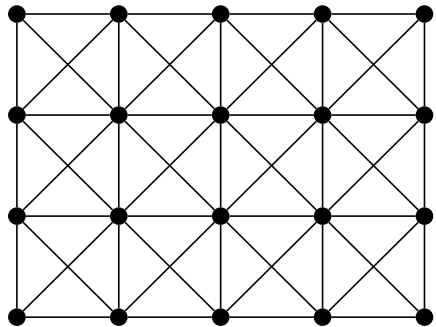

\begin{theorem}[Kle\v{s}\v{c} et al.,\ 2013 \cite{klescetal2013}]For $n \geq 2$, the following holds:
\[ cr(P_n \boxtimes P_2) = n-2. \]\end{theorem}

\begin{theorem}[Ma, 2017 \cite{ma2017}]The following holds:
\[ cr(P_n \boxtimes P_m) = \begin{cases} 4, & \mbox{ for } n=m=3,\\ mn - 4, & \mbox{ for } n > m \geq 3.\end{cases} \] \end{theorem}

\subsubsection{Paths and cycles}

In 2018, Ouyang et al.\ \cite{ouyangetal2018_2} considered the strong product of a path and a cycle:

\begin{theorem}[Ouyang et al.,\ 2018 \cite{ouyangetal2018_2}]\label{thm-strongpc}For $m \geq 1$ and $n \geq 3$, the following holds:
\[ cr(P_m \boxtimes C_n) \leq (m-1)n, \]
with equality for $m = 1$ and $m = 2$.\end{theorem}

\begin{conjecture}[Ouyang et al.,\ 2018 \cite{ouyangetal2018_2}] \hyperref[thm-strongpc]{Theorem~\ref{thm-strongpc}} holds with equality.\end{conjecture}

\section{Hypercubes, meshes, and other recursive constructions}

The graph families in this section have the common property that their number of vertices grow exponentially in terms of their parameters. Many, but not all, of these graph families arose due to their relation to network topologies in VLSI design. Leighton's pioneering work in this area continues to motivate the study of these graphs today \cite{leighton1983, leighton1984}.

\subsection{Hypercubes and related constructions}

The $n$-dimensional hypercube, $Q_n$, contains $2^n$ vertices and $n \cdot 2^{n-1}$ edges. Each vertex is labelled by a different $n$-digit binary number and an edge exists between two vertices when their labels differ in exactly one digit. An example of the $4$-cube, $Q_4$, is displayed in \hyperref[fig-4cube]{Figure \ref{fig-4cube}} in two different drawings.

\begin{figure}[h!]\begin{center}
\begin{tikzpicture}[largegraph,scale=1]
\node (0000) at (-4,-2) {};\node (0001) at (-1,-2) {};\node (0010) at (-4,1) {};\node (0011) at (-1,1) {};
\node (0100) at (-3,-1) {};\node (0101) at (-2,-1) {};\node (0110) at (-3,0) {};\node (0111) at (-2,0) {};
\node (1000) at (4,-2) {};\node (1001) at (1,-2) {};\node (1010) at (4,1) {};\node (1011) at (1,1) {};
\node (1100) at (3,-1) {};\node (1101) at (2,-1) {};\node (1110) at (3,0) {};\node (1111) at (2,0) {};

\draw (0000) -- (0001); \draw (0000) -- (0010); \draw (0000) -- (0100); \draw (0001) -- (0011);
\draw (0001) -- (0101); \draw (0010) -- (0011); \draw (0010) -- (0110); \draw (0011) -- (0111);
\draw (0100) -- (0101); \draw (0100) -- (0110); \draw (0101) -- (0111); \draw (0110) -- (0111);
\draw (1000) -- (1001); \draw (1000) -- (1010); \draw (1000) -- (1100); \draw (1001) -- (1011);
\draw (1001) -- (1101); \draw (1010) -- (1011); \draw (1010) -- (1110); \draw (1011) -- (1111);
\draw (1100) -- (1101); \draw (1100) -- (1110); \draw (1101) -- (1111); \draw (1110) -- (1111);
\draw (0111) -- (1111); \draw (0101) -- (1101); \draw (0001) -- (1001); \draw (0011) -- (1011);

\draw[bend right=90,looseness=0.25] (0000) to (1000);\draw[bend right=90,looseness=0.25] (0100) to (1100);
\draw[bend left=90,looseness=0.25] (0110) to (1110);\draw[bend left=90,looseness=0.25] (0010) to (1010);
\end{tikzpicture} \;\;\;\;\;\;\; \begin{tikzpicture}[largegraph]
\pgfmathsetmacro{\irad}{2.25*(1-4*sin(22.5)*sin(22.5))}
\foreach \n in {0,...,7}{
  \node (a\n) at (90+\n*360/8:2.25) {};
  \node (b\n) at (90+\n*360/8:\irad) {};
}
\foreach \n in {0,...,7}{
  \pgfmathtruncatemacro{\plusone}{mod(\n+1,8)}
  \pgfmathtruncatemacro{\plusthree}{mod(\n+3,8)}
  \draw (a\n) -- (a\plusone);
  \draw (b\n) -- (b\plusthree);
  \draw (a\plusone) -- (b\n);
  \draw (a\n) -- (b\plusone);
}
\end{tikzpicture}\caption{The hypercube $Q_4$ in two different drawings.\label{fig-4cube}}\end{center}\end{figure}
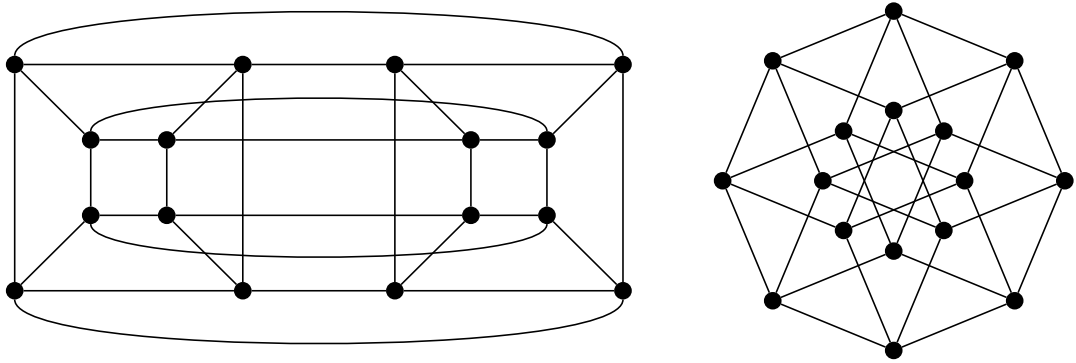

\subsubsection{$n$-cube}

The crossing number of the $n$-dimensional hypercube, often called the $n$-cube, was first considered by Eggleton and Guy \cite{eggletonguy1970} in 1970, who claimed to have discovered an upper bound. By 1973, it was known that their proof contained an error, as was discussed by Erd\H{o}s and Guy \cite{erdosguy1973}. Nonetheless, Erd\H{o}s and Guy conjectured that not only would the upper bound discovered in \cite{eggletonguy1970} be proved to be correct, but that the crossing number would meet this bound exactly. In the following years, some weaker bounds were proved.

\begin{theorem}[Madej, 1991 \cite{madej1991}]For $n\geq 1$, the following holds:
\[ cr(Q_n) \leq \frac{1}{6}4^n - n^22^{n-3} - 3 \cdot 2^{n-4} + \frac{1}{48}(-2)^n. \]\end{theorem}

\begin{theorem}[Faria and De Figueiredo, 2000 \cite{fariadefigueiredo2000}] For $n\geq 1$, the following holds:
\[ cr(Q_n) \leq \frac{165}{1024}4^n - \frac{2n^2-11n+34}{2}2^{n-2}. \]\end{theorem}

Finally, in 2008 the original upper bound from Eggleton and Guy was proved to be correct by Faria et al.\ \cite{fariaetal2008} who established that:

\begin{theorem}[Faria et al.,\ 2008 \cite{fariaetal2008}]For $n\geq 1$, the following holds:
\[ cr(Q_n) \leq \frac{5}{32}4^n - \fl[n^2+1]2^{n-2}. \]
\label{thm:Faria}\end{theorem}

It has since been claimed by Yang et al \cite{yangetal2017_2} that Faria et al. had an error in their proof as well. Indeed, Yang et al.\ \cite{yangetal2017_2} give the following updated upper bounds for the $n$-cube, which implies that the previous upper bound is not tight. At the current time, their paper has not undergone peer review and exists only on ArXiv. We are not in a position to check the following result, however we have verified that their drawing of the 7-cube, provided in their paper, is valid and does improve on the upper bound given in \hyperref[thm:Faria]{Theorem \ref{thm:Faria}} for $n = 7$. This refutes the conjecture by Erd\H{o}s and Guy that equality would hold. In particular, their drawing contains 1744 crossings, compared to the 1760 crossings suggested by \hyperref[thm:Faria]{Theorem \ref{thm:Faria}}. We have also independently discovered our own drawings of the 7-cube with fewer than 1760 crossings.

\begin{theorem}[Yang et al.,\ 2017 \cite{yangetal2017_2}]The following holds:
\[ cr(Q_n) \leq \left\{\begin{array}{ll}\frac{139}{896}4^n - \fl[n^2+1]2^{n-2} + \left(\frac{4}{7}\right)2^{3\fl[n]-n}, & \mbox{ for } 5 \leq n \leq 10,\\
\frac{26695}{172032}4^n - \fl[n^2+1]2^{n-2} - \left(\frac{n^2+2}{3}\right)2^{n-2} + \left(\frac{4}{7}\right)2^{3\fl[n]-n}, & \mbox{ for } n \geq 11.\end{array}\right. \]\end{theorem}

There have also been some lower bounds proved for the $n$-cube. Most notably, in 1993, S\'{y}kora and Vrt'o \cite{sykoravrto1993} showed that:

\begin{theorem}[S\'{y}kora and Vrt'o, 1993 \cite{sykoravrto1993}]For $n\geq 1$, the following holds:
\[ cr(Q_n) > \frac{4^n}{20} - (n+1)2^{n-2}. \]\end{theorem}

The lower bound is trivial for $n \leq 4$. Finally, $Q_1$, $Q_2$ and $Q_3$ are planar, and $Q_4$ is isomorphic to $C_4 \Box C_4$, so it is known from Dean and Richter \cite{deanrichter1995} that $cr(Q_4) = 8$. The crossing number of $Q_n$ for $n \geq 5$ has still not been determined.

\subsubsection{Locally twisted cubes}

A variation of the hypercube is the $n$-dimensional locally twisted cube $LTQ_n$, proposed by Yang et al.\ \cite{yangetal2005}. It is defined as follows: Let $LTQ_2$ be identical to $Q_2$. For $n \geq 3$, $LTQ_n$ is built from two disjoint copies of $LTQ_{n-1}$. In the first copy, augment the labeling of each node by adding 0 to the front, and in the second copy, augment the label of each node by adding 1 to the front. For each node $0x_2x_3\hdots x_n$ in the first copy, add an edge to node $1(x_2+x_n)x_3\hdots x_n$ in the second copy, where the addition is modulo 2. The result is an $n$-regular graph on $2^n$ vertices. An example of the locally twisted cube $LTQ_4$ is displayed in \hyperref[fig-ltq4]{Figure \ref{fig-ltq4}}.

\begin{figure}[h!]\begin{center}
\begin{tikzpicture}[largegraph,scale=0.8]
\node (0000) at (-4,-2) {};\node (0001) at (-1,-2) {};\node (0010) at (-4,1) {};\node (0011) at (-1,1) {};
\node (0100) at (-3,-1) {};\node (0101) at (-2,0) {};\node (0110) at (-3,0) {};\node (0111) at (-2,-1) {};
\node (1000) at (4,-2) {};\node (1001) at (1,-2) {};\node (1010) at (4,1) {};\node (1011) at (1,1) {};
\node (1100) at (3,-1) {};\node (1101) at (2,0) {};\node (1110) at (3,0) {};\node (1111) at (2,-1) {};

\draw (0000) -- (0001);\draw (0000) -- (0010);\draw (0000) -- (0100);\draw (0001) -- (0011);
\draw (0001) -- (0111);\draw (0001) -- (1101);\draw (0010) -- (0011);\draw (0010) -- (0110);
\draw (0011) -- (0101);\draw (0011) -- (1111);\draw (0100) -- (0101);\draw (0100) -- (0110);
\draw (0101) -- (0111);\draw (0110) -- (0111);\draw (1000) -- (1001);\draw (1000) -- (1010);
\draw (1000) -- (1100);\draw (1001) -- (1011);\draw (1001) -- (1111);\draw (1001) -- (0101);
\draw (1010) -- (1011);\draw (1010) -- (1110);\draw (1011) -- (1101);\draw (1011) -- (0111);
\draw (1100) -- (1101);\draw (1100) -- (1110);\draw (1101) -- (1111);\draw (1110) -- (1111);
\draw[bend right=90,looseness=0.25] (0000) to (1000);\draw[bend right=90,looseness=0.25] (0100) to (1100);
\draw[bend left=90,looseness=0.25] (0110) to (1110);\draw[bend left=90,looseness=0.25] (0010) to (1010);
\end{tikzpicture}\caption{The locally twisted cube $LTQ_4$.\label{fig-ltq4}}\end{center}\end{figure}
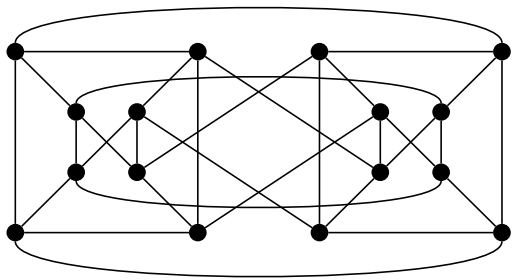

Wang et al.\ \cite{wangetal2017} determined the following bounds for $cr(LTQ_n)$:

\begin{theorem}[Wang et al.,\ 2017 \cite{wangetal2017}]For $n \geq 6$, the following bounds hold:
\[ \frac{4^n}{20} - (n^2+1)2^{n-1} < cr(LTQ_n) \leq \frac{265}{6} 4^{n-4} - \left(n^2 + \frac{15 + (-1)^{n-1}}{6}\right)2^{n-3}. \]\end{theorem}

\subsubsection{Folded hypercube}

The folded hypercube $FQ_n$ is obtained by taking $Q_n$ and adding in edges between all pairs of vertices with complementary labels (that is, their labels differ in all digits). The resulting graph is $(n+1)$-regular. An example of the folded hypercube $FQ_4$ is displayed in \hyperref[fig-fq4]{Figure \ref{fig-fq4}}.

\begin{figure}[h!]\begin{center}
\begin{tikzpicture}[largegraph,scale=0.8]
\node (0000) at (-4,-2) {};\node (0001) at (-1,-2) {};\node (0010) at (-4,1) {};\node (0011) at (-1,1) {};
\node (0100) at (-3,-1) {};\node (0101) at (-2,-1) {};\node (0110) at (-3,0) {};\node (0111) at (-2,0) {};
\node (1000) at (4,-2) {};\node (1001) at (1,-2) {};\node (1010) at (4,1) {};\node (1011) at (1,1) {};
\node (1100) at (3,-1) {};\node (1101) at (2,-1) {};\node (1110) at (3,0) {};\node (1111) at (2,0) {};

\draw (0000) -- (0001); \draw (0000) -- (0010); \draw (0000) -- (0100); \draw (0001) -- (0011);
\draw (0001) -- (0101); \draw (0010) -- (0011); \draw (0010) -- (0110); \draw (0011) -- (0111);
\draw (0100) -- (0101); \draw (0100) -- (0110); \draw (0101) -- (0111); \draw (0110) -- (0111);
\draw (1000) -- (1001); \draw (1000) -- (1010); \draw (1000) -- (1100); \draw (1001) -- (1011);
\draw (1001) -- (1101); \draw (1010) -- (1011); \draw (1010) -- (1110); \draw (1011) -- (1111);
\draw (1100) -- (1101); \draw (1100) -- (1110); \draw (1101) -- (1111); \draw (1110) -- (1111);
\draw (0111) -- (1111); \draw (0101) -- (1101); \draw (0001) -- (1001); \draw (0011) -- (1011);
\draw (0000) -- (1111); \draw (0001) -- (1110); \draw (0010) -- (1101); \draw (0011) -- (1100);
\draw (0100) -- (1011); \draw (0101) -- (1010); \draw (0110) -- (1001); \draw (0111) -- (1000);

\draw[bend right=90,looseness=0.25] (0000) to (1000);\draw[bend right=90,looseness=0.25] (0100) to (1100);
\draw[bend left=90,looseness=0.25] (0110) to (1110);\draw[bend left=90,looseness=0.25] (0010) to (1010);
\end{tikzpicture}\caption{The folded hypercube $FQ_4$.\label{fig-fq4}}\end{center}\end{figure}
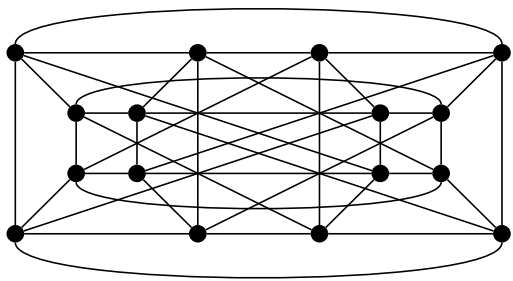

Wang et al.\ \cite{wangetal2015} provided lower and upper bounds for the crossing numbers of these graphs:


\begin{theorem}[Wang et al.,\ 2015 \cite{wangetal2015}]For $n\geq 1$, the following bounds hold:
\[ \frac{1}{20} 4^n (1- (\pi/2(2\lceil n/2\rceil +1))^{-1/2})^{-2} < cr(FQ_n) \leq \frac{11}{32}4^n-(n^2+3n)2^{n-3}. \]\end{theorem}

\subsubsection{Augmented cube}

The augmented cube $AQ_n$ is defined as follows. Let $AQ_1$ be isomorphic to $K_2$, with vertices labelled $0$ and $1$. $AQ_n$ is generated by taking two copies of $AQ_{n-1}$, and prepending 0 to the labels of the first copy and 1 to the labels of the second copy. That is, their new labels in the first copy all begin with 0, and in the second copy all begin with 1. Two vertices are joined by a new edge if and only if their labels differ in only the first position (e.g.\ identical vertices from each copy), or if the labels differ in all positions. The augmented cube was first proposed by Choudum and Sunitha (2002) \cite{choudumsunitha2002}, who also provided a non-recursive definition, as follows. $AQ_n$ is the graph containing vertices labelled with $n$-digit binary numbers, and any two vertices are connected by an edge if and only if there exists an $l$, $1 \leq l \leq n$, such that either (1) the two labels are different in position $l$, and identical in all other positions, or (2) the two labels are identical in the first $l-1$ positions and then different in all subsequent positions. The resulting graph is $(2n-1)$-regular. An example of the augmented cube $AQ_3$ is displayed in \hyperref[fig-aq3]{Figure \ref{fig-aq3}}.

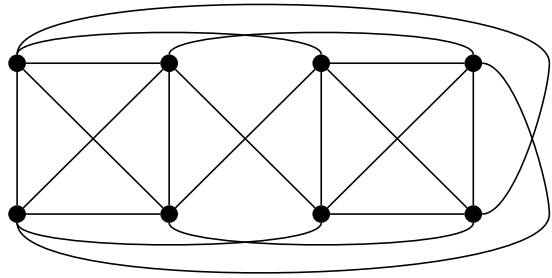
\begin{figure}[h!]\begin{center}
\begin{tikzpicture}[largegraph,scale=0.5]
\node (000) at (-6,2) {};\node (001) at (-2,2) {};\node (010) at (-6,-2) {};\node (011) at (-2,-2) {};
\node (100) at (2,2) {};\node (101) at (6,2) {};\node (110) at (2,-2) {};\node (111) at (6,-2) {};

\draw (000) -- (001); \draw (000) -- (010); \draw (000) -- (011); \draw (001) -- (010); \draw (001) -- (011); \draw (010) -- (011);
\draw (100) -- (101); \draw (100) -- (110); \draw (100) -- (111); \draw (101) -- (110); \draw (101) -- (111); \draw (110) -- (111);
\draw (001) -- (110); \draw (011) -- (100);

\draw[bend left=90,looseness=0.25] (000) to (100);
\draw[bend right=90,looseness=0.25] (010) to (110);
\draw[bend left=90,looseness=0.25] (001) to (101);
\draw[bend right=90,looseness=0.25] (011) to (111);
\draw[bend left=90,looseness=0.35] (000) to (8,2); \draw[out=270,in=0,looseness=0.5] (8,2) to (111);
\draw[bend right=90,looseness=0.35] (010) to (8,-2); \draw[out=90,in=0,looseness=0.5] (8,-2) to (101);
\end{tikzpicture}\caption{The augmented cube $AQ_3$.\label{fig-aq3}}\end{center}\end{figure}

In 2013, Wang et al.\ \cite{wangetal2013} investigated the crossing number of the augmented cube and discovered lower and upper bounds:

\begin{theorem}[Wang et al.,\ 2013 \cite{wangetal2013}]The following hold: $cr(AQ_3) = 4$,\\$cr(AQ_4) \leq 46$, $cr(AQ_5) \leq 328$, $cr(AQ_6) \leq 1848$, $cr(AQ_7) \leq 9112$, and for $n \geq 8$,
\[ \frac{4^n}{5(1 + 2^{2-n})^2} - (4n^2 + 4n + \frac{17}{5})2^{n-1} < cr(AQ_n) < \frac{13}{16}4^n - (2n^2 + \frac{7}{2}n - 6)2^{n-2}. \]
The lower bound is valid for all $n$, but is only meaningful for $n \geq 11$.\end{theorem}

\subsubsection{Cube connected cycle graphs}

The cube connected cycle graph $CCC_n$ is obtained from $Q_n$ by replacing vertices of $Q_n$ with cycles of length $n$, as described in S\'{y}kora and Vrt'o (1993) \cite{sykoravrto1993}. An example of the cube connected cycle graph $CCC_4$ is displayed in \hyperref[fig-ccc4]{Figure \ref{fig-ccc4}}.

\begin{figure}[h!]\begin{center}
\begin{tikzpicture}[largegraph,scale=0.03]
\node (a10) at (100-7,0+7) {};\node (a20) at (100+7,0+7) {};\node (a30) at (100+7,0-7) {};\node (a40) at (100-7,0-7) {};
\node (a11) at (71-7,71+7) {};\node (a21) at (71+7,71+7) {};\node (a31) at (71+7,71-7) {};\node (a41) at (71-7,71-7) {};
\node (a12) at (0-7,100-7) {};\node (a22) at (0-7,100+7) {};\node (a32) at (0+7,100+7) {};\node (a42) at (0+7,100-7) {};
\node (a13) at (-71-7,71-7) {};\node (a23) at (-71-7,71+7) {};\node (a33) at (-71+7,71+7) {};\node (a43) at (-71+7,71-7) {};
\node (a14) at (-100+7,0-7) {};\node (a24) at (-100-7,0-7) {};\node (a34) at (-100-7,0+7) {};\node (a44) at (-100+7,0+7) {};
\node (a15) at (-71+7,-71-7) {};\node (a25) at (-71-7,-71-7) {};\node (a35) at (-71-7,-71+7) {};\node (a45) at (-71+7,-71+7) {};
\node (a16) at (0+7,-100+7) {};\node (a26) at (0+7,-100-7) {};\node (a36) at (0-7,-100-7) {};\node (a46) at (0-7,-100+7) {};
\node (a17) at (71-7,-71+7) {};\node (a27) at (71+7,-71+7) {};\node (a37) at (71+7,-71-7) {};\node (a47) at (71-7,-71-7) {};

\node (b10) at (52-7,0+7) {};\node (b20) at (52+7,0+7) {};\node (b30) at (52+7,0-7) {};\node (b40) at (52-7,0-7) {};
\node (b11) at (37-7,37+7) {};\node (b21) at (37+7,37+7) {};\node (b31) at (37+7,37-7) {};\node (b41) at (37-7,37-7) {};
\node (b12) at (0-7,52-7) {};\node (b22) at (0-7,52+7) {};\node (b32) at (0+7,52+7) {};\node (b42) at (0+7,52-7) {};
\node (b13) at (-37-7,37-7) {};\node (b23) at (-37-7,37+7) {};\node (b33) at (-37+7,37+7) {};\node (b43) at (-37+7,37-7) {};
\node (b14) at (-52+7,0-7) {};\node (b24) at (-52-7,0-7) {};\node (b34) at (-52-7,0+7) {};\node (b44) at (-52+7,0+7) {};
\node (b15) at (-37+7,-37-7) {};\node (b25) at (-37-7,-37-7) {};\node (b35) at (-37-7,-37+7) {};\node (b45) at (-37+7,-37+7) {};
\node (b16) at (0+7,-52+7) {};\node (b26) at (0+7,-52-7) {};\node (b36) at (0-7,-52-7) {};\node (b46) at (0-7,-52+7) {};
\node (b17) at (37-7,-37+7) {};\node (b27) at (37+7,-37+7) {};\node (b37) at (37+7,-37-7) {};\node (b47) at (37-7,-37-7) {};

\foreach \n in {0,...,7}{
  \pgfmathtruncatemacro{\plusone}{mod(\n+1,8)}
  \pgfmathtruncatemacro{\plusthree}{mod(\n+3,8)}
  \draw (a2\n) -- (a3\plusone);
  \draw (b1\n) -- (b4\plusthree);
  \draw (a4\plusone) -- (b2\n);
  \draw (a1\n) -- (b3\plusone);

  \draw (a1\n) -- (a2\n);
  \draw (a2\n) -- (a3\n);
  \draw (a3\n) -- (a4\n);
  \draw (a4\n) -- (a1\n);
  \draw (b1\n) -- (b2\n);
  \draw (b2\n) -- (b3\n);
  \draw (b3\n) -- (b4\n);
  \draw (b4\n) -- (b1\n);
}
\end{tikzpicture}\caption{The cube connected cycle $CCC_4$.\label{fig-ccc4}}\end{center}\end{figure}
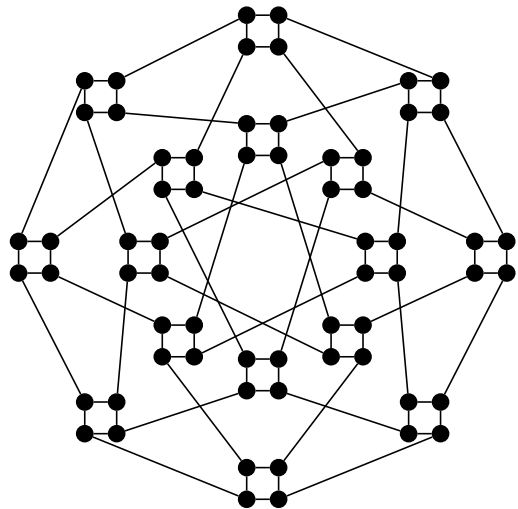

In \cite{sykoravrto1993} S\'{y}kora and Vrt'o establish lower and upper bounds for $cr(CCC_n)$:

\begin{theorem}[S\'{y}kora and Vrt'o, 1993 \cite{sykoravrto1993}]For $n\geq 1$, the following bounds hold:
\[ \frac{4^n}{20} - 3(n+1)2^{n-2} < cr(CCC_n) < \frac{4^n}{6} + 3n^22^{n-3}. \]\end{theorem}

\subsection{Complete mesh of trees}\label{sec-cmot}

The complete (2-dimensional) mesh of trees $M_n$ is defined for any $n$ which is a power of 2, as follows: start with a square $n \times n$ grid. Identify the $n$ vertices in each row and each column with the $n$ leaves of a full, complete, balanced binary tree. What results is a graph on $3n^2 - 2n$ vertices and $4n^2 - 4n$ edges. For a detailed description, see \cite{cimikowski1998}. It is easy to check that $M_2$ is just the cycle on eight vertices, and is hence planar. An example of $M_4$ is displayed in \hyperref[fig-mesh]{Figure \ref{fig-mesh}}.

\begin{figure}[h!]\begin{center}
\begin{tikzpicture}[largegraph]
\foreach \i in {0,...,3}{
  \foreach \j in {0,...,3}{
    \node (\i@\j) at (1.5*\j,1.5*\i) {};
  }
  \node (\i@hroot) at (1.5*1+0.75,1.5*\i+0.5) {};
  \node (\i@hm1) at   (1.5*0+0.75,1.5*\i+0.25) {};
  \node (\i@hm2) at   (1.5*2+0.75,1.5*\i+0.25) {};
  \draw (\i@hm1) -- (\i@hroot) -- (\i@hm2);
  \draw (\i@0) -- (\i@hm1) -- (\i@1);
  \draw (\i@2) -- (\i@hm2) -- (\i@3);
}
\foreach \j in {0,...,3}{
  \node (\j@vroot) at (1.5*\j+0.5,1.5*1+0.75) {};
  \node (\j@vm1) at   (1.5*\j+0.25,1.5*0+0.75) {};
  \node (\j@vm2) at   (1.5*\j+0.25,1.5*2+0.75) {};
  \draw (\j@vm1) -- (\j@vroot) -- (\j@vm2);
  \draw (0@\j) -- (\j@vm1) -- (1@\j);
  \draw (2@\j) -- (\j@vm2) -- (3@\j);
}
\end{tikzpicture}\caption{The complete mesh of trees $M_4$.\label{fig-mesh}}\end{center}\end{figure}
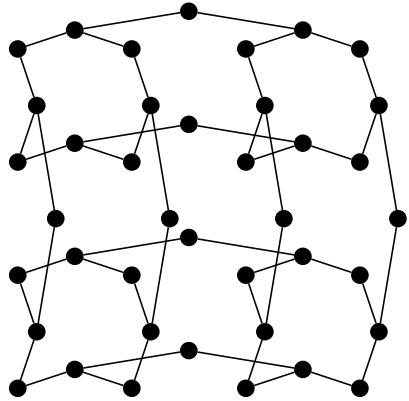

The earliest result on the crossing number of $M_n$ is due to Leighton \cite{leighton1984}, who in 1984 derived a lower bound. The bound is only interesting in an asymptotic sense since it holds trivially for $n \leq 2^{121}$.

\begin{theorem}[Leighton, 1984 \cite{leighton1984}]For $n \geq 2$ and a power of 2, the following holds:
\[ cr(M_n) \geq \frac{n^2\log_2(n) - 121n^2 + 121n}{40}. \]\end{theorem}

In 1996, Cimikowski \cite{cimikowski1996} claimed to have determined an upper bound for the crossing number of $M_n$. He later pointed out there had been an error in his proof, and in 1998 presented an amended upper bound which he conjectured would hold with equality, as well as a practical lower bound, which coincides with the upper bound for the case $n = 4$:

\begin{theorem}[Cimikowski, 1998 \cite{cimikowski1998}]For $n \geq 4$ and a power of 2, the following bounds hold:
\[ \frac{n^2}{4} \leq cr(M_n) \leq n^2\left(\frac{3\log_2(n) - 5}{8}\right) + \frac{n}{2}. \]\label{thm-mesh}\end{theorem}

\begin{conjecture}[Cimikowski, 1998 \cite{cimikowski1998}]The upper bound in \hyperref[thm-mesh]{Theorem~\ref{thm-mesh}}\\holds with equality.\end{conjecture}

In 2003, Cimikowski and Vrt'o \cite{cimikowskivrto2003} published a paper which included the results from \cite{cimikowski1998}, and also gave an alternative lower bound, which is superior to their previous lower bound for $n \geq 512$:

\begin{theorem}[Cimikowski and Vrt'o, 2003 \cite{cimikowskivrto2003}]For $n \geq 2$ and a power of 2, the following holds:
\[ cr(M_n) \geq \frac{5n^2\log_2(n) - 44n^2}{80}. \]\end{theorem}

The crossing number for $M_n$ is only known for $n = 1, 2, 4$. Specifically, $cr(M_1) = cr(M_2) = 0$, and $cr(M_4) = 4$.

\subsection{Butterfly graphs}

The Butterfly graph $BF(r)$ is the graph with $(r+1)2^r$ vertices and $r \cdot 2^{r+1}$ edges, defined as follows: The vertices are labelled $\langle w,i \rangle$ where $i = 0, \hdots, r$ and $w$ is an $r$-bit binary number. Two vertices $\langle w,i \rangle$ and $\langle w',i' \rangle$ are adjacent if and only if $i' = i+1$ and either $w = w'$, or $w$ and $w'$ differ in precisely the $i$-th bit. An example of $BF(3)$ is displayed in \hyperref[fig-BF3]{Figure \ref{fig-BF3}}.

\begin{figure}[h!]\begin{center}
\begin{tikzpicture}[largegraph]
\foreach \n in {0,...,7}{
  \node (a\n) at (\n*1.5,1.5) {};
  \node (b\n) at (\n*1.5,0.5) {};
  \node (c\n) at (\n*1.5,-0.5) {};
  \node (d\n) at (\n*1.5,-1.5) {};
}
\foreach \n in {0,...,7}{
  \draw (a\n) -- (b\n);
  \draw (b\n) -- (c\n);
  \draw (c\n) -- (d\n);
}
\foreach \n in {0,...,3}{
  \pgfmathtruncatemacro{\nextn}{\n+4}
  \draw (a\n) -- (b\nextn);
  \draw (b\n) -- (a\nextn);
}
\foreach \n in {0,1,4,5}{
  \pgfmathtruncatemacro{\nextn}{\n+2}
  \draw(b\n) -- (c\nextn);
  \draw(c\n) -- (b\nextn);
}
\foreach \n in {0,2,4,6}{
  \pgfmathtruncatemacro{\nextn}{\n+1}
  \draw(c\n) -- (d\nextn);
  \draw(d\n) -- (c\nextn);
}
\end{tikzpicture}\caption{The Butterfly graph $BF(3)$.\label{fig-BF3}}\end{center}\end{figure}
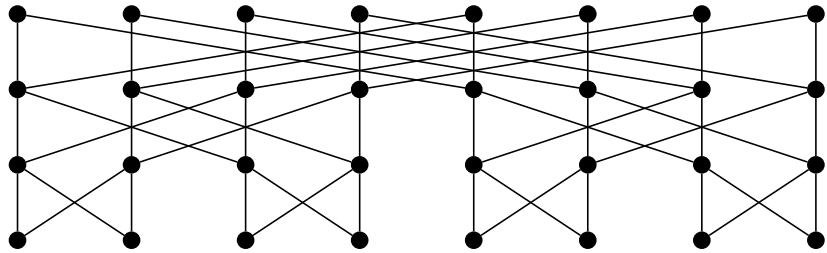

Cimikowski \cite{cimikowski1996} considered the Butterfly graphs in 1996 and proved the following upper bound for their crossing numbers:

\begin{theorem}[Cimikowski, 1996 \cite{cimikowski1996}]For $r\geq 1$, the following holds:
\[ cr(BF(r)) \leq \frac{3}{2}4^r - 3 \cdot 2^r - r \cdot 2^r + 1. \]\end{theorem}

In 2013, Manuel et al.\ \cite{manueletal2013} improved this upper bound, and also provided a lower bound:

\begin{theorem}[Manuel et al.,\ 2013 \cite{manueletal2013}]For $r\geq 3$, the following bounds hold:
\[ \frac{1}{59}4^r - r \cdot 2^r + 2^{r-1} \leq cr(BF(r)) \leq \frac{1}{4}4^r - r \cdot 2^{r-1}, \]
and $cr(BF(3)) = 4$.\end{theorem}

\subsubsection{Wrapped butterfly graphs}

The Wrapped Butterfly graph, denoted $WBF(r)$, is derived from the butterfly network $BF(r)$ by merging the first and last rows into a single row; that is, merging vertex $\langle w,0 \rangle$ with vertex $\langle w,r \rangle$ for all $w$ \cite{cimikowski1996}. Then $WBF(r)$ has $r \cdot 2^r$ vertices and $r \cdot 2^{r+1}$ edges. An example of the Wrapped Butterfly graph $WBF(3)$ is displayed in \hyperref[fig-WBF3]{Figure \ref{fig-WBF3}}.

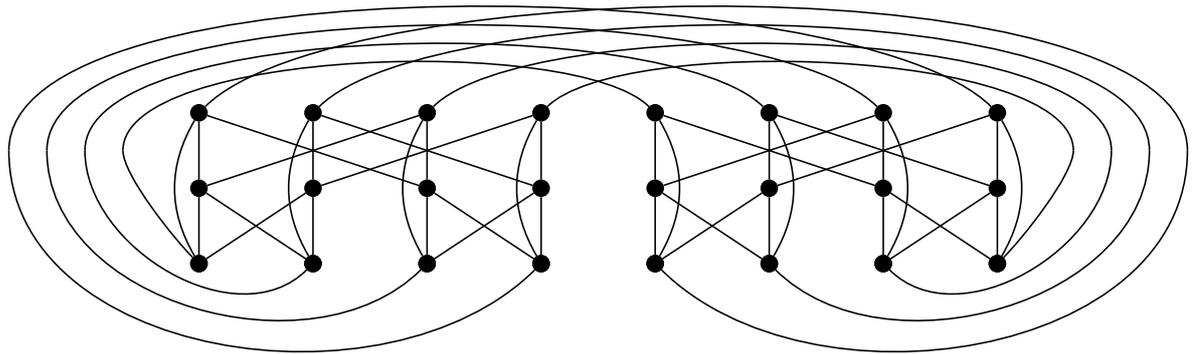
\begin{figure}[h!]\begin{center}
\begin{tikzpicture}[largegraph,scale=0.95]
\foreach \n in {0,...,7}{
  \node (a\n) at (\n*1.5,-1.5) {};
  \node (b\n) at (\n*1.5,0.5) {};
  \node (c\n) at (\n*1.5,-0.5) {};
  \node (d\n) at (\n*1.5,-1.5) {};
}
\foreach \n in {0,...,7}{
  \draw (b\n) -- (c\n);
  \draw (c\n) -- (d\n);
}
\foreach \n in {0,...,3}{
\draw [bend right, looseness=1] (b\n) to (d\n);
}
\foreach \n in {4,...,7}{
\draw [bend left, looseness=1] (b\n) to (d\n);
}

\draw [out=135,in=270,looseness=0.5] (a0) to (-1,0); \draw[out=90,in=135,looseness=0.5] (-1,0) to (b4);
\draw [out=225,in=270,looseness=1] (a1) to (-1.5,0); \draw[out=90,in=135,looseness=0.5] (-1.5,0) to (b5);
\draw [out=225,in=270,looseness=1] (a2) to (-2,0); \draw[out=90,in=135,looseness=0.5] (-2,0) to (b6);
\draw [out=225,in=270,looseness=1] (a3) to (-2.5,0); \draw[out=90,in=135,looseness=0.5] (-2.5,0) to (b7);

\draw [out=45,in=270,looseness=0.5] (a7) to (11.5,0); \draw[out=90,in=45,looseness=0.5] (11.5,0) to (b3);
\draw [out=315,in=270,looseness=1] (a6) to (12,0); \draw[out=90,in=45,looseness=0.5] (12,0) to (b2);
\draw [out=315,in=270,looseness=1] (a5) to (12.5,0); \draw[out=90,in=45,looseness=0.5] (12.5,0) to (b1);
\draw [out=315,in=270,looseness=1] (a4) to (13,0); \draw[out=90,in=45,looseness=0.5] (13,0) to (b0);

\foreach \n in {0,1,4,5}{
  \pgfmathtruncatemacro{\nextn}{\n+2}
  \draw(b\n) -- (c\nextn);
  \draw(c\n) -- (b\nextn);
}
\foreach \n in {0,2,4,6}{
  \pgfmathtruncatemacro{\nextn}{\n+1}
  \draw(c\n) -- (d\nextn);
  \draw(d\n) -- (c\nextn);
}
\end{tikzpicture}\caption{The Wrapped Butterfly graph $WBF(3)$.\label{fig-WBF3}}\end{center}\end{figure}

Cimikowski \cite{cimikowski1996} gave an upper bound for the crossing number of $WBF(r)$:

\begin{theorem}[Cimikowski, 1996 \cite{cimikowski1996}] For $r \geq 1$, the following holds:
\[ cr(WBF(r)) \leq \frac{3}{2}4^r - 3 \cdot 2^r - r \cdot 2^r. \]\end{theorem}

\subsubsection{Benes networks}

The Benes network $B(r)$ is formed by taking two copies of the Butterfly network $BF_r$. For each $w$, vertex $\langle w,r \rangle$ from the first copy is merged with $\langle w,r \rangle$ from the second copy \cite{benes1964}. As an example, $B(3)$ is displayed in \hyperref[fig-benes]{Figure \ref{fig-benes}}.

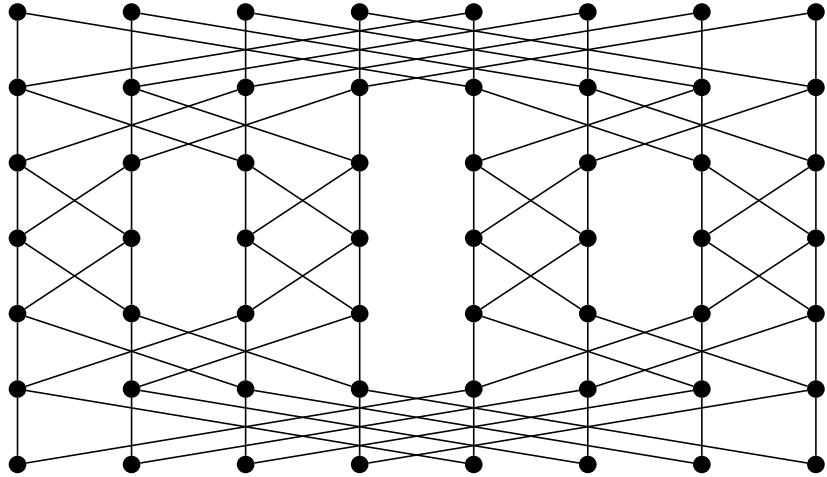
\begin{figure}[h!]\begin{center}
\begin{tikzpicture}[largegraph]
\foreach \n in {0,...,7}{
  \node (a\n) at (\n*1.5,3) {};
  \node (b\n) at (\n*1.5,2) {};
  \node (c\n) at (\n*1.5,1) {};
  \node (d\n) at (\n*1.5,0) {};
  \node (e\n) at (\n*1.5,-1) {};
  \node (f\n) at (\n*1.5,-2) {};
  \node (g\n) at (\n*1.5,-3) {};
}
\foreach \n in {0,...,7}{
  \draw (a\n) -- (b\n);
  \draw (b\n) -- (c\n);
  \draw (c\n) -- (d\n);
  \draw (d\n) -- (e\n);
  \draw (e\n) -- (f\n);
  \draw (f\n) -- (g\n);
}
\foreach \n in {0,...,3}{
  \pgfmathtruncatemacro{\nextn}{\n+4}
  \draw (a\n) -- (b\nextn);
  \draw (b\n) -- (a\nextn);
  \draw (g\n) -- (f\nextn);
  \draw (f\n) -- (g\nextn);
}
\foreach \n in {0,1,4,5}{
  \pgfmathtruncatemacro{\nextn}{\n+2}
  \draw(b\n) -- (c\nextn);
  \draw(c\n) -- (b\nextn);
  \draw(f\n) -- (e\nextn);
  \draw(e\n) -- (f\nextn);
}
\foreach \n in {0,2,4,6}{
  \pgfmathtruncatemacro{\nextn}{\n+1}
  \draw(c\n) -- (d\nextn);
  \draw(d\n) -- (c\nextn);
  \draw(e\n) -- (d\nextn);
  \draw(d\n) -- (e\nextn);
}
\end{tikzpicture}\caption{The Benes network $B(3)$.\label{fig-benes}}\end{center}\end{figure}

In 1996, Cimikowski investigated the Benes network, and determined an upper bound for its crossing number:

\begin{theorem}[Cimikowski, 1996 \cite{cimikowski1996}]For $r \geq 1$, the following holds:
\[ cr(B(r)) \leq 3 \cdot 4^r - 5 \cdot 2^r - 2r \cdot 2^r + 2. \]\end{theorem}

\subsection{Generalized fat trees}

Generalized fat trees, $GFT(h,m,w)$ were introduced in 1995 by Ohring et al.\ \cite{ohringetal1995} as a new network topology. They are a three parameter family of graphs with $\sum_{i=0}^h m^i w^{h-i}$ vertices and $\sum_{i=0}^{h-1} m^{h-i} w^{i+1}$ edges. The definition of generalized fat trees is complicated, and so we refer the interested reader to \cite{ohringetal1995} for the full definition.

%

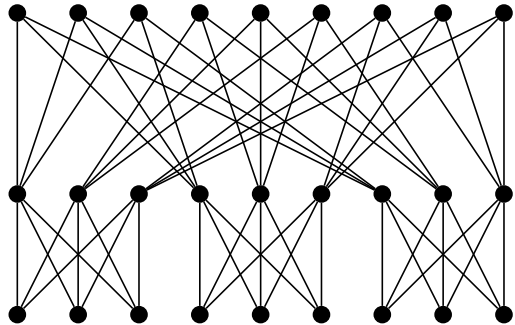
\begin{figure}[h!]\begin{center}
\begin{tikzpicture}[largegraph,scale=0.8]
\foreach \n in {0,...,8}{
  \foreach \l in {0,1,2}{
    \pgfmathtruncatemacro{\x}{-4+\n}
    \pgfmathtruncatemacro{\y}{\l+2^\l}
    \node (\l\n) at (\x,\y) {};
  }
}
\foreach \n in {0,1,2}{
  \foreach \l in {0,1,2}{
    \foreach \m in {0,1,2}{
      \pgfmathtruncatemacro{\nextn}{\n+3*\m}
      \pgfmathtruncatemacro{\nextl}{\l+3*\m}
      \pgfmathtruncatemacro{\nextm}{\n+3*\l}
      \draw (0\nextn) -- (1\nextl);
      \draw (1\nextl) -- (2\nextm);
    }
  }
}
\end{tikzpicture}\caption{The generalized fat tree $GFT(2,3,3)$.\label{fig-gft233}}\end{center}\end{figure}

In 2011, Rajan et al.\ \cite{rajanetal2011} considered the crossing number of generalized fat trees, and determined an upper bound for the special case $GFT(h,3,3)$, as well as a conjecture for an upper bound for the general case:

\begin{theorem}[Rajan et al.,\ 2011 \cite{rajanetal2011}]For $h \geq 1$, the following holds:
\[ cr(GFT(h,3,3)) \leq 3^h + \frac{5}{4}3^{2h} - \frac{1}{4}3^{h+2} - \frac{h}{2}3^{h+1}. \]\end{theorem}

\begin{conjecture}[Rajan et al.,\ 2011 \cite{rajanetal2011}]For $h \geq 2$ and $m,w \geq 1$, the following hold:
\[ cr(GFT(h,m,w)) \leq \left\{\begin{array}{lll}\frac{w^h m(m-1)}{4}\Big(w^{h+2}\frac{1 - \left(\frac{m}{w^2}\right)^h}{w^2-m} - h\Big), & & \mbox{if } m = w,\\ \frac{w^h m(m-1)}{4}\Big(w^{h}h - w\frac{1 - \left(\frac{m}{w}\right)^h}{w-m}\Big), & & \mbox{if } m = w^2,\\ \frac{w^{h+1} m(m-1)}{4} \Big(w^{h+1}\frac{1 - \left(\frac{m}{w^2}\right)^h}{w^2-m} - \frac{1 - \left(\frac{m}{w}\right)^h}{w-m}\Big), & & \mbox{otherwise}.\end{array}\right. \]\end{conjecture}

\subsection{Sierpi\'{n}ski and Sierpi\'{n}ski-like graphs}

In 2005, Klav\v{z}ar and Mohar \cite{klavzarmohar2005} considered Sierpi\'{n}ski graphs, and two of their regularizations. The Sierpi\'{n}ski graph $S(n,k)$ for $n,k \geq 1$ is defined as follows: The vertex set is $\{1,\hdots,k\}^n$, and two different vertices $u = (u_1,\hdots,u_n)$ and $v = (v_1,\hdots,v_n)$ are adjacent if and only if there exists an integer $h \in \{1,\hdots,n\}$, such that

\begin{itemize}\item[(i)] $u_t = v_t$ for $t = 1, \hdots, h-1$;
\item[(ii)] $u_h \neq v_h$; and
\item[(iii)] $u_t = v_h$ and $v_t = u_h$ for $t = h+1, \hdots, n.$\end{itemize}

In this construction, any vertex of the form $(i,i,\hdots,i)$ is called an extreme vertex. From $S(n,k)$, two regularizations can be defined. First, $S^+(n,k)$ is obtained from $S(n,k)$ by adding a single vertex $w$ which is connected to all extreme vertices in $S(n,k)$. Second, $S^{++}(n,k)$ is obtained by taking the union of $k+1$ copies of $S(n-1,k)$, and then connecting the extreme vertices as a $K_{k+1}$. Both regularizations produce a regular graph of degree $k$. Examples of $S(2,4)$, $S^+(2,4)$ and $S^{++}(2,4)$ are displayed in \hyperref[fig-sierpinksi]{Figure \ref{fig-sierpinksi}}. A full definition is given in Klav\v{z}ar and Mohar \cite{klavzarmohar2005}.

\begin{figure}[h!]\begin{center}\begin{tikzpicture}[largegraph,scale=0.7]
\foreach \n in {0,...,3}{
  \foreach \m in {0,...,3}{
    \node (\n\m) at (\m,\n) {};
  }
}
\draw (00) -- (01);\draw (01) -- (02);\draw (02) -- (03);\draw (10) -- (11);\draw (12) -- (13);\draw (20) -- (21);\draw (22) -- (23);\draw (30) -- (31);\draw (31) -- (32);\draw (32) -- (33);
\draw (00) -- (10);\draw (10) -- (20);\draw (20) -- (30);\draw (01) -- (11);\draw (21) -- (31);\draw (02) -- (12);\draw (22) -- (32);\draw (03) -- (13);\draw (13) -- (23);\draw (23) -- (33);
\draw (00) -- (11);\draw (10) -- (01);\draw (02) -- (13);\draw (12) -- (03);\draw (20) -- (31);\draw (30) -- (21);\draw (22) -- (33);\draw (32) -- (23);\draw (11) -- (22);\draw (21) -- (12);
\end{tikzpicture}$\;\;\;\;\;\;\;\;$
\begin{tikzpicture}[largegraph,scale=0.7]
\foreach \n in {0,...,3}{
  \foreach \m in {0,...,3}{
    \node (\n\m) at (\m,\n) {};
  }
}
\node (v) at (1.5,5) {};
\draw (00) -- (01);\draw (01) -- (02);\draw (02) -- (03);\draw (10) -- (11);\draw (12) -- (13);\draw (20) -- (21);\draw (22) -- (23);\draw (30) -- (31);\draw (31) -- (32);\draw (32) -- (33);
\draw (00) -- (10);\draw (10) -- (20);\draw (20) -- (30);\draw (01) -- (11);\draw (21) -- (31);\draw (02) -- (12);\draw (22) -- (32);\draw (03) -- (13);\draw (13) -- (23);\draw (23) -- (33);
\draw (00) -- (11);\draw (10) -- (01);\draw (02) -- (13);\draw (12) -- (03);\draw (20) -- (31);\draw (30) -- (21);\draw (22) -- (33);\draw (32) -- (23);\draw (11) -- (22);\draw (21) -- (12);
\draw (30) -- (v); \draw (33) -- (v); \draw[bend left=90,looseness=0.8] (00) to (v); \draw[bend right=90,looseness=0.8] (03) to (v);
\end{tikzpicture}$\;\;\;\;\;\;\;\;$
\begin{tikzpicture}[largegraph,scale=0.7]
\node (00) at (1,0) {}; \node (01) at (2,0) {}; \node (02) at (1,1) {}; \node (03) at (2,1) {};
\node (10) at (3,0) {}; \node (11) at (4,0) {}; \node (12) at (3,1) {}; \node (13) at (4,1) {};
\node (20) at (4,2) {}; \node (21) at (5,2) {}; \node (22) at (4,3) {}; \node (23) at (5,3) {};
\node (30) at (2,4) {}; \node (31) at (3,4) {}; \node (32) at (2,5) {}; \node (33) at (3,5) {};
\node (40) at (0,2) {}; \node (41) at (1,2) {}; \node (42) at (0,3) {}; \node (43) at (1,3) {};

\draw (00) -- (01); \draw (00) -- (02); \draw (00) -- (03); \draw (01) -- (02); \draw (01) -- (03); \draw (02) -- (03);
\draw (10) -- (11); \draw (10) -- (12); \draw (10) -- (13); \draw (11) -- (12); \draw (11) -- (13); \draw (12) -- (13);
\draw (20) -- (21); \draw (20) -- (22); \draw (20) -- (23); \draw (21) -- (22); \draw (21) -- (23); \draw (22) -- (23);
\draw (30) -- (31); \draw (30) -- (32); \draw (30) -- (33); \draw (31) -- (32); \draw (31) -- (33); \draw (32) -- (33);
\draw (40) -- (41); \draw (40) -- (42); \draw (40) -- (43); \draw (41) -- (42); \draw (41) -- (43); \draw (42) -- (43);
\draw (00) -- (40); \draw (01) -- (10); \draw (11) -- (21); \draw (23) -- (33); \draw (32) -- (42);
\draw (02) -- (30); \draw (31) -- (13); \draw (12) -- (41); \draw (43) -- (22); \draw (20) -- (03);

\end{tikzpicture}
\caption{The Sierpi\'{n}ski graph $S(2,4)$, left, along with the extended Sierpi\'{n}ski graphs $S^+(2,4)$, centre, and $S^{++}(2,4)$, right.\label{fig-sierpinksi}}\end{center}\end{figure}
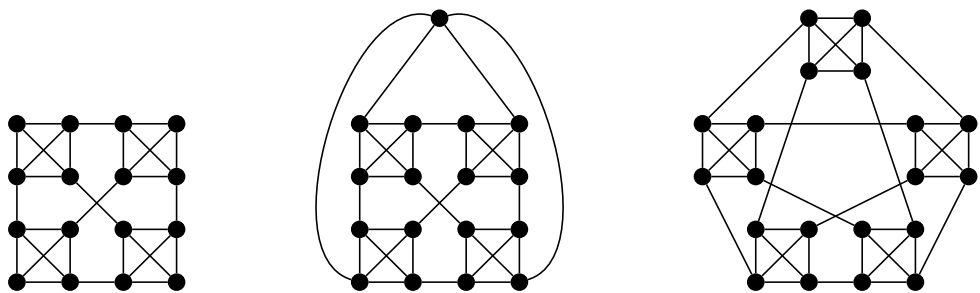

For $S(n,k)$, Klav\v{z}ar and Mohar only determined upper and lower bounds for the case $k = 4$, but for $S^+(n,k)$ and $S^{++}(n,k)$, they determined the crossing numbers in terms of the crossing numbers of complete graphs $K_{k+1}$. Hence, the latter are known precisely only for $k \leq 11$. For $k \leq 3$, all three constructions are planar.

\begin{theorem}[Klav\v{z}ar and Mohar, 2005 \cite{klavzarmohar2005}]For $k \geq 2$, the following hold:
\begin{align*}\frac{3}{16}4^n \leq cr(S(n,4)) & \leq \frac{1}{3}4^n - \frac{12n - 8}{3},\mbox{ for } n \geq 3,\\[0.1em]
cr(S^+(n,k)) & = \frac{k^n - 1}{k-1} cr(K_{k+1}),\mbox{ for } n \geq 1,\\[0.1em]
cr(S^{++}(n,k)) & = \frac{(k+1)k^{n-1} - 2}{k-1} cr(K_{k+1}),\mbox{ for } n \geq 1.\end{align*}\end{theorem}

\subsection{Star maps and pancake graphs}

In 1989, Akers and Krishnamurthy \cite{akerskrishnamurthy1989} proposed two recursive graph families as special cases of Cayley graphs, namely the Star maps (also known as Star graphs) and the pancake graphs. To avoid confusion with stars $S_n = K_{1,n}$, we will use the term Star map and a calligraphy $\mathcal{S}^*$.

The Star map $\mathcal{S}^*_n$ contains $n!$ vertices labelled with the $n!$ permutations on the set of symbols $1, 2, \hdots, n$. An edge $\{i,j\}$ exists if and only if the label for $j$ can be obtained from the label for $i$ by exchanging the first symbol with any other symbol. For example, in $\mathcal{S}^*_3$, the vertex labelled $123$ would be adjacent to $213$ and $321$.

In 2014, L\"{u} et al.\ \cite{luetal2014} considered $\mathcal{S}^*_4$, the smallest non-planar Star map, and determined its crossing number:

\begin{theorem}[L\"{u} et al.,\ 2014 \cite{luetal2014}] For the Star map $\mathcal{S}^*_4$,
\[ cr(\mathcal{S}^*_4) = 8. \]\end{theorem}

The pancake graph $\mathcal{P}_n$ contains $n!$ vertices labelled with the $n!$ permutations on the set of symbols $1, 2, \hdots, n$. An edge $\{i,j\}$ exists if and only if the label for $j$ is a permutation of the label for $i$ such that $i = i_1i_2i_3\cdots i_ki_{k+1}\cdots i_n$ and $j = i_k\cdots i_2i_1i_{k+1}\cdots i_n$ for some $k$ where $2 \leq k \leq n$.

In 2017, Yang et al.\ \cite{yangetal2017} considered $\mathcal{P}_4$, the smallest non-planar pancake graph, and determined its crossing number:

\begin{theorem}[Yang et al.,\ 2017 \cite{yangetal2017}] For the pancake graph $\mathcal{P}_4$,
\[ cr(\mathcal{P}_4) = 6. \]\end{theorem}

Results for larger Star maps or larger pancake graphs are currently unknown. $\mathcal{S}_4$ and $\mathcal{P}_4$ are displayed in \hyperref[fig-starpancake]{Figure \ref{fig-starpancake}}.

\begin{figure}[t!]\begin{center}\begin{tikzpicture}[largegraph,scale=0.6]
\node (00) at (0,1) {}; \node (01) at (1,0) {}; \node (02) at (2,1) {}; \node (03) at (2,2) {}; \node (04) at (1,3) {}; \node (05) at (0,2) {};
\node (10) at (4,1) {}; \node (11) at (5,0) {}; \node (12) at (6,1) {}; \node (13) at (6,2) {}; \node (14) at (5,3) {}; \node (15) at (4,2) {};
\node (20) at (4,5) {}; \node (21) at (5,4) {}; \node (22) at (6,5) {}; \node (23) at (6,6) {}; \node (24) at (5,7) {}; \node (25) at (4,6) {};
\node (30) at (0,5) {}; \node (31) at (1,4) {}; \node (32) at (2,5) {}; \node (33) at (2,6) {}; \node (34) at (1,7) {}; \node (35) at (0,6) {};

\foreach \n in {0,...,3}{
  \foreach \m in {0,...,4}{
    \pgfmathtruncatemacro{\nextm}{\m + 1}
    \draw (\n\m) -- (\n\nextm);
  }
  \draw (\n0) -- (\n5);
}
\draw (01) -- (11); \draw (24) -- (34); \draw (05) -- (30); \draw (13) -- (22);\draw (04) -- (14); \draw (21) -- (31);
\draw[bend right=90,looseness=0.4] (02) to (33); \draw[bend left=90,looseness=0.4] (10) to (25);
\draw (00) -- (20); \draw (03) -- (23); \draw (12) -- (32); \draw (15) -- (35);

\end{tikzpicture}$\;\;\;\;\;\;\;\;\;\;\;\;\;\;\;\;$
\begin{tikzpicture}[largegraph,scale=0.6]
\node (00) at (0,1) {}; \node (01) at (1,0) {}; \node (02) at (2,1) {}; \node (03) at (2,2) {}; \node (04) at (1,3) {}; \node (05) at (0,2) {};
\node (10) at (4,1) {}; \node (11) at (5,0) {}; \node (12) at (6,1) {}; \node (13) at (6,2) {}; \node (14) at (5,3) {}; \node (15) at (4,2) {};
\node (20) at (4,5) {}; \node (21) at (5,4) {}; \node (22) at (6,5) {}; \node (23) at (6,6) {}; \node (24) at (5,7) {}; \node (25) at (4,6) {};
\node (30) at (0,5) {}; \node (31) at (1,4) {}; \node (32) at (2,5) {}; \node (33) at (2,6) {}; \node (34) at (1,7) {}; \node (35) at (0,6) {};

\foreach \n in {0,...,3}{
  \foreach \m in {0,...,4}{
    \pgfmathtruncatemacro{\nextm}{\m + 1}
    \draw (\n\m) -- (\n\nextm);
  }
  \draw (\n0) -- (\n5);
}
\draw (01) -- (11); \draw (24) -- (34); \draw (05) -- (30); \draw (13) -- (22);\draw (04) -- (14); \draw (21) -- (31); \draw (03) -- (20); \draw (15) -- (32);
\draw[bend right=90,looseness=0.4] (02) to (33); \draw[bend left=90,looseness=0.4] (10) to (25);
\draw (00) to (23);
\draw (12) to (35);
\end{tikzpicture}
\caption{The Star map $\mathcal{S}_4$ and the pancake graph $\mathcal{P}_4$.\label{fig-starpancake}}\end{center}\end{figure}
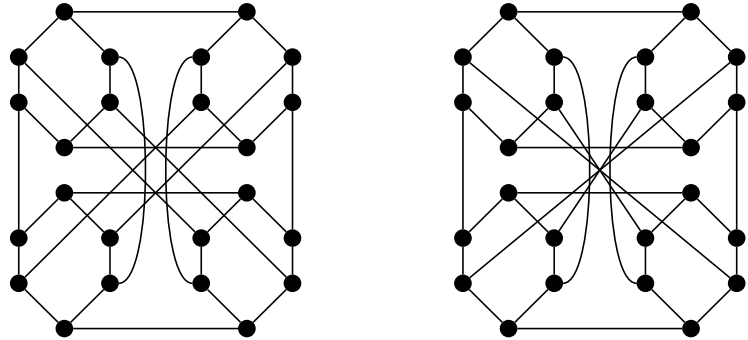

\newpage\subsection*{Acknowledgements}We are greatly indebted to the several anonymous referees whose suggestions significantly improved the flow and quality of this survey.

\null
\vfill
\clearpage
\appendix
\section{Appendices}

\subsection{List of graphs on six vertices}\label{app-g6}
There are 156 graphs up to isomorphism on six vertices, including 112 connected graphs and 44 disconnected graphs. We list them here in the order proposed in Frank Harary's classic textbook, Graph Theory \cite{harary1969}. The graphs are ordered by their number of edges.

\begin{figure}[hb!]
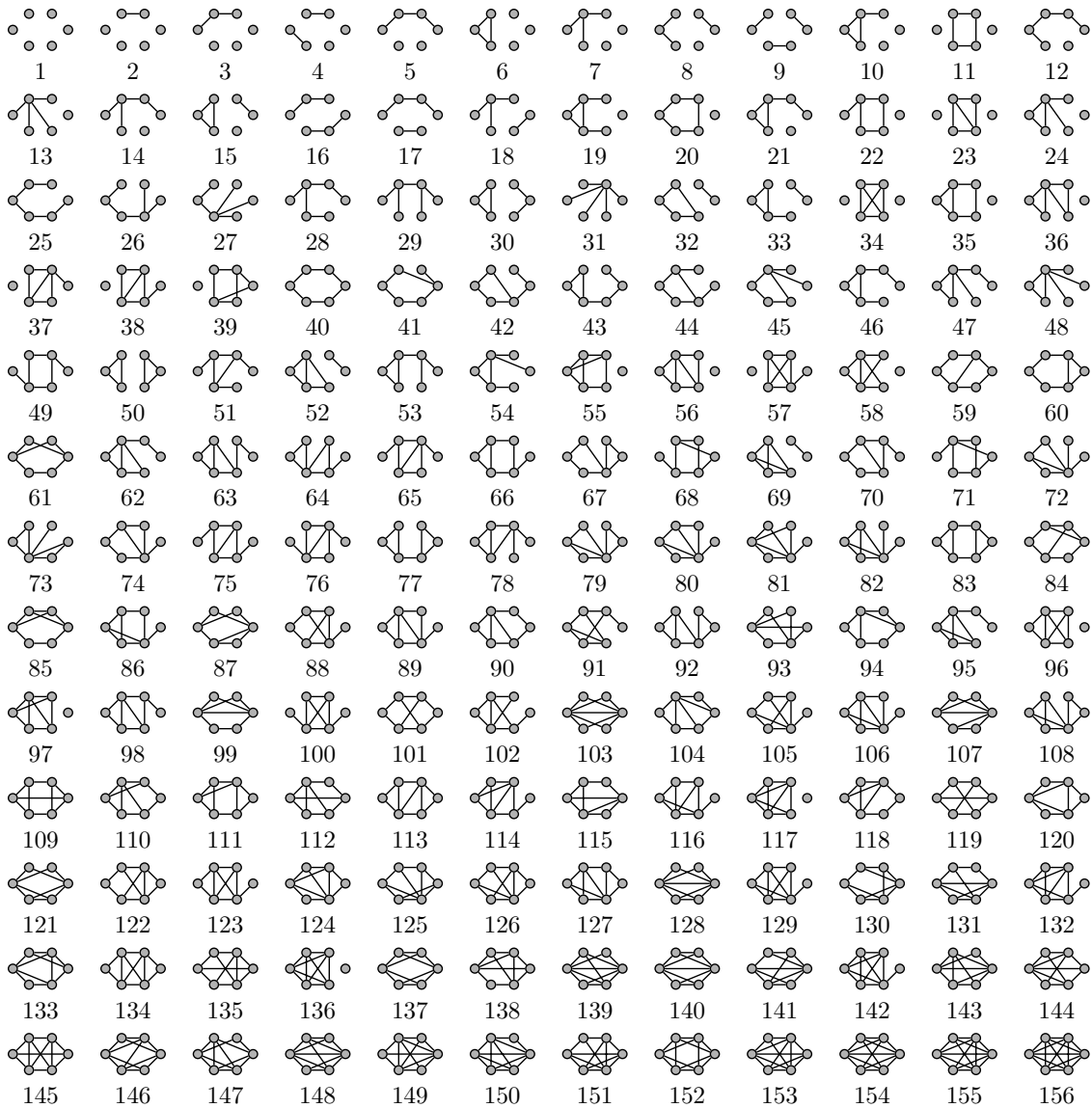
\resizebox{\textwidth}{!}{
\begin{tabular}{cccccccccccc}
\ctikz{\ngraph{6.1}} & \ctikz{\ngraph{6.2}} & \ctikz{\ngraph{6.3}} & \ctikz{\ngraph{6.4}} & \ctikz{\ngraph{6.5}} & \ctikz{\ngraph{6.6}} & \ctikz{\ngraph{6.7}} & \ctikz{\ngraph{6.8}} & \ctikz{\ngraph{6.9}} & \ctikz{\ngraph{6.10}} & \ctikz{\ngraph{6.11}} & \ctikz{\ngraph{6.12}} \\
{\small 1} & {\small 2} & {\small 3} & {\small 4} & {\small 5} & {\small 6} & {\small 7} & {\small 8} & {\small 9} & {\small 10} & {\small 11} & {\small 12}\\
\ctikz{\ngraph{6.13}} & \ctikz{\ngraph{6.14}} & \ctikz{\ngraph{6.15}} & \ctikz{\ngraph{6.16}} & \ctikz{\ngraph{6.17}} & \ctikz{\ngraph{6.18}} & \ctikz{\ngraph{6.19}} & \ctikz{\ngraph{6.20}} & \ctikz{\ngraph{6.21}} & \ctikz{\ngraph{6.22}} & \ctikz{\ngraph{6.23}} & \ctikz{\ngraph{6.24}} \\
{\small 13} & {\small 14} & {\small 15} & {\small 16} & {\small 17} & {\small 18} & {\small 19} & {\small 20} & {\small 21} & {\small 22} & {\small 23} & {\small 24} \\
 \ctikz{\ngraph{6.25}} & \ctikz{\ngraph{6.26}} & \ctikz{\ngraph{6.27}} & \ctikz{\ngraph{6.28}} & \ctikz{\ngraph{6.29}} & \ctikz{\ngraph{6.30}} & \ctikz{\ngraph{6.31}} & \ctikz{\ngraph{6.32}} & \ctikz{\ngraph{6.33}} & \ctikz{\ngraph{6.34}} & \ctikz{\ngraph{6.35}} & \ctikz{\ngraph{6.36}} \\
{\small 25} & {\small 26} & {\small 27} & {\small 28} & {\small 29} & {\small 30} & {\small 31} & {\small 32} & {\small 33} & {\small 34} & {\small 35} & {\small 36}\\
\ctikz{\ngraph{6.37}} & \ctikz{\ngraph{6.38}} & \ctikz{\ngraph{6.39}} & \ctikz{\ngraph{6.40}} & \ctikz{\ngraph{6.41}} & \ctikz{\ngraph{6.42}} & \ctikz{\ngraph{6.43}} & \ctikz{\ngraph{6.44}} & \ctikz{\ngraph{6.45}} & \ctikz{\ngraph{6.46}} & \ctikz{\ngraph{6.47}} & \ctikz{\ngraph{6.48}} \\
{\small 37} & {\small 38} & {\small 39} & {\small 40} & {\small 41} & {\small 42} & {\small 43} & {\small 44} & {\small 45} & {\small 46} & {\small 47} & {\small 48} \\
\end{tabular}}\caption{All 156 graphs on six vertices and their indices, ordered by number of edges.}\end{figure}

\newpage
\subsection{Tables of citations}
\subsubsection{Generalized Petersen graphs}
\label{app-gpref}
\setlength\tabcolsep{1 pt}
\begin{table}[!h]\centering
\caption{Publications where the crossing numbers for $GP(n,k)$ were first proved.}\smallskip\smallskip\label{tab-GP-ref}
\begin{tabular}{|c|c|c|c|c|c|c|c|c|c|c|c|c|c|}\hline
\multicolumn{1}{|r|}{\;\;\; \bf $n$ \;} & \multirow{2}{*}{\: \bf 5 \:} & \multirow{2}{*}{\: \bf 6 \:} & \multirow{2}{*}{\: \bf 7 \:} & \multirow{2}{*}{\: \bf 8 \:} & \multirow{2}{*}{\: \bf 9 \:} & \multirow{2}{*}{\: \bf 10 \:} & \multirow{2}{*}{\: \bf 11 \:} & \multirow{2}{*}{\: \bf 12 \:} & \multirow{2}{*}{\: \bf 13 \:} & \multirow{2}{*}{\: \bf 14 \:} & \multirow{2}{*}{\: \bf 15 \:} & \multirow{2}{*}{\: \bf 16 \:} & \multirow{2}{*}{\: \bf 17 \:} \\
\multicolumn{1}{|c|}{\bf $k$\;\;} & & & & & & & & & & & & & \\
\hline \rule{0pt}{2.5ex} {\bf 1}   &  \cite{exooetal1981}  &  \cite{exooetal1981}  &  \cite{exooetal1981}  &  \cite{exooetal1981}  &  \cite{exooetal1981}  & \cite{exooetal1981}  &  \cite{exooetal1981}  &  \cite{exooetal1981}  &  \cite{exooetal1981}  &  \cite{exooetal1981}  &  \cite{exooetal1981}  &  \cite{exooetal1981}  &  \cite{exooetal1981}\\
\hline \rule{0pt}{2.5ex} {\bf 2}   &  \cite{exooetal1981}  &  \cite{exooetal1981}  &  \cite{exooetal1981}  &  \cite{exooetal1981}  &  \cite{exooetal1981}  &  \cite{exooetal1981}  &  \cite{exooetal1981}  &  \cite{exooetal1981}  &  \cite{exooetal1981}  &  \cite{exooetal1981}  &  \cite{exooetal1981}  &  \cite{exooetal1981}  &  \cite{exooetal1981}\\
\hline \rule{0pt}{2.5ex} {\bf 3}   &    &  \cite{guyharary1967} &  \cite{richtersalazar2002}  &  \cite{fiorini1986}  &  \cite{fiorini1986}  &  \cite{richtersalazar2002}  &  \cite{fiorini1986}  &  \cite{fiorini1986}  &  \cite{richtersalazar2002}  &  \cite{fiorini1986}  &  \cite{fiorini1986}  &  \cite{richtersalazar2002}  &  \cite{fiorini1986}\\
\hline \rule{0pt}{2.5ex} {\bf 4}   &    &    &    &  \cite{guyharary1967}  &  \cite{maetal2005}  &  \cite{sarazin1997}  &  \cite{linetal2009}  &  \cite{fiorini1986}  &  \cite{linetal2009}  &  \cite{linetal2009}  &  \cite{linetal2009}  &  \cite{fiorini1986}  &  \\
\hline \end{tabular}
\end{table}
\setlength\tabcolsep{6 pt}
\newpage
\subsubsection{Cartesian products with 5-vertex graphs}
\label{app-c5}
\begin{table}[!h]\begin{center}\vspace*{-0.7cm}\caption{Publications where the crossing numbers of Cartesian products of 5-vertex graphs with paths, cycles and stars were first proved. Empty cells imply the crossing number is not known and a dash indicates that the graphs are planar.\label{tab-5vcart-ref}}\smallskip\smallskip$\begin{array}{|c|c|c|c|c|}\hline
\rule{0pt}{2.3ex} i         & G^5_i                  & cr(G^5_i \Box P_n) & cr(G^5_i \Box C_n) & cr(G^5_i \Box S_n) \\
\hline 1  & \ctikz{\ngraph{5.1}}   &   -    &  -   & \cite{bokal2007}                    \\
\hline 2  & \ctikz{\ngraph{5.2}}   & \cite{klesc1991} & \cite{klesc1991}   & \cite{huangzhao2008}, \cite{ho2008}\\
\hline 3  & \ctikz{\ngraph{5.3}}   & \cite{klesc2001_2}    & \cite{klesc2005}   & \cite{bokal2007_2}          \\
\hline 4  & \ctikz{\ngraph{5.4}}   & \cite{klesc2001_2}  & \cite{klesc2001_2}  & \cite{klesc2009}   \\
\hline 5  & \ctikz{\ngraph{5.5}}   & \cite{klesc2001_2} & \cite{klesc2001_2}&                                                   \\
\hline 6  & \ctikz{\ngraph{5.6}}   & \cite{klesc2001_2} & \cite{klesc2005}  & \cite{klesc2009}                              \\
\hline 7  & \ctikz{\ngraph{5.7}}   & \cite{klesc2001_2} & \cite{klesc2005} & \unv{\cite{zhanghuang2011}}    \\
\hline\end{array}$\end{center}\end{table}

\newpage
\subsubsection{Cartesian products of paths with 6-vertex graphs}
\label{app-c6p}
\begin{table}[!h]\begin{center}\vspace*{-0.25cm}\caption{Publications where the crossing numbers of Cartesian products of 6-vertex graphs with paths were first proved. \label{tab-6vcart_path-ref}}\smallskip\smallskip$\begin{array}{|c|c|c||c|c|c||c|c|c|}\hline \rule{0pt}{2.3ex} i & G^6_i 	& cr(G^6_i \Box P_n) 	& i & G^6_i 	& cr(G^6_i \Box P_n) 	& i & G^6_i 	& cr(G^6_i \Box P_n)	\\
\hline 25  &\ctikz{\ngraph{6.25}}	& \text{-}        	                &61&\ctikz{\ngraph{6.61}}  	& \cite{klescpetrillova2013_2}  	&89&\ctikz{\ngraph{6.89}}  	    & \cite{klescetal2014}              \\	
\hline 26  &\ctikz{\ngraph{6.26}}	& \cite{klescpetrillova2013_2}   	&62&\ctikz{\ngraph{6.62}}  	& \cite{petrillova2019}          	&90&\ctikz{\ngraph{6.90}}  	    & \cite{clancyetaltoappear}         \\
\hline 27  &\ctikz{\ngraph{6.27}}	& \cite{klescpetrillova2013_2}   	&63&\ctikz{\ngraph{6.63}}  	& \cite{petrillova2019}             &91&\ctikz{\ngraph{6.91}}  	    & \cite{petrillova2015}             \\
\hline 28  &\ctikz{\ngraph{6.28}}	& \cite{klescpetrillova2013_2}   	&64&\ctikz{\ngraph{6.64}}  	& \cite{klescpetrillova2013_2}  	&93&\ctikz{\ngraph{6.93}} 	    &\unv{ \cite{wanghuang2008_2}}	    \\
\hline 29  &\ctikz{\ngraph{6.29}}	& \cite{klescpetrillova2013_2}   	&65&\ctikz{\ngraph{6.65}}  	& \cite{klescetal2014}          	&94&\ctikz{\ngraph{6.94}}  	    & \cite{klescpetrillova2013_2}	    \\
\hline 31  &\ctikz{\ngraph{6.31}}	& \cite{bokal2007}               	&66&\ctikz{\ngraph{6.66}}  	& \cite{klescpetrillova2013_2}  	&103&\ctikz{\ngraph{6.103}} 	& \cite{ouyangetal2014}	            \\
\hline 40  &\ctikz{\ngraph{6.40}}	& \text{-}        	                &68&\ctikz{\ngraph{6.68}}  	& \cite{petrillova2015}         	&104&\ctikz{\ngraph{6.104}} 	& \cite{klescpetrillova2013_2}	    \\
\hline
\end{array}$\end{center}\end{table}
\newpage
\subsubsection{Cartesian products of cycles with 6-vertex graphs}
\vspace{-0.5cm}
\label{app-c6c}
\begin{table}[!h]\begin{center}\caption{Publications where the crossing numbers of Cartesian products of 6-vertex graphs with cycles were first proved.\label{tab-6vcart_cycle-ref}}\smallskip\smallskip\scalebox{0.82}{$\begin{array}{|c|c|c||c|c|c|}\hline \rule{0pt}{2.3ex} i & G^6_i & cr(G^6_i \Box C_n) & i & G^6_i & cr(G^6_i \Box C_n) \\
\hline 25     & \ctikz{\ngraph{6.25}}  & \text{-}                  &     60     & \ctikz{\ngraph{6.60}}  & \cite{clancyetaltoappear} \\
\hline 26     & \ctikz{\ngraph{6.26}}  & \cite{drazenskaklesc2007} &     63     & \ctikz{\ngraph{6.63}}  & \cite{clancyetaltoappear} \\
\hline 27     & \ctikz{\ngraph{6.27}}  & \cite{drazenskaklesc2007} &     64     & \ctikz{\ngraph{6.64}}  & \cite{clancyetaltoappear} \\
\hline 28     & \ctikz{\ngraph{6.28}}  & \cite{drazenskaklesc2007} &     66     & \ctikz{\ngraph{6.66}}  & \cite{clancyetaltoappear} \\
\hline 29     & \ctikz{\ngraph{6.29}}  & \cite{drazenskaklesc2007} &     67     & \ctikz{\ngraph{6.67}}  & \cite{drazenskaklesc2011} \\
\hline 40     & \ctikz{\ngraph{6.40}}  & \cite{richtersalazar2001} &     70     & \ctikz{\ngraph{6.70}}  & \cite{clancyetaltoappear} \\
\hline 41     & \ctikz{\ngraph{6.41}}  & \cite{drazenskaklesc2011} &     75     & \ctikz{\ngraph{6.75}}  & \cite{clancyetaltoappear} \\
\hline\end{array}$}\end{center}\end{table}
\vspace{-0.5cm}
\subsubsection{Cartesian products of stars with 6-vertex graphs}
\label{app-c6s}
\begin{table}[!h]\begin{center}\vspace*{-0.5cm}\caption{Publications where the crossing numbers of Crossing numbers of Cartesian products of 6-vertex graphs with stars were first proved. \label{tab-6vcart_star-ref}}\smallskip\smallskip\scalebox{0.82}{$\begin{array}{|c|c|c||c|c|c|}\hline \rule{0pt}{2.3ex} i & G^6_i & cr(G^6_i \Box S_n)  & i & G^6_i & cr(G^6_i \Box S_n) \\
\hline 25     & \ctikz{\ngraph{6.25}}  & \cite{bokal2007}          & 77     & \ctikz{\ngraph{6.77}}  & \cite{klescschrotter2013}    \\
\hline 26     & \ctikz{\ngraph{6.26}}  & \cite{bokal2007_2}        & 79     & \ctikz{\ngraph{6.79}}  & \cite{klescschrotter2013}    \\
\hline 27     & \ctikz{\ngraph{6.27}}  & \cite{klescschrotter2013} & 80     & \ctikz{\ngraph{6.80}}  & \cite{klescschrotter2013}    \\
\hline 28     & \ctikz{\ngraph{6.28}}  & \cite{bokal2007_2}        & 85     & \ctikz{\ngraph{6.85}}  & \unv{\cite{zhangetal2008_2}} \\
\hline 29     & \ctikz{\ngraph{6.29}}  & \cite{bokal2007_2}        & 93     & \ctikz{\ngraph{6.93}}  & \unv{\cite{luhuang2010}}     \\
\hline 31     & \ctikz{\ngraph{6.31}}  & \cite{klescschrotter2013} & 94     & \ctikz{\ngraph{6.94}}  & \cite{wangetal2013_2}        \\
\hline 43     & \ctikz{\ngraph{6.43}}  & \cite{klescschrotter2013} & 104    & \ctikz{\ngraph{6.104}} & \cite{klescschrotter2013}    \\
\hline
\end{array}$}\end{center}\end{table}
\newpage
\subsubsection{Join products with 5-vertex graphs}
\vspace{-0.5cm}
\label{app-j5}
\begin{table}[!h]\begin{center}\caption{Publications where the crossing numbers of joins of connected 5-vertex graphs with discrete graphs, paths and cycles were first proved. Empty cells imply that the crossing number has not yet been determined.\label{tab-5vjoin-ref}}\smallskip\smallskip
$\begin{array}{|c|c|c|c|c|}\hline
\rule{0pt}{2.3ex} i         & G^5_i                      & cr(G^5_i + D_n)    & cr(G^5_i + P_{n-1})          & cr(G^5_i + C_n)   \\
\hline 1  & \ctikz{\ngraph{5.1}}       & \cite{kleitman1971}    & \cite{klesc2007}  & \cite{klesc2007}\\
\hline 2  & \ctikz{\ngraph{5.2}}       & \cite{huangzhao2008}, \cite{ho2008}  & \unv{\cite{lietal2008}}& \cite{klescetal2017} \\
\hline 3  & \ctikz{\ngraph{5.3}}       & \cite{stas2019_2} & \cite{stas2019_2} & \cite{stas2019_2}\\
\hline 4  & \ctikz{\ngraph{5.4}}       & \cite{klescstastoappear} & \cite{klescstastoappear} & \cite{klescstastoappear}\\
\hline 5  & \ctikz{\ngraph{5.5}}       & \cite{bereznystas2017} & \cite{staspetrillova2018}                   & \cite{staspetrillova2018}\\
\hline 6  & \ctikz{\ngraph{5.6}}       & \cite{stas2019_4}      & \cite{stas2019_4}& \unv{\cite{yueetal2013}}\\
\hline 7  & \ctikz{\ngraph{5.7}}       & \cite{dinghuang2018} & \unv{\cite{zhenghuang2011}} & \cite{klescschrotter2012}\\
\hline\end{array}$\end{center}\end{table}

\newpage
\subsubsection{Join products with 6-vertex graphs}
\vspace{-0.5cm}
\label{app-j6}
\begin{table}[!h]\centering\vspace*{-0.3cm}\caption{Publications where the crossing numbers of joins of particular 6-vertex graphs with discrete graphs, paths and cycles were first proved.\label{tab-6vjoin-ref}}\smallskip\smallskip
\scalebox{0.82}{$\begin{array}{|c|c|c|c|c|}\hline
\rule{0pt}{2.3ex} i & G^6_i & cr(G^6_i + D_n) & cr(G^6_i + P_{n-1}) & cr(G^6_i + C_n)\\
\hline 20 &  \ctikz{\ngraph{6.20}}      & \cite{bereznystas2019}              &                &                   \\
\hline 25 &  \ctikz{\ngraph{6.25}}      & \cite{kleitman1971}              & \cite{klesc2007} & \cite{klesc2007} \\
\hline 31 &  \ctikz{\ngraph{6.31}}      & \unv{\cite{meihuang2007}}      & \unv{\cite{lietal2008_2}} & \cite{klescetal2017} \\
\hline 35 &  \ctikz{\ngraph{6.35}}      & \cite{stas2020_3}             &                 &                  \\
\hline 40 &  \ctikz{\ngraph{6.40}}      & \cite{kleitman1971}         & \cite{klesc2007} & \cite{klesc2007} \\
\hline 41 &  \ctikz{\ngraph{6.41}}      & \cite{bereznystas2019}         & \cite{bereznystas2019} & \\
\hline 44 &  \ctikz{\ngraph{6.44}}      & \cite{stas2018} & \unv{\cite{lietal2008_2}} &                          \\
\hline\end{array}$}\end{table}
\null
\vfill
\clearpage
\subsection{Glossary of symbols}\label{app-glossary}

\hyperref[tab-glossary]{Table \ref{tab-glossary}} contains an alphabetized list of symbols which are used in this survey. Care has been taken to not use any symbols for multiple purposes.

\begin{table}[h!]\centering\caption{List of symbols used in this survey, with descriptions.}\smallskip\smallskip\label{tab-glossary}
\renewcommand*{\arraystretch}{1.5}\resizebox{15cm}{!}{\begin{tabular}{|l|l||l|l|}\hline {\bf Symbol} & {\bf Description} & {\bf Symbol} & {\bf Description} \\
\hline $AQ_n$                  & 	 Augmented cube                 &	 $\mathcal{P}_n$          &	 Pancake graph                        \\
\hline $B(r)$                  & 	 Benes network                  &	 $Q_n$                    &	 $n$-dimensional hypercube            \\
\hline $BF(r)$                 & 	 Butterfly graph                &	 $S_n$                    &	 Star on $n+1$ vertices               \\
\hline $C_n$                   & 	 Cycle on $n$ vertices          &	 $\mathcal{S}_n$          &	 Sunlet graph                         \\
\hline $\overline{C}_n$        & 	 Complement of $C_n$            &	 $\mathcal{S}^*_n$        &	 Star map                             \\
\hline $CCC_n$                 & 	 Cube connected cycle graph     &	 $S(n,k)$                 &	 Sierpi\'{n}ski graph                 \\
\hline $Ci_n(L)$               & 	 Circulant graph                &	 $S^+(n,k)$               &	 Extended Sierpi\'{n}ski graph        \\
\hline
\end{tabular}}\end{table}

\newpage
\subsection{Incorrect results}\label{app-error}

Throughout the history of research into crossing numbers, it is unfortunately common for papers to contain errors. Indeed, the inaugural result relating to crossing numbers, Zarankiewicz's proof for the crossing number of complete bipartite graphs, was later found to have an error and the result is still not resolved to this day. In this section we attempt to detail any publications which contain incorrect proofs or results. In some cases, the flaws were found and either simply noted, or in some cases corrected. In other cases, the flaws have not been widely recognised, in which case we provide a counterexample demonstrating the incorrect result. The intention is not to disparage the authors, but rather to ensure that researchers do not use the flawed results as a basis for future proofs, and to provide references containing the corrected results when such results exist.

The incorrect results are listed in order of appearance. When the result is known to be incorrect due to a later publication, references to the publications are given. When the result is known to be incorrect due to the existence of a better drawing, we show such a counterexample as obtained by QuickCross \cite{clancyetal2019}. When the result is known to be incorrect because a drawing does not exist with the proposed number of crossings, we use Crossing Number Web Compute \cite{chimaniwiederasite2016} to determine the true crossing number for a minimal counterexample.

\begin{itemize}\item In 1955, Zarankiewicz \cite{zarankiewicz1955} claimed to have determined that the crossing number of complete bipartite graphs $K_{m,n}$. The flaw in the proof was noted and communicated privately, and subsequently described by Guy \cite{guy1969} in 1969. The result is now known to be correct for $\min\{m,n\} \leq 6$, and the special cases $m \leq 8$, $n \leq 10$; see \hyperref[sec-complete-bipartite]{Section \ref{sec-complete-bipartite}}.

\item In 1986, Fiorini \cite{fiorini1986} claimed to have determined that the crossing number of the generalized Petersen graph $GP(10,3)$ was equal to four. This was shown to be false in 1992 by McQuillan and Richter \cite{mcquillanrichter1992} who conjectured that the true value would be six, which was finally proved by Richter and Salazar \cite{richtersalazar2002} in 2002. Fiorini (along with co-author Gauci) also corrected his earlier proofs in a paper published in 2003 \cite{fiorinigauci2003}; see \hyperref[sec-GP]{Section \ref{sec-GP}}.

\item In 1996, Cimikowski \cite{cimikowski1996} claimed to have determined an upper bound for the crossing number of the complete mesh of trees $M_n$, but in 1998 he reported the error himself and provided a corrected upper bound \cite{cimikowski1998}; see \hyperref[sec-cmot]{Section \ref{sec-cmot}}.

\item In 2005, Wang and Huang \cite{wanghuang2005} determined the crossing number of the Cartesian product of four six-vertex graphs with paths. For the first of these graphs, $G^6_{131}$ (which they label $G_1$), they claimed that $cr(G^6_{131} \Box P_n) = 3n - 1$ for $n \geq 1$. This is not correct. Indeed, we used Crossing Number Web Compute to show that $cr(G^6_{131} \Box P_1) = 4$, rather than 2 as would be suggested by \cite{wanghuang2005}.

\item In 2007, Dra\v{z}ensk\'{a} and Kle\v{s}\v{c} \cite{drazenskaklesc2007} determined the crossing numbers of $G^6_i \Box C_n$ for $i = 26, 27, 28$, and claimed that the results held for $n \geq 3$. However, in their proof they used the general formula for $cr(K_{1,4} \Box C_n)$ which only holds for $n \geq 6$. Hence, the values they gave for $n = 3, 4, 5$ were incorrect. We have used Crossing Number Web Compute to show that $cr(G^6_{26} \Box C_3) = 1$, $cr(G^6_{26} \Box C_4) = 2$, $cr(G^6_{26} \Box C_5) = 4$, $cr(G^6_{27} \Box C_3) = 2$, $cr(G^6_{27} \Box C_4) = 4$, $cr(G^6_{27} \Box C_5) = 8$, $cr(G^6_{28} \Box C_3) = 1$, $cr(G^6_{28} \Box C_4) = 2$, and $cr(G^6_{28} \Box C_5) = 4$.

\item In 2008, He \cite{he2008} considered the crossing number of $K_{4,n}$ with two edges deleted. In particular, if the 4 vertices are denoted $y_1$, $y_2$, $y_3$, $y_4$ and the $n$ vertices are denoted $x_1$, $x_2, \hdots, x_n$, then they remove edges $\{x_1,y_1\}$ and $\{x_1, y_2\}$. He claimed to have proved the crossing number is equal to $Z(4,n) - 2\fl[n] + 2$ for $n \geq 4$. The result appears to be incorrect for odd values of $n$. In particular, in \hyperref[fig-err1]{Figure \ref{fig-err1}} we display a drawing of $K_{4,5}$ with two edges removed as described above, with only four crossings, rather than the six suggested by \cite{he2008}.

\begin{figure}[h!]\begin{center}\begin{tikzpicture}[largegraph,scale=0.8]
\foreach \i in {0,...,3} \node (\i) at (90-90*\i:3) {};
\foreach \i in {4,...,7} \node (\i) at (90-90*\i:1) {};
\node (8) at (270:2) {};
\draw (0) -- (1) -- (2) -- (3) -- (0);
\draw (4) -- (5) -- (6) -- (7) -- (4);
\draw (6) -- (8) -- (2);
\foreach \i in {0,...,3}{
  \pgfmathtruncatemacro{\next}{4+mod(\i+1,4)}
  \pgfmathtruncatemacro{\prev}{4+mod(\i+3,4)}
  \draw (\prev) -- (\i) -- (\next);
}
\end{tikzpicture}\caption{A drawing of $K_{4,5}$ minus two edges, with only four crossings.}\label{fig-err1}\end{center}\end{figure}
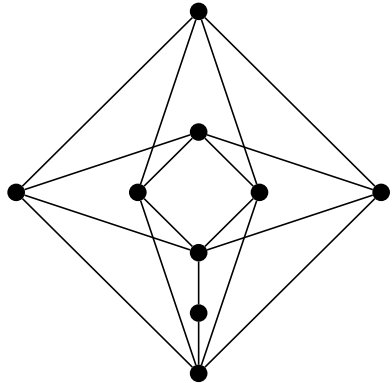

\item In 2011, He \cite{he2011} considered the crossing number of $G^5_3 \Box S_n$ and $G^5_4 \Box S_n$ and claimed that both were equal to $4\fl[n]\fl[n-1] + \fl[n]$ for $n \geq 1$. These results are contradicted by Bokal (2007) \cite{bokal2007_2} and Kle\v{s}\v{c} (2009) \cite{klesc2009} respectively, who showed that each has crossing number equal to $3\fl[n]\fl[n-1] + \fl[n]$. Tests with QuickCross have confirmed that it is possible to find drawings with the latter number of crossings, e.g.\ see \hyperref[fig-err2]{Figure \ref{fig-err2}} in which drawings of $G^5_3 \Box S_3$ and $G^5_4 \Box S_3$ are drawn with four crossings, rather than the five suggested by \cite{he2011}.

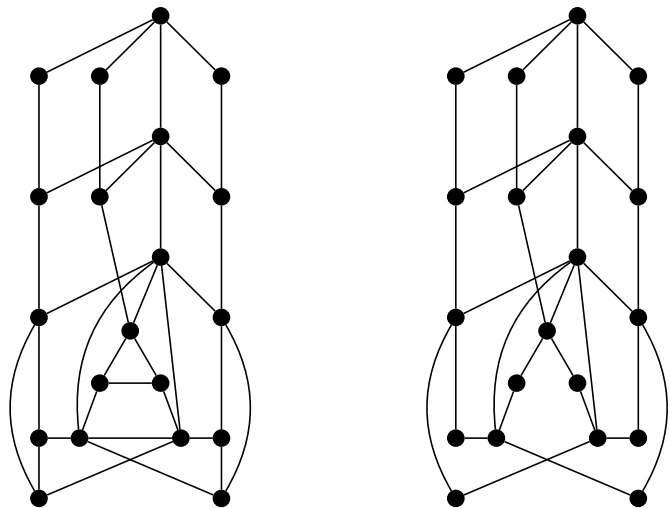
\begin{figure}[h!]\begin{center}\begin{tikzpicture}[largegraph,scale=0.8]
\begin{scope}[shift={(0,4)}]\node (0) at (2,1) {};\node (1) at (0,0) {};\node (2) at (1,0) {};\node (3) at (3,0) {};\end{scope}
\begin{scope}[shift={(0,2)}]\node (4) at (2,1) {};\node (5) at (0,0) {};\node (6) at (1,0) {};\node (7) at (3,0) {};\end{scope}
\begin{scope}[shift={(0,0)}]\node (8) at (2,1) {};\node (9) at (0,0) {};\node (11) at (3,0) {};\end{scope}
\node (12) at (0,-2) {};\node (18) at (2/3,-2) {};\node (19) at (7/3,-2) {};\node (13) at (3,-2) {};\node (14) at (0,-3) {};\node (15) at (3,-3) {};
\pgfmathsetmacro{\rad}{1/sqrt(3)}
\begin{scope}[shift={(1.5,-0.8)}]\node (10) at (90:\rad) {};\node (16) at (210:\rad) {};\node (17) at (330:\rad) {};\end{scope}
\foreach \i in {1,2,3} \draw (0) -- (\i);
\foreach \i in {5,6,7} \draw (4) -- (\i);
\foreach \i in {9,10,11} \draw (8) -- (\i);
\draw (1) -- (5) -- (9) -- (12) -- (14);\draw (2) -- (6) -- (10);\draw (0) -- (4) -- (8);\draw (3) -- (7) -- (11) -- (13) -- (15);
\draw (12) -- (18) -- (19) -- (13);\draw (14) -- (19);\draw (15) -- (18);\draw (8) -- (19);\draw[bend right] (8) to (18);
\draw (10) -- (16) -- (17) -- (10);\draw (16) -- (18);\draw (17) -- (19);\draw[bend right] (9) to (14);\draw[bend left] (11) to (15);
\end{tikzpicture} \;\;\;\;\;\;\;\;\;\;\;\;\;\;\;\;\;\begin{tikzpicture}[largegraph,scale=0.8]
\begin{scope}[shift={(0,4)}]\node (0) at (2,1) {};\node (1) at (0,0) {};\node (2) at (1,0) {};\node (3) at (3,0) {};\end{scope}
\begin{scope}[shift={(0,2)}]\node (4) at (2,1) {};\node (5) at (0,0) {};\node (6) at (1,0) {};\node (7) at (3,0) {};\end{scope}
\begin{scope}[shift={(0,0)}]\node (8) at (2,1) {};\node (9) at (0,0) {};\node (11) at (3,0) {};\end{scope}
\node (12) at (0,-2) {};\node (18) at (2/3,-2) {};\node (19) at (7/3,-2) {};\node (13) at (3,-2) {};\node (14) at (0,-3) {};\node (15) at (3,-3) {};
\pgfmathsetmacro{\rad}{1/sqrt(3)}
\begin{scope}[shift={(1.5,-0.8)}]\node (10) at (90:\rad) {};\node (16) at (210:\rad) {};\node (17) at (330:\rad) {};\end{scope}
\foreach \i in {1,2,3} \draw (0) -- (\i);
\foreach \i in {5,6,7} \draw (4) -- (\i);
\foreach \i in {9,10,11} \draw (8) -- (\i);
\draw (1) -- (5) -- (9) -- (12);\draw (2) -- (6) -- (10);\draw (0) -- (4) -- (8);\draw (3) -- (7) -- (11) -- (13);\draw (12) -- (18);
\draw (19) -- (13);\draw (14) -- (19);\draw (15) -- (18);\draw (8) -- (19);\draw[bend right] (8) to (18);\draw (16) -- (10) -- (17);
\draw (16) -- (18);\draw (17) -- (19);\draw[bend right] (9) to (14);\draw[bend left] (11) to (15);
\end{tikzpicture}\caption{Drawings of $G^5_3 \Box S_3$ and $G^5_4 \Box S_3$, with only four crossings.}\label{fig-err2}\end{center}\end{figure}

\item In 2015, Cruz and Japson \cite{cruzjapson2015} considered the crossing number of $GP(17,5)$ and claimed it was equal to 14. This result is contradicted by Newcombe (2019) \cite{newcombe2019}, who presented a drawing of $GP(17,5)$ with only 13 crossings.

\item In 2015, Li \cite{li2015} considered the join product of $G^5_{16}$ with discrete graphs, path graphs and cycles. In particular, they claim that $cr(G^5_{16} + P_{n-1}) = 4\fl[n]\fl[n-1] + \fl[n] + n + 1$ for $n \geq 2$. However, in Figure \ref{fig-err3} we display a drawing of $G^5_{16} + P_1$ with only three crossings, rather than four as suggested by \cite{li2015}.

\begin{figure}[h!]\begin{center}\begin{tikzpicture}[largegraph,scale=0.5]
\node (1) at (1,4) {};
\node (2) at (4,7) {};
\node (3) at (4,5) {};
\node (4) at (4,3) {};
\node (5) at (3,1) {};
\node (6) at (5,1) {};
\node (7) at (7,4) {};

\draw (1) -- (2) -- (7) -- (3) -- (1) -- (4) -- (7) --(6) -- (5) -- (1);
\draw (2) -- (3) -- (4) -- (5);
\draw (5) -- (6) -- (3) -- (5);
\draw (4) -- (6);
\draw [bend right=90] (1) to (6);
\draw [bend right=90] (5) to (7);
\end{tikzpicture}\caption{A drawing of $G^5_{16} + P_1$, with only three crossings.}\label{fig-err3}\end{center}\end{figure}

\item In 2016, Hsieh and Lin \cite{hsiehlin2016} claimed to have determined the crossing number of the join product of various path powers with discrete graphs and path graphs. All of the claimed results appear to be incorrect, as they all rely on an intermediate result, Lemma 5 of \cite{hsiehlin2016}, which claims that $cr(P^{m-1}_m + D_n) = cr(P^{m-2}_m + D_n)$ for $3 \leq m \leq 6$ and $n \geq 1$. This result is incorrect. In \cite{hsiehlin2016}, $P^k_m$ is defined as the $k$-th power on the path graph with $m$ vertices, rather than $m$ edges. The minimal counterexample can be seen by considering graphs on four vertices. By Hsieh and Lin's definition, $P^3_4 = K_4$ and $P^2_4 = K_4 \setminus e$. However, we know from Ho \cite{ho2009} that $cr(K_4 + D_n) = 2\fl[n]\fl[n-1] + n$, and from Kle\v{s}\v{c} and Schr\"{o}tter \cite{klescschrotter2011} that $cr((K_4 \setminus e) + D_n) = 2\fl[n]\fl[n-1] + \fl[n].$

\item In 2016, Vijaya et al.\ \cite{vijayaetal2016} considered the join product of $G^6_{133}$ with discrete graphs, paths and cycles. They provided detailed proofs for the cases of discrete graphs and paths, but omitted the proof for cycles, for which they claimed that $cr(G^6_{133} + C_n) = Z(6,n) + n + 2\fl[n] + 5$, for $n \geq 3$. This appears to be incorrect. In particular, in \hyperref[fig-err3]{Figure \ref{fig-err4}} we display a drawing of $G^6_{133} + C_3$ which has fifteen crossings, rather than sixteen as suggested by \cite{vijayaetal2016}.

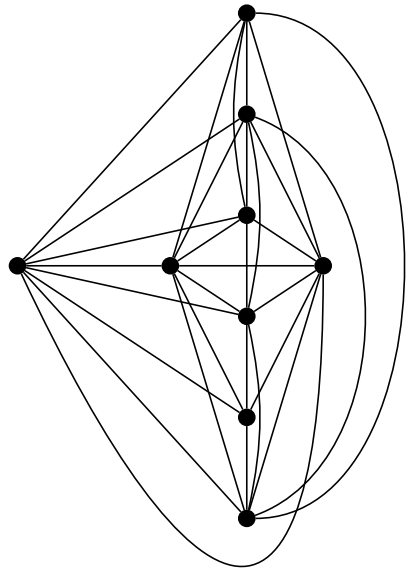
\begin{figure}[h!]\begin{center}\begin{tikzpicture}[largegraph,scale=0.67]
\foreach \i in {0,...,5}{
  \node (\i) at (0,10-2*\i) {};
}
\foreach \j in {6,7,8}{
  \node (\j) at (-4.5+3*\j-18,5) {};
  \foreach \i in {0,...,5}{
    \draw (\j) -- (\i);
  }
}
\draw (6) -- (7) -- (8);
\draw[out=295,in=270,looseness=3.5] (6) to (8);
\draw (0) -- (1) -- (2) -- (3) -- (4) -- (5);
\draw[bend right=12] (0) to (2);
\draw[bend left=12] (1) to (3);
\draw[bend left=12] (3) to (5);
\draw[bend left=70] (1) to (5);
\draw[bend left=90] (0) to (5);
\end{tikzpicture}\caption{A drawing of $G^6_{133} + C_3$, with only fifteen crossings.}\label{fig-err4}\end{center}\end{figure}

\item In 2016, Zhou and Li \cite{zhouli2016} considered the join product of $G^6_{111}$ with discrete graphs, paths and cycles. They claimed that $cr(G^6_{111} + D_n) = Z(6,n) + 2\fl[n]$, $cr(G^6_{111} + P_{n-1}) = Z(6,n) + 2\fl[n] + 1$ and $cr(G^6_{111} + C_n) = Z(6,n) + 2\fl[n] + 3$. These appear to be incorrect, and Crossing Number Web Compute shows that the minimal counterexamples are $cr(G^6_{111} + D_1) = 1$, $cr(G^6_{111} + P_1) = 4$ and $cr(G^6_{111} + C_3) = 13$, rather than the zero, three and eleven crossings suggested by \cite{zhouli2016} respectively.
\end{itemize}

The proof files for the crossing numbers determined by Crossing Number Web Compute are available as follows:

$cr(G^6_{131} \Box P_1) = 4$: \href{http://crossings.uos.de/job/KN4vnYbb797WV0C5TQdTNA}{http://crossings.uos.de/job/KN4vnYbb797WV0C5TQdTNA}

$cr(G^6_{111} + D_1) = 1$: \href{http://crossings.uos.de/job/MuRSG12mxlzDmJz4cc9eFA}{http://crossings.uos.de/job/MuRSG12mxlzDmJz4cc9eFA}

$cr(G^6_{111} + P_1) = 4$: \href{http://crossings.uos.de/job/dhr5Vl\_8pWDL-jAARvBNzw}{http://crossings.uos.de/job/dhr5Vl\_8pWDL-jAARvBNzw}

$cr(G^6_{111} + C_3) = 13$: \href{http://crossings.uos.de/job/SHHUcOtUyGpQFNY-Nqmwxg}{http://crossings.uos.de/job/SHHUcOtUyGpQFNY-Nqmwxg}


\addcontentsline{toc}{section}{References}
\begin{thebibliography}{100}

\bibitem{abregoetal2015}Bernado M. \'{A}brego, Oswin Aichholzer, Silvia Fern\'{a}ndez-Merchant, Thomas Hackl, J\"{u}rgen Pammer, Alexander Pilz, Pedro Ramos, Gelasio Salazar and Birgit Vogtunhuber,  \newblock All Good Drawings of Small Complete Graphs,  \newblock In: {\em Proc. 31st European Workshop on Computational Geometry (EuroCG)} (2015), pp. 57--60.
\bibitem{ackerman2019}Eyal Ackerman,  \newblock On topological graphs with at most four crossings per edge,  \newblock {\em Comp. Geom.}, 85 (2019), 101574.
\bibitem{adamssonrichter2004}Jay Adamsson and R. Bruce Richter,  \newblock Arrangements, circular arrangements and the crossing number of $C_7 \times C_n$,  \newblock {\em J. Combin. Th. Series B}, 90(1) (2004), 21--39.
\bibitem{aignerziegler2010} Martin Aigner and G\"{u}nter M. Ziegler,  \newblock {\em Proofs from THE BOOK,} Springer-Verlag Berlin Heidelberg (2010), pp. 261--268.
\bibitem{ajtaietal1982} Mikl\'{o}s Ajtai, Va\v{s}ek Chv\'{a}tal, Monroe M. Newborn and Endre Szemer\'{e}di,  \newblock Crossing-free subgraphs,  \newblock {\em North-Holland Math. Stud.}, 60(C) (1982), 9--12.
\bibitem{akerskrishnamurthy1989} Sheldon B. Akers and Balakrishnan Krishnamurthy,  \newblock A group-theoretic model for symmetric interconnection networks,  \newblock {\em IEEE Trans. Comput.}, 38 (1989), 555--566.
\bibitem{albertsonetal2010} Michael O. Albertson, Daniel W. Cranston and Jacob Fox,  \newblock Crossings, colorings and cliques,  \newblock {\em Elec. J. Combin.}, 16, Note 45 (2010).
\bibitem{andersonetal1996} Mark Anderson, R. Bruce Richter and Peter Rodney,  \newblock The Crossing Number of $C_6 \times C_6$,  \newblock {\em Congr. Numer.}, 118 (1996), 97--107.
\bibitem{andersonetal1996_2} Mark Anderson, R. Bruce Richter and Peter Rodney,  \newblock The Crossing Number of $C_7 \times C_7$,  \newblock {\em Congr. Numer.}, 125 (1996), 97--117.
\bibitem{asano1986} Kouhei Asano,  \newblock The crossing number of $K_{1,3,n}$ and $K_{2,3,n}$,  \newblock {\em J. Graph Th.}, 10 (1986), 1--8.
\bibitem{balogh2019} J\'ozsef Balogh, Bernard Lidick\'y and Gelasio Salazar, \newblock Closing in On Hill's Conjecture, \newblock {\em SIAM J. Discr. Math.}, 33(3) (2019), 1261--1276.
\bibitem{barattoth2010} J\'{a}nos B\'{a}rat and G\'{e}za T\'{o}th,  \newblock Towards the Albertson conjecture,  \newblock {\em Elec. J. Combin.}, 17, Note 73 (2010).
\bibitem{beinekeringeisen1980} Lowell W. Beineke and Richard D. Ringeisen,  \newblock On the crossing numbers of products of cycles and graphs of order four,  \newblock {\em J. Graph Th.}, 4(2) (1980), 145--155.
\bibitem{beinekewilson2010} Lowell W. Beineke and Robin Wilson,  \newblock The Early History of the Brick Factory Problem,  \newblock {\em Math. Int.}, 32(2) (2010), 41--48.
\bibitem{benes1964} Vaclav E. Bene\v{s},  \newblock Permutation groups, complexes, and rearrangeable connecting networks,  \newblock {\em Bell System Tech. J.}, 43(4) (1964), 1619--1640.
\bibitem{bereznystas2017} \v{S}tefan Bere\v{z}n\'{y} and Michal Sta\v{s},  \newblock On the crossing number of the join of five vertex graph $G$ with the discrete graph $D_n$, {\em Acta Elec. Inf.}, 17(3) (2017), 27--32.
\bibitem{bereznystas2019} \v{S}tefan Bere\v{z}n\'{y} and Michal Sta\v{s},  \newblock Cyclic permutations and crossing numbers of join products of two symmetric graphs of order six, {\em Carpathian J. Math.}, 35(2) (2019), 137--146.
\bibitem{bereznystas2020} \v{S}tefan Bere\v{z}n\'{y} and Michal Sta\v{s},  \newblock On the crossing number of join of the wheel on six vertices with the discrete graph, {\em Carpathian J. Math.}, 36(3) (2020), 383--392.
\bibitem{bokal2007} Drago Bokal,  \newblock On the crossing number of Cartesian products with paths,  \newblock {\em J. Comb. Th. Series B}, 97(3) (2007), 381--384.
\bibitem{bokal2007_2} Drago Bokal,  \newblock On the crossing numbers of Cartesian products with trees,  \newblock {\em J. Graph Th.}, 56 (20007), 287--300.
\bibitem{bokaletal2019} Drago Bokal, Zden\u{e}k Dvo\u{r}\'{a}k, Petr Hlin\u{e}n\'{y}, Jes\'{u}s Lea\~{o}s, Bojan Mohar and Tilo Wiedera, \newblock Bounded degree conjecture holds precisely for $c$-crossing critical graphs with $c\leq 12$, (2019), \newblock \arxiv{1903.05363}.
\bibitem{chialee2015} Gek L. Chia and Chan L. Lee,  \newblock Crossing Numbers of Nearly Complete Graphs and Nearly Complete Bipartite Graphs,  \newblock {\em Ars Combin.}, 121 (2015), 437--446.
\bibitem{chimanietal2014} Markus Chimani, Carsten Gutwenger, Michael J\"{u}nger, Gunnar W Klau, Karsten Klein and Petra Mutzel,  \newblock The Open Graph Drawing Framework (OGDF),  \newblock {\em Handbook of Graph Drawing and Visualization}, (2011), pp. 543--569.
\bibitem{chimaniwiedera2016} Markus Chimani and Tilo Wiedera,  \newblock An ILP-based Proof System for the Crossing Number Problem,  \newblock In: {\em LIPIcs-Leibniz International Proceedings in Informatics}, vol. 57. Schloss Dagstuhl-Leibniz-Zentrum fuer Informatik, (2016).
\bibitem{chimaniwiederasite2016} Markus Chimani and Tilo Wiedera,  \newblock Crossing Number Web Compute,  \newblock \href{http://crossings.uos.de}{http://crossings.uos.de}, (2016).
\bibitem{choudumsunitha2002} Sheshayya A. Choudum and Vadivel Sunitha,  \newblock Augmented Cubes,  \newblock {\em Networks}, 40 (2002), 71--84.
\bibitem{christianetal2013} Robin Christian, R. Bruce Richter and Gelasio Salazar,  \newblock Zarankiewicz's conjecture is finite for each fixed $m$,  \newblock {\em J. Comb. Th. Ser. B}, 103(2) (2013), 237--247.
\bibitem{cimikowski1996} Robert Cimikowski,  \newblock Topological Properties of Some Interconnected Network Graphs,  \newblock {\em Congr. Numer.}, 121 (1996), 19--32.
\bibitem{cimikowski1998} Robert Cimikowski,  \newblock Crossing number bounds for the mesh of trees,  \newblock {\em Congr. Numer.}, 129 (1998), 107--116.
\bibitem{cimikowskivrto2003} Robert Cimikowski and Imrich Vrt'o,  \newblock Improved bounds for the crossing number of the mesh of trees,  \newblock {\em J. Interconnect. Net.}, 4(1) (2003), 17--35.
\bibitem{clancyetal2019} Kieran Clancy, Michael Haythorpe and Alex Newcombe,  \newblock An effective crossing minimisation heuristic based on star insertion,  \newblock {\em J. Graph Alg. App.}, 23(2) (2019), 135--166.
\bibitem{clancyetaltoappear} Kieran Clancy, Michael Haythorpe and Alex Newcombe,  \newblock On the Crossing Numbers of Cartesian Products of Small Graphs with Paths, Cycles and Stars,  \newblock{J. Comb. Math. Comb. Comp.}, to appear.
\bibitem{coxeter1950} Harold S.M. Coxeter,  \newblock Self-Dual Configurations and Regular Graphs,  \newblock {\em Bull. Amer. Math. Soc.}, 56 (1950), 413--455.
\bibitem{cruzjapson2015} Caselyn D.C. Cruz and Danica M. Japson,  \newblock The cross
ing number of the generalized Petersen graph $P(17,K)$, where $K \leq \lfloor \frac{17}{2} \rfloor$,  \newblock {\em Undergraduate Thesis}, Polytechnic University of the Phillipines, (2015). $\hfill \Asterisk$
\bibitem{deklerketal2006} Etienne de Klerk, John Maharry, Dmitrii V. Pasechnik, R. Bruce Richter and Gelasio Salazar,  \newblock Improved bounds for the crossing numbers of $K_{m,n}$ and $K_n$,  \newblock {\em SIAM J. Discr. Math.}, 20 (2006), 189--202.
\bibitem{deklerketal2007} Etienne de Klerk, Dmitrii V. Pasechnik and Alexander Schrijver,  \newblock Reduction of symmetric semidefinite programs using the regular $*$-representation,  \newblock {\em Math. Prog. Ser. B}, 109 (2007), 613--624.
\bibitem{deanrichter1995}Alice M. Dean and R. Bruce Richter,  \newblock The crossing number of $C_4 \times C_4$,  \newblock {\em J. Graph Th.}, 19(1) (1995), 125--129.
\bibitem{dinghuang2018} Zongpeng Ding and Yuanqiu Huang,  \newblock The crossing numbers of join of some graphs with $n$ isolated vertices,  \newblock {\em Disc. Math. Graph Th.}, 38 (2018), 899--909.
\bibitem{dingetal2018} Zongpeng Ding, Yuanqiu Huang and Zhangdong Ouyang,  \newblock The crossing numbers of Cartesian products of paths with some graphs,  \newblock {\em Ars Combin.}, 141 (2018), 101--110.
\bibitem{drazenska2011} Em\'{i}lia Dra\v{z}ensk\'{a},  \newblock The crossing number of $G \Box C_n$ for the graph $G$ on six vertices,  \newblock {\em Math. Slovaca}, 61(5) (2011), 675--686.
\bibitem{drazenskaklesc2007} Em\'{i}lia Dra\v{z}ensk\'{a} and Mari\'{a}n Kle\v{s}\v{c},  \newblock The crossing numbers of products of cycles with 6-vertex trees,  \newblock {\em Tatra Mt. Math. Publ.}, 36 (2007), 109--119.
\bibitem{drazenskaklesc2008} Em\'{i}lia Dra\v{z}ensk\'{a} and Mari\'{a}n Kle\v{s}\v{c},  \newblock The crossing numbers of products of the graph $K_{2,2,2}$ with stars,  \newblock {\em Carpathian J. Math.}, 24(3) (2008), 327--331.
\bibitem{drazenskaklesc2011} Em\'{i}lia Dra\v{z}ensk\'{a} and Mari\'{a}n Kle\v{s}\v{c},  \newblock On the crossing numbers of $G \Box C_n$ for graphs $G$ on six vertices,  \newblock {\em Disc. Math. Graph Th.}, 31(2) (2011), 239--252.
\bibitem{eggletonguy1970} Roger B. Eggleton and Richard K. Guy,  \newblock The crossing number of the $n$-cube,  \newblock {\em Notices Amer. Math. Soc.}, 17(5) (1970), 757.
\bibitem{erdosguy1973} Paul Erd\H{o}s and Richard K. Guy,  \newblock Crossing Number Problems,  \newblock {\em Amer. Math. Monthly}, 80(1) (1973), 52--58.
\bibitem{exooetal1981} Geoffrey Exoo, Frank Harary and Jerald Kabell,  \newblock The crossing numbers of some generalized Petersen graphs,  \newblock {\em Math. Scand.}, 48 (1981), 184--188.
\bibitem{fariadefigueiredo2000}Luerbio Faria, Celina M.H. de Figueiredo,  \newblock On Eggleton and Guy's Conjectured Upper Bound for the Crossing Number of the $n$-cube,  \newblock {\em Math. Slovaca}, 50(3) (2000), 271--287.
\bibitem{fariaetal2008}Luerbio Faria, Celina M.H. de Figueiredo, Ondrej S\'{y}kora, and Imrich Vrt'o,  \newblock An improved upper bound on the crossing number of the hypercube,  \newblock {\em J. Graph Th.}, 59(2) (2008), 145--161.
\bibitem{fiorini1986} Stanley Fiorini,  \newblock On the crossing number of generalized Petersen graphs,  \newblock {\em North-Holland Math. Stud.} 123 (1986), 225--241.
\bibitem{fiorinigauci2003} Stanley Fiorini and John Baptist Gauci,  \newblock The crossing number of the generalized Petersen graph $P(3k,k)$,  \newblock {\em Math. Bohemica}, 128(4) (2003), 337--347.
\bibitem{foleyetal2002} Abbie Foley, Rachel Krieger, Adrian Riskin and Isabelle Stanton,  \newblock The Crossing Numbers of Some Twisted Toroidal Grid Graphs,  \newblock {\em Bull. ICA}, 36 (2002), 80--88.
\bibitem{gareyjohnson1983} Michael R. Garey and David S. Johnson,  \newblock Crossing number is NP-complete,  \newblock {\em SIAM J. Algebraic Disc. Meth.}, 4(3) (1983), 312--316.
\bibitem{gaucixuereb2019} John Baptist Gauci and Cheryl Zerafa Xuereb,  \newblock A note on isomorphic generalized Petersen graphs with an application to the crossing number of $GP[3k-1,k]$ and $GP[3k+1,k]$,  \newblock {\em Disc. Math. Letters}, 2 (2019), 44--46.
\bibitem{gethneretal2017} Ellen Gethner, Leslie Hogben, Bernard Lidick\'{y}, Florian Pfender, Amanda Ruiz and Michael Young,  \newblock On Crossing Numbers of Complete Tripartite and Balanced Complete Multipartite Graphs,  \newblock {\em J. Graph Th.}, 84(4) (2017), 552--565.
\bibitem{glebskysalazar2004} Lev Yu Glebsky and Gelasio Salazar,  \newblock The crossing number of $C_m \times C_n$ is as conjectured for $n \geq m(m+1)$,  \newblock {\em J. Graph Th.}, 47(1) (2004), 53--72.
\bibitem{gutwenger2010} Carsten Gutwenger,  \newblock {\em Application of SPQR-Trees in the Planarization Approach for Drawing Graphs},  \newblock PhD Thesis, Technical University of Dortmund, (2010).
\bibitem{guy1960} Richard K. Guy,  \newblock A combinatorial problem,  \newblock {\em Bull. Malayan Math. Soc.}, 7 (1960), 68--72.
\bibitem{guy1969} Richard K. Guy,  \newblock The decline and fall of Zarankiewicz's theorem,  \newblock {\em Proof techniques in Graph Theory}, F. Harary, ed., Academic Press, New York (1969), pp. 63--69.
\bibitem{guy1972} Richard K. Guy,  \newblock Crossing numbers of graphs,  \newblock In: {\em Graph Theory and Applications}, volume 303 of {\em Lecture Notes in Mathematics}, Spinger-Verlag, Berlin-Heidelberg-New York, (1972), pp.111--124.
\bibitem{guyharary1967} Richard K. Guy and Frank Harary,  \newblock On the M\"{o}bius ladders,  \newblock {\em Canad. Math. Bull.} 10 (1967), 493--496.
\bibitem{guyhill1973} Richard K. Guy and Anthony Hill,  \newblock The crossing number of the complement of a circuit,  \newblock {\em Disc. Math.}, 5 (1973), 335--344.
\bibitem{harary1969} Frank Harary,  \newblock {\em Graph Theory},  \newblock Addison-Wesley, Reading, MA, (1969).
\bibitem{harary1963} Frank Harary and Anthony Hill, \newblock On the Number of Crossings in a Complete Graph, \newblock {\em Proceedings of the Edinburgh Mathematical Society}, 13(4) (1963), 333--338.
\bibitem{hararykainen1993} Frank Harary and Paul C. Kainen,  \newblock The Cube of a Path is Maximal Planar,  \newblock {\em Bull. ICA}, 7 (1993), 55--56.
\bibitem{hararyetal1999} Frank Harary, Paul C. Kainen and Adrian Riskin,  \newblock Every Graph of Cyclic Bandwidth 3 is Toroidal,  \newblock {\em Bull. ICA}, 27 (1999), 81--84.
\bibitem{hararykainen1973} Frank Harary, Paul C. Kainen and Allen J. Schwenk,  \newblock Toroidal graphs with arbitrarily high crossing numbers,  \newblock {\em Nanta. Math.}, 6(1) (1973), 58--67.
\bibitem{harborth1970} Heiko Harborth, \newblock \"{U}ber die Kreuzungsahl vollst\"{a}ndiger, $n$-geteilter Graphen, \newblock {\em Mathematische Nachrichten}, 48 (1971), 179--188.
\bibitem{haythorpenewcombe2019} Michael Haythorpe and Alex Newcombe,  \newblock On the Crossing Number of the Cartesian Product of a Sunlet Graph and a Star Graph,  \newblock {\em Bull. Austral. Math. Soc.}, 100 (2019), 5--12.
\bibitem{he2008} Peiling He, The Crossing Number of the Complete Bipartite Graph by Deleting Two Edges $e_1$ and $e_2$, {\em J. Hengyang Norm. Uni.}, 29(6) (2008), 23--25. $\hfill \Asterisk$
\bibitem{hehuang2007} Peiling He and Yuanqiu Huang,  \newblock On the Crossing Numbers of Cartesian Product of a 6-vertice Graph and Paths,  \newblock {\em J. Hunan Inst. Hum. Sci. Tech.}, 2007(6) (2007), 3--5. $\hfill \Asterisk$
\bibitem{hehuang2007_2} Peiling He and Yuanqiu Huang,  \newblock The crossing number of $W_4 \times S_n$,  \newblock {\em J. Zhengzhou Uni. (Nat. Sci. Ed.)}, 39(4) (2007), 14--18. $\hfill \Asterisk$
\bibitem{heetal2011} Peiling He, Zhijun Luo and Yuanqiu Huang,  \newblock Crossing Number of Several Complete Bipartite Graphs by Deleting One Edge,  \newblock {\em J. Henan Norm. Uni. Nat. Sci.}, 39(2) (2011), 24--26. $\hfill \Asterisk$
\bibitem{he2011} Xioanian He,  \newblock The Crossing Number of Cartesian Products of Stars with 5-vertex Graphs II,  \newblock {\em Chin. Quart. J. Math.}, 26(4) (2011), 563--567.
\bibitem{hehuang2005} Xiaonian He and Yuanqiu Huang,  \newblock The Crossing Number of a Class of Cartesian Products,  {\em J. Jishou Uni. Nat. Sci. Ed.}, 26(1) (2005), 8--11. $\hfill \Asterisk$
\bibitem{ho2004} Pak Tung Ho,  \newblock The crossing number of $K_{1,5,n}$, $K_{2,4,n}$ and $K_{3,3,n}$,  \newblock {\em Int. J. Pure Appl. Math.}, 17(4) (2004), 491--515.
\bibitem{ho2007} Pak Tung Ho,  \newblock The crossing number of $C(3k+1;\{1,k\})$,  \newblock {\em Disc. Math.}, 307(22) (2007), 2771--2774.
\bibitem{ho2008} Pak Tung Ho,  \newblock The crossing number of $K_{1,m,n}$,  \newblock {\em Disc. Math.}, 308(24) (2008), 5996--6002.
\bibitem{ho2008_2} Pak Tung Ho,  \newblock The crossing number of $K_{2,2,2,n}$,  \newblock {\em Far East J. Appl. Math.}, 30(1) (2008), 43--69.
\bibitem{ho2009} Pak Tung Ho,  \newblock On the crossing number of some complete multipartite graphs,  \newblock {\em Utilitas Math.}, 79 (2009), 125--143.
\bibitem{ho2013} Pak Tung Ho,  \newblock The Crossing Number of $K_{2,4,n}$,  \newblock {\em Ars Comb.}, 109 (2013), 527--537.
\bibitem{hsiehlin2016} Sunyuan Hsieh and Chengchian Lin, \newblock The Crossing Number of Join Products of kth Power of Path $P_m$ with Isolated Vertices and Path $P_n$, In: {\em Computer Symposium (ICS) 2016 International}, IEEE, (2016), pp. 62--67.
\bibitem{huangwang2010} Yuanqiu Huang and Jing Wang, \newblock Survey of the crossing number of graphs,  \newblock {\em J. East China Normal Uni. Nat. Sci.}, 2010(3) (2010), 68--80.
\bibitem{huangzhao2006} Yuanqiu Huang and Tinglei Zhao,  \newblock On the Crossing Number of the Complete Tripartite $K_{1,6,n}$,  \newblock {\em Acta Math. Appl. Sinica}, 29(6) (2006), 1046--1053.
\bibitem{huangzhao2006_2} Yuanqiu Huang and Tinglei Zhao,  \newblock On the Crossing Number of the Complete Tripartite $K_{1,8,n}$,  \newblock {\em Acta Math. Sinica}, 26(7) (2006), 1115--1122.
\bibitem{huangzhao2008} Yuanqiu Huang and Tinglei Zhao,  \newblock The crossing number of $K_{1,4,n}$,  \newblock {\em Disc. Math.}, 308(9) (2008), 1634--1638.
\bibitem{isaacs1975} Rufus Isaacs,  \newblock Infinite Families of Nontrivial Trivalent Graphs Which Are Not Tait Colorable,  \newblock {\em Amer. Math. Monthly}, 82 (1975), 221--239.
\bibitem{jendrolscerbova1982} Stanislav Jendrol' and M\'{a}ria \v{S}cerbov\'{a},  \newblock On the crossing numbers of $S_m \times P_n$ and $S_m \times C_n$. {\em Casopis pro pestov\'{a}ni matematiky}, 107 (1982), 225--230.
\bibitem{klavzarmohar2005} Sandi Klav\v{z}ar and Bojan Mohar,  \newblock Crossing Numbers of Sierpi\'{n}ski-Like Graphs,  \newblock {\em J. Graph Th.}, 50(3) (2005), 186--198.
\bibitem{kleitman1971} Daniel J. Kleitman,  \newblock The crossing number of $K_{5,n}$,  \newblock {\em J. Combin. Th.}, 9 (1971), 315--323.
\bibitem{klesc1991} Mari\'{a}n Kle\v{s}\v{c},  \newblock On the crossing numbers of Cartesian products of stars and paths or cycles,  \newblock {\em Math. Slovaca}, 41(2) (1991), 113--120.
\bibitem{klesc1994} Mari\'{a}n Kle\v{s}\v{c},  \newblock The crossing numbers of products of paths and stars with 4-vertex graphs,  \newblock {\em J. Graph Th.}, 18(6) (1994), 605--614.
\bibitem{klesc1995} Mari\'{a}n Kle\v{s}\v{c},  \newblock The crossing numbers of certain Cartesian products,  \newblock {\em Disc. Math. Graph Th.}, 15(1) (1995), 5-10.
\bibitem{klesc1996} Mari\'{a}n Kle\v{s}\v{c},  \newblock The crossing number of $K_{2,3} \times P_n$ and $K_{2,3} \times S_n$,  \newblock {\em Tatra Mountains Math. Publ.}, 9 (1996), 51--56.
\bibitem{klesc1999} Mari\'{a}n Kle\v{s}\v{c},  \newblock The crossing number of $K_5 \times P_n$,  \newblock {\em Tatra Mountains Math. Pub.}, 18 (1999), 63--68.
\bibitem{klesc1999_2} Mari\'{a}n Kle\v{s}\v{c},  \newblock The crossing numbers of products of a 5-vertex graph with paths and cycles,  \newblock {\em Disc. Math. Graph Th.}, 19(1) (1999), 59--69.
\bibitem{klesc2001} Mari\'{a}n Kle\v{s}\v{c},  \newblock On the crossing numbers of products of stars and graphs of order five,  \newblock {\em Graphs Comb.}, 17(2) (2001), 289--294.
\bibitem{klesc2001_2} Mari\'{a}n Kle\v{s}\v{c},  \newblock The crossing numbers of Cartesian products of paths with 5-vertex graphs,  \newblock {\em Disc. Math.}, 233(1--3) (2001), 353--359.
\bibitem{klesc2002} Mari\'{a}n Kle\v{s}\v{c},  \newblock The crossing number of $K_{2,3} \times C_3$,  \newblock {\em Disc. Math}, 251(1--3) (2002), 109--117.
\bibitem{klesc2005} Mari\'{a}n Kle\v{s}\v{c},  \newblock Some crossing numbers of products of cycles,  \newblock {\em Disc. Math. Graph Th.}, 25(1--2) (2005), 197--210.
\bibitem{klesc2007} Mari\'{a}n Kle\v{s}\v{c},  \newblock The join of graphs and crossing numbers,  \newblock {\em Elec. Notes Disc. Math.}, 28 (2007), 349--355.
\bibitem{klesc2009} Mari\'{a}n Kle\v{s}\v{c},  \newblock On the crossing numbers of Cartesian products of stars and graphs on five vertices,  \newblock In: Fiala J., Kratochv\'{\i}l J., Miller M. (eds) {\em Combinatorial Algorithms. IWOCA 2009. Lecture Notes in Computer Science, vol 5874}, Springer, Berlin, Heidelberg, (2009), pp. 324--333.
\bibitem{klesc2010} Mari\'{a}n Kle\v{s}\v{c},  \newblock The crossing numbers of join of the special graph on six vertices with path and cycle,  \newblock {\em Disc. Math.}, 310(9) (2010), 1475--1481.
\bibitem{klesckocurova2007} Mari\'{a}n Kle\v{s}\v{c} and Anna Koc\'{u}rov\'{a},  \newblock The crossing numbers of products of 5-vertex graphs with cycles,  \newblock {\em Disc. Math.}, 307(11-12) (2007), 1395--1403.
\bibitem{klesckravecova2008} Mari\'{a}n Kle\v{s}\v{c} and Daniela Kravecov\'{a},  \newblock The crossing number of $P^2_5 \Box C_n$,  \newblock {\em Creat. Math. Inf.}, 17(3) (2008), 431--438.
\bibitem{klesckravecova2012} Mari\'{a}n Kle\v{s}\v{c} and Daniela Kravecov\'{a},  \newblock The crossing number of $P^2_n \Box C_3$,  \newblock {\em Disc. Math.}, 312(14) (2012), 2096--2101.
\bibitem{klescetal2011} Mari\'{a}n Kle\v{s}\v{c}, Daniela Kravecov\'{a} and Jana Petrillov\'{a},  \newblock The Crossing Numbers of Join of Special Graphs,  \newblock In: {\em Electrical Engineering and Informatics II: Proceedings of the Faculty of Electrical Engineering and Informatics of the Technical University of Ko\v{s}ice (2011)}, (2011), pp. 522--527.
\bibitem{klescetal2014} Mari\'{a}n Kle\v{s}\v{c}, Daniela Kravecov\'{a} and Jana Petrillov\'{a},  \newblock On the crossing numbers of Cartesian products of paths with special graphs,  \newblock {\em Carpathian J. Math.}, 30(3) (2014), 317--325.
\bibitem{klescpetrillova2013} Mari\'{a}n Kle\v{s}\v{c} and Jana Petrillov\'{a},  \newblock On the optimal drawings of products of paths with graphs,  \newblock {\em Acta Elec. Inf.}, 13(3) (2013), 56--61.
\bibitem{klescpetrillova2013_2} Mari\'{a}n Kle\v{s}\v{c} and Jana Petrillov\'{a},  \newblock The crossing numbers of products of paths with graphs of order six,  \newblock {\em Disc. Math. Graph Th.}, 33(3) (2013), 571--582.
\bibitem{klescetal2013} Mari\'{a}n Kle\v{s}\v{c}, Jana Petrillov\'{a} and Mat\'{u}\v{s} Valo,  \newblock Minimal number of crossings in strong product of paths,  \newblock {\em Carpathian J. Math.}, 29(1) (2013), 27--32.
\bibitem{klescetal2017} Mari\'{a}n Kle\v{s}\v{c}, Jana Petrillov\'{a} and Mat\'{u}\v{s} Valo,  \newblock On the crossing numbers of Cartesian products of wheels and trees,  \newblock {\em Disc. Math. Graph Th.}, 37(2) (2017), 399--413.
\bibitem{klescetal1996} Mari\'{a}n Kle\v{s}\v{c}, R. Bruce Richter and Ian Stobert,  \newblock The crossing number of $C_5 \times C_n$,  \newblock {\em J. Graph Th.}, 22(3) (1996), 239--243.
\bibitem{klescschrotter2011} Mari\'{a}n Kle\v{s}\v{c} and \v{S}tefan Schr\"{o}tter,  \newblock The crossing numbers of join products of paths with graphs of order four,  \newblock {\em Disc. Math. Graph Th.}, 31(2) (2011), 321--331.
\bibitem{klescschrotter2012} Mari\'{a}n Kle\v{s}\v{c} and \v{S}tefan Schr\"{o}tter,  \newblock The Crossing Numbers of Join of Paths and Cycles with Two Graphs of Order Five,  \newblock In: {\em Mathematical Modeling and Computational Science}, Springer, Berlin, Heidelberg, (2012), pp. 160--167.
\bibitem{klescschrotter2013} Mari\'{a}n Kle\v{s}\v{c} and \v{S}tefan Schr\"{o}tter,  \newblock On the crossing numbers of Cartesian products of stars and graphs of order six,  \newblock {\em Disc. Math. Graph Th.}, 33(3) (2013), 583--597.
\bibitem{klescstastoappear} Mari\'{a}n Kle\v{s}\v{c} and Michal Sta\v{s},  \newblock Cyclic permutations in determining crossing numbers,  \newblock {\em Disc. Math. Graph Th.}, to appear.
\bibitem{klescvalo2012} Mari\'{a}n Kle\v{s}\v{c} and Mat\'{u}\v{s} Valo,  \newblock Minimum crossings in join of graphs with paths and cycles,  \newblock {\em Acta Elec. Inf.}, 12(3) (2012), 32--37.
\bibitem{kravecova2012} Daniela Kravecov\'{a},  \newblock The crossing number of $P^2_n \Box C_4$,  \newblock {\em Acta Elec. Inf.}, 12(3) (2012), 42--46.
\bibitem{kravecova2012_2} Daniela Kravecov\'{a},  \newblock The crossing number of $P_5^2 \times P_n$,  \newblock {\em Creat. Math. Inf.}, 28(1) (2012), 49--56.
\bibitem{leighton1983} Frank T. Leighton,  \newblock {\em Complexity Issues in VLSI: Optimal Layouts for the Shuffle-Exchange Graph and Other Networks}  \newblock MIT Press, Cambridge, MA, (1983).
\bibitem{leighton1984} Frank T. Leighton,  \newblock New Lower Bound Techniques for VLSI,  \newblock {\em Math. Sys. Th.}, 17(1) (1984), 47--70.
\bibitem{lietal2008} Bo Li, Jing Wang and Yuanqiu Huang,  \newblock On the crossing number of the join of some 5-vertex graphs and $P_n$,  \newblock {\em Int. J. Math. Comb.}, 2 (2008), 70--77. $\hfill \Asterisk$
\bibitem{lietal2008_2} Bo Li, Jing Wang and Yuanqiu Huang,  \newblock On the crossing number of the join of some 6-vertex graphs and $P_n$,  \newblock {\em J. Jishou Uni. Nat. Sci. Ed. (China)}, 29(6) (2008), 29--35. $\hfill \Asterisk$
\bibitem{lietal2008_3} Bo Li, Lixi Zhang and Yuanqiu Huang,  \newblock The crossing number of Cartesian product of star with a 6-vertex graph,  \newblock {\em J. Hunan Uni. Arts Sci. Nat. Sci. Ed.}, 20(2) (2008), 6--11. $\hfill \Asterisk$
\bibitem{lietal2009} Bo Li, Lixi Zhang and Yuanqiu Huang,  \newblock Crossing Number of Products of Several 6-Vertex Graphs with Star,  \newblock {\em J. Shantou Uni. Nat. Sci. Ed.}, 24(4) (2009), 4--13. $\null\hfill \Asterisk$
\bibitem{li2014} Liping Li,  \newblock The Crossing Numbers of Join of a 5-Vertex Graph with Path and Cycle,  \newblock {\em Math. Prac. Th.}, 44(11) (2014), 203--211. $\hfill \Asterisk$
\bibitem{li2013} Min Li,  \newblock Crossing Numbers of Join of the Graph on Five Vertices with $n$ Isolated Vertices and Paths,  \newblock {\em J. Hubei Uni. Arts Sci.}, 34(11) (2013), 15--17. $\hfill \Asterisk$
\bibitem{li2013_2} Min Li,  \newblock The Crossing Numbers of Join of Some 5-Vertex Graphs with Paths and Cycles,  \newblock {\em J. Henan Norm. Uni. Nat. Sci. Ed.}, 41(4) (2013), 40--44. $\hfill \Asterisk$
\bibitem{li2015} Min Li,  \newblock The crossing numbers of the join of a 5-vertex graph with vertex, path and cycle,  \newblock {\em J. Yangzhou Uni. Nat. Sci. Ed.}, 18(1) (2015), 4--8. $\hfill \Asterisk$
\bibitem{linetal2005} Xiaohui Lin, Yuansheng Yang, Jianguo L\"{u} and Xin Hao,  \newblock The crossing number of $C(mk; \{1,k\})$,  \newblock {\em Graphs Combin.}, 21(1) (2005), 89--96.
\bibitem{linetal2006} Xioahui Lin, Yuansheng Yang, Jianguo L\"{u} and Xin Hao,  \newblock The Crossing Number of $C(n;\{1,[n/2]-1\})$,  \newblock {\em Util. Math.}, 71 (2006), 245--255.
\bibitem{linetal2009} Xiaohui Lin, Yuansheng Yang, Wenping Zheng, Lei Shi and Weiming Lu,  \newblock The crossing numbers of generalized Petersen graphs with small order,  \newblock {\em Disc. Appl. Math.}, 157(5) (2009), 1016--1023.
\bibitem{liuetal2012} Wei Liu, Zihan Yuan and Shengxiang L\"{u},  \newblock The Crossing Number of $c(7,2) \setminus \{e_1,e_2\} \times P_n$,  \newblock {\em J. Hunan Inst. Sci. Tech.}, 25(1) (2012), 16--19. $\hfill \Asterisk$
\bibitem{luetal2014} Bo L\"{u}, Xirong Xu, Yuansheng Yang, Ke Zhang and Baigong Zheng,  \newblock Crossing number of Star graph $S_4$,  \newblock {\em J. Dalian Uni. Tech.}, 54(4) (2014), 469--476.
\bibitem{luhuang2008} Shengxiang L\"{u} and Yuanqiu Huang,  \newblock The crossing number of $K_5 \times S_n$,  \newblock {\em J. Math. Res. Exp.}, 28(3) (2008), 445--459.
\bibitem{luhuang2010} Shengxiang L\"{u} and Yuanqiu Huang,  \newblock On the crossing number of $K_{2,4} \times S_n$,  \newblock {\em J. Sys. Sci. Math. Sci.}, 30(7) (2010), 929--935. $\hfill \Asterisk$
\bibitem{luhuang2011} Shengxiang L\"{u} and Yuanqiu Huang,  \newblock The crossing number of $K_5 \setminus e \times S_n$,  \newblock {\em J. Hunan Uni. Art. Sci. Nat. Sci.}, 23(1) (2011), 1--5. $\hfill \Asterisk$
\bibitem{ma2017} Dengju Ma,  \newblock The crossing number of the strong product of two paths,  \newblock {\em Austral. J. Combin.}, 68(1) (2017), 35--47.
\bibitem{maetal2005} Dengju Ma, Han Ren and Junjie Lu,  \newblock Crossing Number of the Generalized Petersen Graph $G(2m+1,m)$,  \newblock {\em J. East China Norm. Uni. (Nat. Sci. Ed.)}, 2005(1) (2005), 34--39.
\bibitem{maetal2005_2} Dengju Ma, Han Ren and Junjie Lu,  \newblock The crossing number of the circular graph $C(2m+2,m)$,  \newblock {\em Disc. Math.}, 304 (2005), 88--93.
\bibitem{madej1991} Tom Madej,  \newblock Bounds for the crossing number of the $N$-cube,  \newblock {\em J. Graph Th.}, 15(1) (1991), 81--97.
\bibitem{manueletal2013} Paul D. Manuel, Bharati Rajan, Indra Rajasingh and P. Vasanthi Beula,  \newblock Improved bounds on the crossing number of butterfly network,  \newblock {\em Disc. Math. Th. Comp. Sci.}, 15(2) (2013), 87--94.
\bibitem{mcquillanetal2015} Dan McQuillan, Shengjun Pan and R. Bruce Richter,  \newblock On the crossing number of $K_{13}$,  \newblock {\em J. Combin. Th. Series B}, 115 (2015), 224--235.
\bibitem{mcquillanrichter1992} Dan McQuillan and R. Bruce Richter,  \newblock On the crossing numbers of certain generalized Petersen graphs,  \newblock {\em Disc. Math.}, 104(3) (1992), 311--320.
\bibitem{meihuang2007} Hanfei Mei and Yuanqiu Huang,  \newblock The Crossing Number of $K_{1,5,n}$,  \newblock {\em Int. J. Math. Combin.}, 1(1) (2007), 33--44. $\hfill \Asterisk$
\bibitem{montaron2005} Bernard Montaron,  \newblock An Improvement of the Crossing Number Bound,  \newblock {\em J. Graph Th.}, 50(1) (2005), 43--54.
\bibitem{nahas2003} Nagi H. Nahas,  \newblock On the Crossing Number of $K_{m,n}$,  \newblock {\em Elec. J. Combin.}, 10, Note 8, (2003).
\bibitem{newcombe2019} Alex Newcombe,  \newblock {\em An efficient heuristic for crossing minimisation and its applications},  \newblock PhD Thesis, Flinders University, (2019).
\bibitem{norin2013} Sergey Norin,  \newblock Presentation at the BIRS Workshop on geometric and topological graph theory (13w5091),  \newblock October 1, (2013). https://www.birs.ca/events/2013/5-day-workshops/13w5091/videos/watch/201310011538-Norin.html.  Accessed:  07-06-2020.
\bibitem{ohringetal1995} Sabine R. Ohring, Maximilian Ibel, Sajal K. Das and Mohan J. Kumar,  \newblock On generalized fat trees.  In: {\em Parallel Processing Symposium, 1995. Proceedings, 9th international}, IEEE, (1995), pp. 37--44.
\bibitem{oporowskizhao2009} Bogdan Oporowski and David Zhao,  \newblock Coloring graphs with crossings,  \newblock {\em Disc. Math.}, 309(6) (2009), 2948--2951.
\bibitem{ouyangetal2014} Zhangdong Ouyang, Jing Wang and Yuanqiu Huang,  \newblock On the crossing number of $K_{1,1,m} \Box P_n$,  \newblock {\em Sci. Sinica Math.}, 44(12) (2014), 1337--1342.
\bibitem{ouyangetal2014_2} Zhangdong Ouyang, Jing Wang and Yuanqiu Huang,  \newblock The crossing number of the Cartesian product of paths with complete graphs,  \newblock {\em Disc. Math.}, 328 (2014), 71--78.
\bibitem{ouyangetal2018} Zhanhgdong Ouyang, Jing Wang and Yuanqiu Huang,  \newblock The crossing number of join of the generalized Petersen graph $P(3,1)$ with path and cycle,  \newblock {\em Disc. Math. Graph Th.}, 38(2) (2018), 351--370.
\bibitem{ouyangetal2018_2} Zhangdong Ouyang, Jing Wang and Yuanqiu Huang,  \newblock The Strong Product of Graphs and Crossing Numbers,  \newblock {\em Ars Combin.}, 137 (2018), 141--147.
\bibitem{pachetal2006} J\'{a}nos Pach, Rados Radoi\v{c}i\'{c}, G\'{a}bor Tardos and G\'{e}za T\'{o}th,  \newblock Improving the crossing lemma by finding more crossings in sparse graphs,  \newblock {\em Disc. Comp. Geom.}, 36(4) (2006), 527--552.
\bibitem{pachtoth1997} J\'{a}nos Pach and G\'{e}za T\'{o}th,  \newblock Graphs drawn with few crossings per edge,  \newblock {\em Combinatorica}, 17(3) (1997), 427--439.
\bibitem{panrichter2007} Shengjun Pan and R. Bruce Richter,  \newblock The crossing number of $K_{11}$ is 100,  \newblock {\em J. Graph Th.}, 56(2) (2007), 128--134.
\bibitem{petrillova2015} Jana Petrillov\'{a},  \newblock On the optimal drawings of Cartesian products of special 6-vertex graphs with path,  \newblock {\em Math. Model. Geom.}, 3(3) (2015), 19--28.
\bibitem{petrillova2019} Jana Petrillov\'{a},  \newblock On the optimal drawings of the products of special graphs.  \newblock In: {\em Proceedings of the 22nd Czech-Japan Seminar on Data Analysis and Decision Making} (2019), 135--144.
\bibitem{pinontoanrichter2004} Benny Pinontoan and R. Bruce Richter,  \newblock Crossing numbers and sequences of graphs I: general tiles,  \newblock {\em Austral. J. Combin.}, 30 (2004), 197--206.
\bibitem{qian2007} Chunhua Qian,  \newblock {\em The crossing number of a 5-vertex graph and the star graph},  \newblock Masters Thesis, Hunan Normal University, (2007). $\hfill \Asterisk$
\bibitem{qianhuang2007} Chunhua Qian and Yuanqiu Huang,  \newblock The Crossing Number of the Complete Quadruple Graph,  \newblock {\em J. Hunan Uni. Art. Sci. Nat. Sci. Ed.}, 19(4) (2007), 13--17. $\null\hfill \Asterisk$
\bibitem{rajanetal2011} Bharati Rajan, Indra Rajasingh and P. Vasanthi Beulah,  \newblock On the crossing number of generalized fat trees,  \newblock In: {\em International Conference on Informatics Engineering and Information Science}, Springer, Berlin, Heidelberg, (2011), pp. 440--448.
\bibitem{rajanetal2012} Bhatari Rajan, Indra Rajasingh and P. Vasanthi Beulah,  \newblock Crossing Number of Join of Triangular Snake with $mK_1$, Path and Cycle,  \newblock {\em Int. J. Comp. Appl.}, 44(17) (2012), 20--22.
\bibitem{richtersalazar2001} R. Bruce Richter and Gelasio Salazar,  \newblock The crossing number of $C_6 \times C_n$,  \newblock {\em Austral. J. Combin.}, 23 (2001), 135--143.
\bibitem{richtersalazar2002} R. Bruce Richter and Gelasio Salazar,  \newblock The crossing number of P(N,3),  \newblock {\em Graph. Combin.}, 18(2) (2002), 381--394.
\bibitem{richterthomassen1995} R. Bruce Richter and Carsten Thomassen,  \newblock Intersections of curve systems and the crossing number of $C_5 \times C_5$,  \newblock {\em Disc. Comput. Geom.}, 13(1) (1995), 149--159.
\bibitem{ringeisenbeineke1978} Richard D. Ringeisen and Lowell W. Beineke,  \newblock The crossing number of $C_3 \times C_n$,  \newblock {\em J. Combin. Th. Series B}, 24(2) (1978), 134--136.
\bibitem{saaty1964} Thomas L. Saaty,  \newblock The minimum number of intersections in complete graphs,  \newblock {\em Proc. Nat. Acad. Sci. USA}, 52 (1964), 688--690.
\bibitem{salazar1997} Gelasio Salazar,  \newblock {\em Crossing Numbers of Certain Families of Graphs},  \newblock PhD Thesis, Carleton University, (1997).
\bibitem{salazar2000} Gelasio Salazar,  \newblock A lower bound for the crossing number of $C_m \times C_n$,  \newblock {\em J. Graph Th.}, 35(3) (2000), 222--226.
\bibitem{salazar2005} Gelasio Salazar,  \newblock On the crossing numbers of loop networks and generalized Petersen graphs,  \newblock {\em Disc. Math.}, 302(1--3) (2005), 243--253.
\bibitem{salazarugalde2004} Gelasio Salazar and Edgardo Ugalde,  \newblock An Improved Bound for the Crossing Number of $C_m \times C_n$: a Self-Contained Proof Using Mostly Combinatorial Arguments,  \newblock {\em Graphs Combin.}, 20 (2004), 247--253.
\bibitem{sarazin1997} Marko L. Sara\v{z}in,  \newblock The crossing number of the generalized Petersen graph $P(10,4)$ is four,  \newblock {\em Math. Slovaca}, 47(2) (1997), 189--192.
\bibitem{schaefersurvey} Marcus Schaefer,  \newblock The graph crossing number and its variants: A survey,  \newblock {\em Elec. J. Combin.}, DS21 version 4, (2020).
\bibitem{schaefer2018} Marcus Schaefer,  \newblock {\em Crossing numbers of graphs.}  \newblock CRC Press, (2018).
\bibitem{shahrokhietal1995} Farhad Shahrokhi, Ondrej S\'{y}kora, L\'{a}szl\'{o} A. Sz\'{e}kely and Imrich Vrt'o,  \newblock Crossing Numbers of Meshes,  \newblock In: {\em International Symposium on Graph Drawing}, Springer, Berlin, Heidelberg, (1995), pp. 463--471.
\bibitem{shahrokhietal1998} Farhad Shahrokhi, Ondrej S\'{y}kora, L\'{a}szl\'{o} A. Sz\'{e}kely and Imrich Vrt'o,  \newblock Intersection of Curves and Crossing Number of $C_m \times C_n$ on Surfaces,  \newblock {\em Disc. Comp. Geom.}, 19 (1998), 237--247.
\bibitem{shanthinibabujee2016} Natarajan Shanthini and Jayapal Baskar Babujee,  \newblock Crossing Number for $S_4 + S_N$,  \newblock {\em Asian J. Res. Soc. Sci. Hum.}, 6(9) (2016), 2030--2041. $\hfill \Asterisk$
\bibitem{shanthinibabujee2016_2} Natarajan Shanthini and Jayapal Baskar Babujee,  \newblock Crossing numbers for join of the tripartite graph with $nK_1$ and path,  \newblock {\em Int. J. Pure Appl. Math.}, 109(5) (2016), 13--22. $\hfill \Asterisk$
\bibitem{stas2018} Michal Sta\v{s},  \newblock Cyclic permutations: crossing numbers of the join products of graphs,  \newblock In: {\em Proc. APLIMAT 2018: 17th Conf. Appl. Math. (Slovak University of Technology, Bratislava, 2018)}, (2018), pp. 979--987.
\bibitem{stas2018_2} Michal Sta\v{s},  \newblock Determining crossing numbers of graphs of order six using cycle permutations,  \newblock {\em Bull. Austral. Math. Soc}, 98 (2018), 353--362.
\bibitem{stas2019} Michal Sta\v{s},  \newblock Determining Crossing Number of Join of the Discrete Graph with Two Symmetric Graphs of Order Five,  \newblock {\em Symmetry}, 11(2) (2019), 123.
\bibitem{stas2019_2} Michal Sta\v{s},  \newblock Determining crossing number of one graph of order five using cyclic permutations,  \newblock In: {\em APLIMAT 2019: Proceedings of the 18th Conference on Applied Mathematics (Slovak University of Technology, Bratislava, 2019)}, (2019), pp. 1126--1133.
\bibitem{stas2019_3} Michal Sta\v{s},  \newblock Alternative proof on the crossing number of $K_{2,3,n}$,  \newblock {\em Math. Model. Geom.}, 7(1) (2019), 13--20.
\bibitem{stas2019_4} Michal Sta\v{s},  \newblock Alternative proof on the crossing number of $K_{1,4,n}$,  \newblock In: {\em Proceedings of the 22nd Czech-Japan Seminar on Data Analysis and Decision Making} (2019), 165--174.
\bibitem{stas2020} Michal Sta\v{s},  \newblock On the Crossing Numbers of the Joining of a Specific Graph on Six Vertices with the Discrete Graph,  \newblock {\em Symmetry}, 12(1) (2020), 1--12.
\bibitem{stas2020_2} Michal Sta\v{s},  \newblock  Join Products $K_{2,3} + C_n$,  \newblock {\em Mathematics}, 8(6) (2020), 1--9.
\bibitem{stas2020_3} Michal Sta\v{s},  \newblock On the crossing numbers of join products of five graphs of order six with the discrete graph,  \newblock {\em Opuscula}, 40(3) (2020), 383--397.
\bibitem{stas2020_4} Michal Sta\v{s},  \newblock On the crossing number of join product of the discrete graph with special graphs of order five,  \newblock {\em Elec. J. Graph Th. Appl.}, 8(2) (2020), 339--351.
\bibitem{stas2020_5} Michal Sta\v{s},  \newblock Determining Crossing Numbers of the Join Products of two Specific Graphs of Order Six With the Discrete Graph,  \newblock {\em Filomat}, 34:9 (2020), 2829--2846.
\bibitem{staspetrillova2018} Michal Sta\v{s} and Jana Petrillov\'{a},  \newblock On the join products of two special graphs on five vertices with the path and the cycle,  \newblock {\em Math. Model. Geom.}, 6(2) (2018), 1--11.
\bibitem{stasvaliska2021} Michal Sta\v{s} and Juraj Valiska,  \newblock On the crossing numbers of join products of $W_4 + P_n$ and $W_4 + C_n$,  \newblock {\em Opuscula Math.}, 41(1) (2021), 95--112.
\bibitem{su2011} Zhenhua Su,  \newblock The Crossing Numbers of a 6-vertex Graph $G$ Product with $S_n$,  \newblock {\em J. Math. Stud.}, 44(4) (2011), 411--416. $\hfill \Asterisk$
\bibitem{su2017} Zhenhua Su,  \newblock The Crossing Number of Cartesian Products of the Special Graph on Six Vertices with Stars,  \newblock {\em Math. Prac. Th.}, 47(12) (2017), 182--188. $\hfill \Asterisk$
\bibitem{su2017_2} Zhenhua Su,  \newblock The crossing number of the join product of $C_6 + 3K_2$ with $P_n$ and $C_n$,  \newblock {\em Oper. Res. Trans.}, 21(3) (2017), 23--34.
\bibitem{su2017_3} Zhenhua Su,  \newblock The Crossing Number of the Join Product of $K_{1,1,1,2} + P_n$  \newblock {\em Acta Math. Appl. Sinica}, 40(3) (2017), 345--354.
\bibitem{su2017_4} Zhenhua Su,  \newblock The Crossing Numbers of the Join of a 6-vertex Graphs $H$ with $P_n$ and $C_n$,  \newblock {\em J. Shanxi Norm. Uni.}, 31(3) (2017), 14--19. $\hfill \Asterisk$
\bibitem{su2018} Zhenhua Su,  \newblock Crossing numbers of join products of $K_{1,1,3} \vee C_n$,  \newblock {\em Comp. Eng. Appl.}, 54(9) (2018), 57--61.
\bibitem{su2018_2} Zhenhua Su,  \newblock The crossing numbers of $K_5 + P_n$,  \newblock {\em J. Math. Res. App.}, 38(4) (2018), 331--340.
\bibitem{suhuang2008} Zhenhua Su and Yuanqiu Huang,  \newblock Crossing Numbers of Some Classes of Cartesian Product Graphic,  \newblock {\em J. Jishou Uni. Nat. Sci. Ed.}, 29(6) (2008), 25--28. $\hfill \Asterisk$
\bibitem{suhuang2012} Zhenhua Su and Yuanqiu Huang,  \newblock Crossing Number of $K_{2,3} \vee P_n$,  \newblock {\em Appl. Math. China Ser. A}, 27(4) (2012), 488--492.
\bibitem{suhuang2012_2} Zhenhua Su and Yuanqiu Huang,  \newblock Crossing Numbers of $\{P_6^2 + e\} \times S_n$,  \newblock {\em J. Jishou Uni. Nat. Sci. Ed.}, 33(4) (2012), 20--24. $\hfill \Asterisk$
\bibitem{suhuang2013} Zhenhua Su and Yuanqiu Huang,  \newblock The Crossing Number of $S_5 \vee C_n$,  \newblock {\em J. Math. Stud.}, 46(4) (2013), 413--417. $\hfill \Asterisk$
\bibitem{suhuang2014} Zhenhua Su and Yuanqiu Huang,  \newblock Crossing Number of join of three 5-vertex graphs with $P_n$,  \newblock {\em App. Math. China}, 29(2) (2014), 245--252.
\bibitem{sykoravrto1993} Ondrej S\'{y}kora and Imrich Vrt'o,  \newblock On crossing numbers of hypercubes and cube connected cycles,  \newblock {\em BIT Num. Math.}, 33(2) (1993), 232--237.
\bibitem{tangetal2007} Ling Tang, Shengxiang L\"{u} and Yuanqiu Huang,  \newblock The Crossing Number of Cartesian Products of Complete Bipartite Graphs $K_{2,m}$ with Paths $P_n$,  \newblock {\em Graphs Combin.}, 23 (2007), 659--666.
\bibitem{vijayaetal2016} N. Vijaya, Bharati Rajan and Ibrahim Venkat,  \newblock Crossing numbers of join of a graph on six vertices with a path and a cycle,  \newblock {\em Int. J. Advance. Soft. Comp. Appl.}, 8(2) (2016), 41--51. $\hfill \Asterisk$
\bibitem{vrto2014} Imrich Vrt'o,  \newblock Bibliography on crossing numbers,  \newblock ftp://ftp.ifi.savba.sk/pub/ imrich/crobib.pdf, (2014).
\bibitem{wangetal2013} Guoqing Wang, Haoli Wang, Yuansheng Yang, Xuezhi Yang and Wenping Zheng,  \newblock An upper bound for the crossing number of augmented cubes,  \newblock {\em Int. J. Comp. Math.}, 90(2) (2013), 183--227.
\bibitem{wangetal2017} Haoli Wang, Xirong Xu, Yuansheng Yang, Bao Liu, Wenping Zheng and Guoqing Wang,  \newblock An upper bound for the crossing number of locally twisted cubes,  \newblock {\em Ars. Comb.}, 131 (2017), 87--106.
\bibitem{wangetal2015} Haoli Wang, Yuansheng Yang, Yan Zhou, Wenping Zheng and Guoqing Wang,  \newblock The crossing number of folded hypercubes,  \newblock {\em Util. Math.}, 98 (2015), 393--408.
\bibitem{wanghuang2005} Jing Wang and Yuanqiu Huang,  \newblock The Crossing Numbers of Cartesian Products of Paths with 6-Vertex Graphs,  \newblock {\em J. Jishou Uni. Nat. Sci. Ed.}, 26(2) (2005), 9--13. $\null\hfill \Asterisk$
\bibitem{wanghuang2008} Jing Wang and Yuanqiu Huang,  \newblock Crossing number of the complete tripartite graph $K_{1,10,n}$,  \newblock {\em Appl. Math. China Ser. A.}, 23(3) (2008), 349--356.
\bibitem{wanghuang2008_2} Jing Wang and Yuanqiu Huang,  \newblock The crossing number of $K_{2,4} \times P_n$,  \newblock {\em Acta Math. Sci., Ser. A., Chin. Ed.}, 28 (2008), 251--255. $\hfill \Asterisk$
\bibitem{wanghuang2008_3} Jing Wang and Yuanqiu Huang,  \newblock The Crossing Number of the Circulant Graph $C(3k-1;\{1,k\})$,  \newblock {\em Int. J. Math. Combin.}, 3 (2008), 79--84. $\hfill \Asterisk$
\bibitem{wangetal2019} Jing Wang, Zhangdon Ouyang and Yuanqiu Huang,  \newblock The Crossing Number of the Hexagonal Graph $H_{3,n}$,  \newblock {\em Disc. Math. Graph Th.}, 39 (2019), 547--554.
\bibitem{wangetal2013_2} Jing Wang, Lixi Zhang and Yuanqiu Huang,  \newblock On the Crossing Number of the Cartesian Product of a 6-Vertex Graph with $S_n$,  \newblock {\em Ars Combin.}, 109 (2013), 257--266.
\bibitem{wangma2017} Junshuai Wang and Dengju Ma,  \newblock The Crossing Number of $V_8 \times P_n$,  \newblock {\em J. Nantong Uni. Nat. Sci}, 16(3) (2017), 79--85. $\hfill \Asterisk$
\bibitem{wanghuangtoappear} Yuxi Wang and Yianqui Huang,  \newblock The Crossing Number of Cartesian Product of 5-Wheel with any Tree,  \newblock {\em Disc. Math. Graph Th.}, to appear.
\bibitem{watkins1969} Mark E. Watkins,  \newblock A Theorem on Tait Colorings with an Application to the Generalized Petersen Graphs,  \newblock {\em J. Combin. Th.}, 6 (1969), 152--164.
\bibitem{winterbach2005} Wynand Winterbach,  \newblock {\em The crossing number of a graph in the plane},  \newblock Master's Thesis, University of Stellenbosch, (2005).
\bibitem{woodall1993} Douglas R. Woodall,  \newblock Cyclic-order graphs and Zarankiewicz's crossing-number conjecture,  \newblock {\em J. Graph Th.}, 17(6) (1993), 657--671.
\bibitem{xiaohuang2003} Wenbing Xiao and Yuanqiu Huang,  \newblock The Crossing Numbers of Some Classes of Cartesian Product Graphic,  \newblock {\em J. Nat. Sci. Hunan Norm. Uni.}, 26(4) (2003), 3--7. $\null\hfill \Asterisk$
\bibitem{xiaoetal2004} Wenbing Xiao, Hongzhuan Wang and Yuanqiu Huang,  \newblock The Crossing Number of Star $S_n$ Cartesian Product and 6 Vevtex Graph,  \newblock {\em J. Hunan Uni. Arts Sci. Nat. Sci. Ed.}, 16(4) (2004), 15--17. $\hfill \Asterisk$
\bibitem{yangetal2005} Xiaofan Yang, David J. Evans and Graham M. Megson,  \newblock The locally twisted cubes,  \newblock {\em Int. J. Comput. Math.}, 82 (2005), 401--413.
\bibitem{yang2007} Xiwu Yang,  \newblock {\em The Crossing Number of Flower Snark and $K_m^- \Box P_n$},  \newblock Masters Thesis, Dalian University of Technology, (2007). $\hfill \Asterisk$
\bibitem{yangetal2004} Yuansheng Yang, Xiaohui Lin, Jianguo Lu and Xin Hao,  \newblock The crossing number of $C(n$; \{1, 3\}$)$,  \newblock {\em Disc. Math.}, 289(1--3) (2004), 107--118.
\bibitem{yangetal2017} Yuansheng Yang, Bo Lv, Baigong Zheng, Xirong Xu and Ke Zhang,  \newblock The crossing number of pancake graph $P_4$ is six,  \newblock {\em Ars Combin.}, 131 (2017), 43--53.
\bibitem{yangetal2017_2} Yuansheng Yang, Guoqing Wang, Haoli Wang and Yan Zhou,  \newblock The Erd\H{o}s and Guy's conjectured equality on the crossing number of hypercubes,  \newblock \arxiv{1201.4700}, submitted 2017.
\bibitem{yangzhao2001} Yuansheng Yang and Chengye Zhao,  \newblock The crossing number of $C_n(1,k)$,  \newblock {\em National symposium on software technology held by Chinese Computer Institute 2001}, (2001), 134--136. $\hfill \Asterisk$
\bibitem{yuan2011} Xiuhua Yuan,  \newblock On the crossing numbers of $K_{2,3} \vee C_n$,  \newblock {\em J. East China Norm. Uni. Nat. Sci.}, 2011(5) (2011), 21--24.
\bibitem{yuanhuang2007} Zihan Yuan and Yuanqiu Huang,  \newblock Crossing Number of $C_m + \{e_1\}$ and $C_m + \{e_1\} + \{e_2\}$ with $P_n$,  \newblock {\em J. Jishou Uni. Nat. Sci. Ed.}, 28(3) (2007), 16--18. $\hfill \Asterisk$
\bibitem{yuanhuang2007_2} Zihan Yuan and Yuanqiu Huang,  \newblock The Crossing Number of Cartesian Product of a 3-connected Graph with Six Vertices and a Path,  \newblock {\em Math. Th. Appl.}, 27(2) (2007), 49--51. $\hfill \Asterisk$
\bibitem{yuanhuang2010} Zihan Yuan and Yuanqiu Huang,  \newblock Crossing numbers of $K_{(1,1,2,2)} \times P_n$,  \newblock {\em Appl. Math. China Ser. A}, 25(1) (2010), 75--84.
\bibitem{yuanhuang2011} Zihan Yuan and Yuanqiu Huang,  \newblock The Crossing Number of $C(12,2) \times P_n$,  \newblock {\em Adv. Math. China}, 40(5) (2011), 587--594.
\bibitem{yuanhuang2011_2} Zihan Yuan and Yuanqiu Huang,  \newblock The Crossing Number of Petersen Graph $P(4,1)$ with Paths $P_n$,  \newblock {\em Oper. Res. Trans.}, 15(3) (2011), 95--106. $\hfill \Asterisk$
\bibitem{yuanhuangetal2008} Zihan Yuan and Yuanqiu Huang and Jinwang Liu,  \newblock The crossing number of $C(7,2) \times P_n$,  \newblock {\em Adv. Math. China}, 37(2) (2008), 245--253. $\hfill \Asterisk$
\bibitem{yuanetal2013} Zihan Yuan, Yuanqiu Huang and Jinwang Liu,  \newblock The Crossing Number of Cartesian Product of Circulant Graph (9,2) with Path $P_n$,  \newblock {\em Acta Math. Appl. Sinica}, 36(2) (2013), 350--362.
\bibitem{yuanetal2008} Zihan Yuan, Tang Lien, Yuanqiu Huang and Jinwang Liu,  \newblock The Crossing Number of $C(8,2) \Box P_n$,  \newblock {\em Graphs Combin.}, 24(6) (2008), 597--604.
\bibitem{yuanetal2009} Zihan Yuan, Jing Wang and Yuanqiu Huang,  \newblock The Crossing Number of Cartesian Product of $P_n$ with Circulant Graph $C(10,2)$,  \newblock {\em Acta Math. Appl. Sinica}, 32(6) (2009), 1133--1144.
\bibitem{yueetal2014} Weijun Yue, Yuanqiu Huang and Zhangdong Ouyang,  \newblock On crossing numbers of join of $W_4 + C_n$. {\em Comp. Eng. Appl.}, 50(18) (2014), 79--84.
\bibitem{yueetal2013} Weijun Yue, Yuanqiu Huang and Ling Tang,  \newblock The crossing number of join products of three 6-vertex graphs with $C_n$,  \newblock {\em J. Hunan Uni. Arts Sci. Nat. Sci. Ed.}, 25(4) (2013), 1--7. $\hfill \Asterisk$
\bibitem{zarankiewicz1955} Kazimierz Zarankiewicz,  \newblock On a problem of P. Turan concerning graphs,  \newblock {\em Fund. Math.}, 41(1) (1955), 137--145.
\bibitem{zhanghuang2011} Lixi Zhang and Yuanqiu Huang,  \newblock On the crossing number of Cartesian products of a 5-vertex graph $G_7$ with star,  \newblock {\em J. Hunan Uni. Arts Sci. Nat. Sci. Ed.}, 23(4) (2011), 3--6. $\hfill \Asterisk$
\bibitem{zhangetal2008} Lixi Zhang, Bo Li and Yuanqiu Huang,  \newblock Crossing Number of the Cartesian Product of a 6-Vertex Graph with $S_n$,  \newblock {\em J. Jishou Uni. Nat. Sci. Ed.}, 29(5) (2008), 23--29. $\hfill \Asterisk$
\bibitem{zhangetal2008_2} Lixi Zhang, Bo Li and Yuanqiu Huang,  \newblock On the crossing number of products of a graph of order 6 and the star,  \newblock {\em J. Hunan Uni. Arts Sci. Nat. Sci. Ed.}, 20(1) (2008), 16--19. $\hfill \Asterisk$
\bibitem{zhang2013} Yulan Zhang,  \newblock The Crossing Numbers of $S_4 \vee C_n$,  \newblock {\em J. Gansu Lianhe Uni. Nat. Sci. Ed}, 27(2) (2013), 31--34. $\hfill \Asterisk$
\bibitem{zhenghuang2011} Dunyoung Zheng and Yuanqiu Huang,  \newblock Crossing Number of the Join Graph of a 5-Vertex Graph and Path,  \newblock {\em J. Shantou Uni. Nat. Sci. Ed.}, 26(4) (2011), 11--17. $\null\hfill \Asterisk$
\bibitem{zhengetal2007} Wenping Zheng, Xiaohui Lin, Yuansheng Yang and Chong Cui,  \newblock On the crossing number of $K_m \Box P_n$,  \newblock {\em Graphs Combin.}, 23(3) (2007), 327--336.
\bibitem{zhengetal2008} Wenping Zheng, Xiaohui Lin, Yuansheng Yang and Chengrui Deng,  \newblock On the crossing number of $K_m \Box C_n$ and $K_{m,l} \Box P_n$,  \newblock {\em Disc. Appl. Math.}, 156 (2008), 1892--1907.
\bibitem{zhengetal2008_2} Wenping Zheng, Xiaohui Lin, Yuansheng Yang and Chengrui Deng,  \newblock The Crossing Number of Kn\"{o}del Graph $W_{3,n}$,  \newblock {\em Util. Math.}, 75 (2008), 211--224.
\bibitem{zhengetal2009} Wenping Zheng, Xiaohui Lin, Yuansheng Yang and Yang Gui,  \newblock On the Crossing Numbers of the $k$-th Power of $P_n$,  \newblock {\em Ars Combin.}, 92 (2009), 397--409.
\bibitem{zhengetal2008_3} Wenping Zheng, Xiaohui Lin, Yuansheng Yang and Xiwu Yang,  \newblock The crossing number of flower snarks and related graphs,  \newblock {\em Ars Combin.}, 86 (2008), 57--64.
\bibitem{zhengetal2011} Wenping Zheng, Xiaohui Lin, Yuansheng Yang and Xiwu Yang,  \newblock The Crossing Numbers of Cartesian Product of Cone Graph $C_m + \overline{K_l}$ with Path $P_n$,  \newblock {\em Ars Combin.}, 98 (2011), 433--445.
\bibitem{zhouhuang2007} Zhidong Zhou and Yuanqiu Huang,  \newblock The Crossing Number of $K_{3,3} \times P_n$,  \newblock {\em J. Hunan Norm. Uni. Nat. Sci.}, 30(1) (2007), 31--34. $\hfill \Asterisk$
\bibitem{zhouetal2013} Zhidong Zhou, Yuanqiu Huang, Xiaoduo Peng and Juan Ouyang,  \newblock The Crossing Number of the Joint Graph of a Small Graph and a Path or a Cycle,  \newblock {\em J. Sys. Sci. Math. Sci.}, 33(2) (2013), 206--216.
\bibitem{zhouli2016} Zhidong Zhou and Long Li,  \newblock On the crossing number of the joint graph,  \newblock {\em J. Shaoyang Uni. Nat. Sci. Ed.}, 13(3) (2016), 16--24. $\hfill \Asterisk$
\bibitem{zhouli2016_2} Zhidong Zhou and Long Li,  \newblock On the crossing numbers of the joins of the special graph on six vertices with $nK_1$, $P_n$ or $C_n$,  \newblock {\em Oper. Res. Trans.}, 20(4) (2016), 115--126.
\bibitem{zhouli2017} Zhidong Zhou and Long Li,  \newblock On the Crossing Number of Join Product of Path and a Graph of Order Six,  \newblock {\em Math. Appl.}, 30(1) (2017), 72--77.



\end{thebibliography}
\end{document}